\newif\iffinal
\finalfalse	
\finaltrue	

\documentclass[letterpaper,12pt,reqno]{amsart}
\RequirePackage[utf8]{inputenc}
\usepackage[portrait,margin=2.5cm]{geometry}
\iffinal\else\usepackage[notref,notcite]{showkeys}\fi
\usepackage{mathrsfs,soul}
\usepackage{hyperref}
\usepackage[foot]{amsaddr}
\usepackage{amssymb,amsthm,amsfonts,amsbsy,latexsym,dsfont}
\usepackage[textsize=small]{todonotes}       
\usepackage{graphicx}
\usepackage[numeric,initials,nobysame]{amsrefs}
\usepackage{upref,setspace}

\usepackage{tikz}
\usetikzlibrary{calc,intersections,through,backgrounds,shapes.geometric}
\usetikzlibrary{graphs}

\usepackage{pdfsync}

\usepackage{caption}
\usepackage{subcaption}
\usepackage{fourier}

\usepackage{enumerate}
\newenvironment{enumeratei}{\begin{enumerate}[\upshape (i)]}{\end{enumerate}}
\newenvironment{enumeratea}{\begin{enumerate}[\upshape (a)]}{\end{enumerate}}

\numberwithin{equation}{section}
\numberwithin{figure}{section}
\numberwithin{table}{section}

\newtheorem{thm}{Theorem}[section]
\newtheorem{lem}[thm]{Lemma}
\newtheorem{cor}[thm]{Corollary}
\newtheorem{prop}[thm]{Proposition}

\newtheorem{ass}[thm]{Assumption}
\newtheorem*{ass*}{Assumption}
\newtheorem*{theorem*}{Theorem}

\newtheorem{lemma}[thm]{Lemma}

\newtheorem*{thm3.9*}{Theorem 3.9*}

\theoremstyle{definition}
\newtheorem{rem}{Remark}

\renewcommand{\leq}{\leqslant}
\renewcommand{\geq}{\geqslant}

\newcommand{\ind}{\mathds{1}}
\newcommand{\eps}{\varepsilon}

\newcommand{\set}[1]{\left\{#1\right\}}

\newcommand{\probc}{\stackrel{\mathrm{P}}{\longrightarrow}}

\newcommand{\weakc}{\stackrel{\mathrm{w}}{\longrightarrow}}
\newcommand{\convas}{\stackrel{\mathrm{a.s.}}{\longrightarrow}}

\def\qed{ \hfill $\blacksquare$}

\newcommand{\cA}{\mathcal{A}}\newcommand{\cB}{\mathcal{B}}\newcommand{\cC}{\mathcal{C}}
\newcommand{\cD}{\mathcal{D}}\newcommand{\cE}{\mathcal{E}}\newcommand{\cF}{\mathcal{F}}
\newcommand{\cG}{\mathcal{G}}\newcommand{\cH}{\mathcal{H}}

\newcommand{\cM}{\mathcal{M}}\newcommand{\cN}{\mathcal{N}}\newcommand{\cO}{\mathcal{O}}
\newcommand{\cP}{\mathcal{P}}
\newcommand{\cS}{\mathcal{S}}\newcommand{\cT}{\mathcal{T}}\newcommand{\cU}{\mathcal{U}}
\newcommand{\cV}{\mathcal{V}}

\newcommand{\vM}{\mathbf{M}}
\newcommand{\vP}{\mathbf{P}}

\newcommand{\vX}{\mathbf{X}}\newcommand{\vZ}{\mathbf{Z}}

\newcommand{\ve}{\mathbf{e}}
\newcommand{\vh}{\mathbf{h}}

\newcommand{\vp}{\mathbf{p}}\newcommand{\vq}{\mathbf{q}}
\newcommand{\vt}{\mathbf{t}}\newcommand{\vu}{\mathbf{u}}
\newcommand{\vw}{\mathbf{w}}\newcommand{\vx}{\mathbf{x}}

\newcommand{\mvgamma}{\boldsymbol{\gamma}}

\newcommand{\mvpi}{\boldsymbol{\pi}}

\newcommand{\mvxi}{\boldsymbol{\xi}}

\newcommand{\bG}{\mathbb{G}}

\newcommand{\bN}{\mathbb{N}}
\newcommand{\bR}{\mathbb{R}}
\newcommand{\bT}{\mathbb{T}}

\newcommand{\dA}{\mathds{A}}\newcommand{\dB}{\mathds{B}}

\newcommand{\dV}{\mathds{V}}

\newcommand{\sS}{\mathscr{S}}

\newcommand{\sP}{\mathfrak{P}}
\newcommand{\sF}{\mathfrak{F}}

\DeclareMathOperator{\E}{\mathbb{E}}
\DeclareMathOperator{\pr}{\mathbb{P}}

\DeclareMathOperator{\var}{Var}

\DeclareMathOperator{\dis}{dis}

\DeclareMathOperator{\GH}{GH}
\DeclareMathOperator{\GHP}{GHP}
\DeclareMathOperator{\diam}{diam}
\DeclareMathOperator{\height}{ht}
\DeclareMathOperator{\ord}{ord}
\DeclareMathOperator{\dist}{dist}
\DeclareMathOperator{\exec}{exc}

\DeclareMathOperator{\scl}{scl}
\DeclareMathOperator{\er}{er}
\DeclareMathOperator{\mass}{mass}
\DeclareMathOperator{\con}{con}
\DeclareMathOperator{\nr}{nr}
\DeclareMathOperator{\dfs}{dfs}
\DeclareMathOperator{\per}{per}
\DeclareMathOperator{\jt}{jt}
\DeclareMathOperator{\partition}{ptn}

\newcommand{\sss}{\scriptscriptstyle}

\newcommand{\erdos}{Erd\H{o}s-R\'enyi }
\newcommand{\ldown}{l^2_{\downarrow}}

\begin{document}

\title[Continuum limit of inhomogeneous random graphs]{Continuum limit of critical inhomogeneous random graphs}

\date{}

\author[Bhamidi]{Shankar Bhamidi$^1$}
\address{$^1$Department of Statistics and Operations Research, 304 Hanes Hall, University of North Carolina, Chapel Hill, NC 27599}
\author[Sen]{Sanchayan Sen$^2$}
\address{$^2$Department of Mathematics and Computer Science, Eindhoven University of Technology, the Netherlands}
\author[Wang]{Xuan Wang$^3$}
\address{$^3$Databricks, 160 Spear Street, San Francisco, CA 94105}
\email{bhamidi@email.unc.edu, sanchayan.sen1@gmail.com, xuanwang9527@gmail.com}

\maketitle
\begin{abstract}
The last few years have witnessed tremendous interest in understanding the structure as well as the behavior of dynamics for inhomogeneous random graph models to gain insight into real-world systems. In this study we analyze the maximal components at criticality of one famous class of such models, the rank-one inhomogeneous random graph model \cite{norros2006conditionally}, \cite[Section 16.4]{BBSJOR07}. Viewing these components as measured random metric spaces, under finite moment assumptions for the weight distribution, we show that the components in the critical scaling window with distances scaled by $n^{-1/3}$ converge in the Gromov-Haussdorf-Prokhorov metric to rescaled versions of the limit objects identified for the \erdos random graph components at criticality in \cite{addario2012continuum}. A key step is the construction of connected components of the random graph through an appropriate tilt of a fundamental class of random trees called $\vp$-trees \cite{camarri2000limit,aldous2004exploration}.  This is the first step in rigorously understanding the scaling limits of objects such as the minimal spanning tree and other strong disorder models from statistical physics \cite{braunstein2003optimal} for such graph models. By asymptotic equivalence \cite{janson2010asymptotic}, the same results are true for the Chung-Lu model \cite{chung2002average,chung2002connected,chung2006complex} and the Britton-Deijfen-Martin-L\"{o}f model \cite{britton2006generating}. A crucial ingredient of the proof of independent interest are tail bounds for the height of $\vp$-trees. The techniques developed in this paper form the main technical bedrock for the general program developed in \cite{SBSSXW-universal} for proving universality of the continuum scaling limits in the critical regime for a wide array of other random graph models including the configuration model and inhomogeneous random graphs with general kernels \cite{BBSJOR07}.

\end{abstract}

\section{Introduction}
\label{sec:intro}
Motivated by applications and empirical data, the last few years have seen a tremendous interest in formulating and studying a wide array of random graph models, estimating the parameters in the model from data, and studying dynamic processes such as epidemics on such models to gain insight into real-world systems, see e.g. \cite{albert2002statistical,newman2003structure,van2009random,durrett2007random,BBSJOR07,chung2006complex} and the references therein. In such random graph models different vertices have different propensities for connecting to other vertices. To fix ideas consider the main model analyzed in this work:

\noindent{\bf Rank-one model:} This version of the model was introduced by Norros and Reittu \cite{norros2006conditionally,BBSJOR07} (with a variant arising in the work of Aldous in the construction of the standard multiplicative coalescent \cite{aldous1997brownian}), and is sometimes referred to as the Norros-Reittu model. We construct a random graph on the vertex set $[n] = \set{1,2,\ldots, n}$ as follows. Each vertex $i\in [n]$ has an associated weight $w_i\geq 0$. Think of this as the propensity of the vertex to form friendships (form edges) in a network.  Write $\vw = (w_i)_{i\in [n]}$ for the vector of weights and let $l_n = \sum_{i=1}^n w_i$ be the sum of these weights.
The weights actually form a triangular array $\vw = \vw^{\sss(n)} = (w_i^{\sss(n)} : i \in [n])$, but we will omit $n$ in the notation.
Taking the weight sequence $\vw$ as input, the random graph is constructed as follows. Define probabilities $p_{ij} := 1-\exp(-w_i w_j/l_n)$. Construct the random graph $\cG_n^{\nr}(\vw)$ by putting an edge $\set{i,j}$ between vertices $i,j$ with probability $p_{ij}$, independently across edges.

This is an important example of the general class of inhomogeneous random graphs analyzed by Bollob\'as, Janson and Riordan \cite{BBSJOR07}. In Section \ref{sec:mot} below, we describe in more detail why this model forms the key to understanding the critical regime in a wide array of other models including the configuration model and inhomogeneous random graph models modulated via a general kernel.  At this stage let it suffice to say this model is also closely related to two other famous models of inhomogeneous random graphs (and in fact shown to be asymptotically equivalent in a number of settings \cite{janson2010asymptotic}).

\begin{enumeratea}
	\item {\bf Chung-Lu model \cite{chung2002average,chung2002connected,chung2006complex}:} Given the set of weights $\vw$ as above, here one attaches edges independently with probability
	\[p_{ij}:= \max\set{\frac{w_i w_j}{l_n},1}.\]
	\item {\bf Britton-Deijfen-Martin-L\"{o}f model \cite{britton2006generating}: } Here one attaches edges independently with probability
	\[p_{ij}:= \frac{w_i w_j}{l_n+ w_i w_j}.\]
\end{enumeratea}

These models are inhomogeneous in the sense that different vertices have different proclivity to form edges. Further, assume that the empirical distribution of weights $F_n  = n^{-1} \sum_{i=1}^n \delta_{w_i}$ satisfies
\begin{equation}
\label{eqn:fn-to-f}
	F_n \weakc F , \qquad \mbox{as } n\to\infty,
\end{equation}
 for a limiting cumulative distribution function $F$, where $\weakc$ denotes weak convergence. Then by \cite[Theorem 3.13]{norros2006conditionally} as $n\to\infty$, the degree distribution converges in the sense that for $k\geq 0$, writing $N_k(n)$ for the number of vertices with degree $k$,
 \[\frac{N_k(n)}{n} \probc \E\left(e^{-W} \frac{W^k}{k!}\right),\qquad k\geq 0, \]
where $W\sim F$ and $\probc$ denotes convergence in probability. Thus, asymptotically one can get any desired tail behavior for the degree distribution by choosing the weight sequence appropriately.  Also note that the \erdos random graph $\cG^{\text{er}}(n, \lambda/n)$ is a special case of the Norros-Reittu model where all the weights $w_i \equiv -n\log(1-\lambda/n)\approx \lambda$.

Aside from applications, inhomogeneous random graph models have sparked a lot of interest in the statistical physics community, in particular in understanding how the network structure affects weak and strong disorder models of flow, e.g. first passage percolation, minimal spanning tree \cite{braunstein2003optimal} etc. In the next section we describe what is known about how the network transitions from the subcritical to the supercritical regime and then describe these conjectures from statistical physics in more detail.

\subsection{Connectivity and phase transition}
\label{sec:intro-conn-ph-trans}
The main aim of this study is the structural properties of the maximal components in the critical regime. In order to define the critical regime, we first recall known connectivity properties of the model. Let $W\sim F$ where $F$ as before denotes the limiting weight distribution as in \eqref{eqn:fn-to-f}. Assume that $0< \E(W^2)< \infty$ and further
\begin{equation}
\label{eqn:convg-sum-nu}
	\frac{\sum_{i=1}^n w_i^2}{\sum_{i=1}^n w_i} \to \frac{\E(W^2)}{\E(W)}, \qquad \mbox{as } n\to\infty.
\end{equation}
Define the parameter
\begin{equation}
\label{eqn:nu-def}
\nu = \frac{\E(W^2)}{\E(W)}.
\end{equation}
 Write $\cC_n^{\sss(i)}$ for the $i$-th largest component of $\cG_n^{\nr}(\vw)$ (breaking ties arbitrarily) and write $|\cC_n^{\sss(i)}|$ for the number of vertices in this component.
Then by \cite[Theorem 3.1 and Section 16.4]{BBSJOR07} (also see \cite{chung2002connected,norros2006conditionally,britton2006generating}), the phase transition is given in the following results:
\begin{enumeratei}
	\item {\bf Subcritical regime:} If $\nu <1$, then $|\cC_n^{\sss(1)}|/n \probc 0$ as $n \to \infty$.
	\item {\bf Supercritical regime:} If $\nu > 1$, then
	\[\frac{|\cC_n^{\sss(1)}|}{n} \probc {\rho(F)} > 0,\]
	where ${\rho(F)}$ is the survival probability of an associated approximating branching process.
	\item {\bf Critical regime:} Thus $\nu=1$ corresponds to the critical regime. This is the regime of interest for this paper. Similar to the \erdos random graph, we will study the entire critical scaling window, by working with weights $(1+\lambda/n^{1/3}) w_i$, {for fixed $\lambda \in \bR$,} so that the connection probability is
		\begin{equation}
		\label{eqn:pij-def}
		p_{ij}(\lambda):= 1-\exp\left(-\left(1+\frac{\lambda}{n^{1/3}}\right)\frac{w_i w_j}{l_n}\right).
		\end{equation}
We write $\cG_n^{\nr}(\vw, \lambda)$ for the corresponding random graph on the vertex set $[n]$. In this critical regime (with $\nu=1$), assuming 
$\sum_{i=1}^n w_i^3/n \to \E(W^3)< \infty$ as $n \to \infty$,
it is known \cite{SBVHJL10,Hofs09a} that for any fixed $i$, the $i$-th largest component scales like $|\cC_n^{\sss(i)}(\lambda)|\sim n^{2/3}$.  We shall give a precise description of this result in Section \ref{sec:descp}.
\end{enumeratei}

\subsection{Motivation and outline}
\label{sec:mot}
 Let us now informally describe the motivations behind our work. The main aim of this paper is to study the maximal components $\cC_n^{\sss(i)}(\lambda)$ for the rank-one model above in the critical scaling window and show that \emph{these maximal components viewed as metric spaces with edge lengths rescaled by $n^{-1/3}$ converge to random fractals related to the continuum random tree}. A natural question is why one should focus on this particular class of random graph models. In the following paragraphs, we give two major motivations of our work:

\begin{enumeratea}
	\item {\bf Universality for random graph processes at criticality:} The nature of emergence and scaling limits of {\bf component sizes} of maximal components in the critical regime of the Norros-Reittu model (and the closely related multiplicative coalescent) have recently been observed in a number of other random graph models including the configuration model \cite{nachmias2010critical,joseph2010component} as well as a general class of dynamic random graph processes called bounded-size rules (see e.g. \cite{spencer2007birth,bhamidi2012augmented}). The first step in understanding the metric structure of the maximal components for these models in the critical regime is the rank-one model. This paper forms the main technical bedrock for proving continuum scaling limits in the critical regime for a wide array of other random graph models in \cite{SBSSXW-universal} including the configuration model and inhomogeneous random graphs with general kernels \cite{BBSJOR07}. Using the technical tools developed in this paper, in particular Section \ref{sec:ptree-dfs-connected-irg}, we paraphrase the main theme of  \cite{SBSSXW-universal} as follows:
	\emph{\begin{quote}
		A general program for proving universality for the metric structure of the maximal components in the critical regime is proven using connections established in this paper between connected components in the rank-one model and tilted versions of $\vp$-trees (Proposition \ref{prop:distrib-gone-tildg} and Theorem \ref{thm:conv-condition-on-connected}). This general program coupled with model specific analysis show that maximal components in the critical regime of a number of fundamental random graph models, including the configuration model (under moment conditions),  and inhomogeneous random graph models with finite type space all satisfy analogous results to Theorem \ref{thm:inhom-random-graph} with distances in maximal components scaling like $n^{1/3}$.
	\end{quote}}
	  Thus it is not just the main results in this paper that are of interest, rather the {\bf proof section} of this paper sets out the tools required to prove this general program of universality.
	\item {\bf Scaling limits at criticality and the minimal spanning tree:} Our second main motivation was to rigorously understand predictions from statistical physics (see e.g. \cite{braunstein2003optimal} and the references therein) which predict that most inhomogeneous random graph models in the critical regime satisfy a remarkably universal behavior in the following sense.  Assume that the limiting degree distribution has finite third moments, then distances in the maximal components in the critical regime scale like $n^{1/3}$. Further consider the minimal spanning tree on the giant component in the supercritical regime where each edge is assigned a Uniform$[0,1]$ edge length. It is conjectured that (graph) distances in this object also scale like $n^{1/3}$. For the case of the \erdos random graph, this entire program has been carried forth in \cite{addario2013scaling}. Proving this conjecture for general inhomogenous random graphs rigorously turns out to be technically quite challenging and in particular requires a number of non-obvious assumptions, see e.g. Assumption \ref{ass:wt-seq}(d). To strengthen convergence in the stronger $l^4$ metric we needed to derive tail bounds for the height of $\vp$-trees (Theorem \ref{thm:ht-p-tree}) which are of independent interest.
	
\end{enumeratea}

\noindent{\bf Organization of the paper:}  In Section \ref{sec:prelim}, we start with the appropriate spaces and topology for convergence of a collection of metric spaces and define the Gromov-Haussdorf-Prokhorov metric. We state our main results in Section \ref{sec:results} and defer a precise description of the limit objects to Section \ref{sec:descp}. We discuss our main results and their relevance as well as the technical challenges in extending these results in Section \ref{sec:disc}. Starting from Section \ref{sec:proofs-comp} we prove the main results.

\section{Notation and preliminaries}
\label{sec:prelim}
We introduce some basic notation in Section \ref{sec:notation}. In Section \ref{sec:topology-prelim} we define relevant notions of convergence of measured metric spaces. In Section \ref{sec:graph-tree-prelim} we recall graph theoretic definitions required in the sequel. 

\subsection{Notation and conventions}\label{sec:notation}

For any set $A$, we write $|A|$ for its cardinality and $\ind\set{A}$ for the associated indicator function. Given two intervals $A, B \subset \bR$, we write $C(A,B)$ for the space of continuous functions $f: A \to B$, equipped with the $L_\infty$-norm $\|f\|_\infty := \sup_{x \in A}|f(x)|$. We write $D(A,B)$ for the space of RCLL (right-continuous-left-limit) functions $f: A \to B$, equipped with the Skorohod topology. We use the standard Landau notation of $o(\cdot)$, $O(\cdot)$ and the corresponding \emph{order in probability} notation $o_P(\cdot)$ and $O_P(\cdot)$.

We use $\pr(\cdot)$ to denote the canonical probability measure whose meaning should be clear from the context and $\E[ \cdot ]$ to denote the corresponding expectation operator. 
 The abbreviation ``\emph{i.i.d.}'' stands for ``independent and identically distributed''. As already mentioned in the introduction, we use $\probc$, $\weakc$ and $\convas$ to denote convergence in probability, weak convergence and almost-sure convergence. We say a sequence of events $E_n$, $n \in \bN$, occur \emph{with high probability} if $\pr(E_n) \to 0$ as $n \to \infty$.
We use notation such as $K_{\ref{thm:ht-p-tree}}$ and $K_{\ref{cor:tightness-l-uniform-integrable}}$ to denote absolute constants; here the objects in the subscript refer to the corresponding results where these constants are first introduced; e.g. $K_{\ref{thm:ht-p-tree}}$ is the constant in Theorem \ref{thm:ht-p-tree} and $K_{\ref{cor:tightness-l-uniform-integrable}}$ is the constant in Corollary \ref{cor:tightness-l-uniform-integrable}. Local constants are denoted by $C_1, C_2, \cdots$ or $B_1, B_2, \cdots$.


\subsection{Topology on the space of measured metric spaces}
\label{sec:topology-prelim}

We mainly follow \cite{EJP2116,addario2013scaling,burago2001course}. All metric spaces under consideration will be measured compact metric spaces. Let us recall the Gromov-Haussdorf distance $d_{\GH}$ between metric spaces.  Fix two metric spaces $X_1 = (X_1,d_1)$ and $X_2 = (X_2, d_2)$. For a subset $C\subseteq X_1 \times X_2$, the distortion of $C$ is defined as
\begin{equation}
	\label{eqn:def-distortion}
	\dis(C):= \sup \set{|d_1(x_1,y_1) - d_2(x_2, y_2)|: (x_1,x_2) , (y_1,y_2) \in C}.
\end{equation}
A correspondence $C$ between $X_1$ and $X_2$ is a measurable subset of $X_1 \times X_2$ such that for every $x_1 \in X_1$ there exists at least one $x_2 \in X_2$ such that $(x_1,x_2) \in C$ and vice-versa. The Gromov-Haussdorf distance between the two metric spaces  $(X_1,d_1)$ and $(X_2, d_2)$ is defined as
\begin{equation}
\label{eqn:dgh}
	d_{\GH}(X_1, X_2) = \frac{1}{2}\inf \set{\dis(C): C \mbox{ is a correspondence between } X_1 \mbox{ and } X_2}.
\end{equation}

We will need a metric that also keeps track of associated measures on the corresponding spaces. A compact measured metric space $(X, d , \mu)$ is a compact metric space $(X,d)$ with an associated finite measure $\mu$ on the Borel sigma algebra $\cB(X)$. Given two compact measured metric spaces $(X_1, d_1, \mu_1)$ and $(X_2,d_2, \mu_2)$ and a measure $\pi$ on the product space $X_1\times X_2$, the discrepancy of $\pi$ with respect to $\mu_1$ and $\mu_2$ is defined as
\begin{equation}
	\label{eqn:def-discrepancy}
	D(\pi;\mu_1, \mu_2):= ||\mu_1-\pi_1|| + ||\mu_2-\pi_2||
\end{equation}
where $\pi_1, \pi_2$ are the marginals of $\pi$ and $||\cdot||$ denotes the total variation of signed measures. Then the Gromov-Haussdorf-Prokhorov distance between $X_1$ and $X_2$ is defined as
\begin{equation}
\label{eqn:dghp}
	d_{\GHP}(X_1, X_2):= \inf\set{ \max\left(\frac{1}{2} \dis(C),~D(\pi;\mu_1,\mu_2),~\pi(C^c)\right) },
\end{equation}
where the infimum is taken over all correspondences $C$ and measures $\pi$ on $X_1 \times X_2$.

Denote by $\sS$ the collection of all measured metric spaces $(X,d,\mu)$. The function $d_{\GHP}$ is a pseudometric on $\sS$, and defines an equivalence relation $X \sim Y \Leftrightarrow d_{\GHP}(X,Y) = 0$ on $\sS$. Let $\bar \sS := \sS / \sim $ be the space of isometry equivalent classes of measured compact metric spaces and $\bar d_{\GHP}$ the induced metric. Then by \cite{EJP2116}, $(\bar \sS, \bar d_{\GHP})$ is a complete separable metric space. To ease notation, we will continue to use $(\sS, d_{\GHP})$ instead of $(\bar \sS, \bar d_{\GHP})$ and $X = (X, d, \mu)$ to denote both the metric space and the corresponding equivalence class.

Since we will be interested in not just one metric space but an
infinite collection of metric spaces, the relevant space of interest is a subset of $\sS^{\bN}$. For a fixed measured compact metric space $(X,d,\mu)$, define the diameter as $\diam(X)  := \sup_{x,y \in X} d(x,y)$ and the total mass as $\mass(X) := \mu(X)$. Then the relevant space for our study will be
\begin{equation}
\label{eqn:def-space-m}
	\mathcal{M}:=\set{(X_1, X_2,\hdots): X_i=(X_i, d_i, \mu_i)\in \sS,
	\ \sum_{i=1}^{\infty}(\diam(X_i)^4+\mass(X_i)^4)<\infty }.
\end{equation}
The two relevant topologies on this space are the following:
\begin{enumeratei}
	\item {\bf Product topology:} We use $\cT_1$ for the product topology inherited using $d_{\GHP}$ on a single co-ordinate.
	\item {\bf $l^4$ metric \cite{addario2012continuum}: } We use $\cT_2$ for the topology on $\cM$ induced by the distance
	\begin{equation}
	\label{eqn:l4-dist-def}
		\dist((X_1, X_2,\hdots), (X_1', X_2',\hdots)):=\left(\sum_{i=1}^{\infty}d_{\GHP}(X_i, X_i')^4\right)^{1/4}.
	\end{equation}
\end{enumeratei}

The aim of this paper is to study the limits of connected components of random graphs viewed as measured metric spaces.  In order to state our results, both the metric and the corresponding measure need to be re-scaled appropriately. To make this precise, we introduce the scaling operator $\scl(\alpha,\beta)$, for fixed constants $\alpha, \beta \in (0,\infty)$, as follows:
\begin{align*}
	\scl(\alpha, \beta) : \sS &\to \sS, \qquad \scl(\alpha, \beta)[(X , d , \mu)]:= (X, d',\mu'),
\end{align*}
where $d'(x,y) := \alpha d(x,y)$ for all $x,y \in X$, and $\mu'(A) := \beta \mu(A)$ for $A \subset X$. For simplicity, we write the output of the above scaling operator as $\scl(\alpha, \beta) X$. Using the definition of $d_{\GHP}$, it is easy to check that for $X \in \sS$ and for fixed $\alpha, \beta> 0$,
\begin{equation*}
	d_{\GHP}( \scl(\alpha, \beta) X, X ) \leq |\alpha-1| \cdot \diam(X) + |\beta - 1| \cdot \mass(X).
\end{equation*}
Note that $\diam(\cdot)$ and $\mass(\cdot)$ are both continuous functions on $(\sS, d_{\GHP})$. Using this fact and the above bound we have the following proposition:
\begin{prop}
	\label{prop:scale-operator}
Let $\set{\alpha_n:n\geq 1}$ and $\set{\beta_n: n\geq 1}$ be two sequences of positive numbers.	Further assume $\lim_{n \to \infty} \alpha_n = \alpha >0$ and $\lim_{n \to \infty}\beta_n = \beta >0$. Let $\set{X_n: n \geq 1} \subset \sS$ be a sequence of metric spaces such that $X_n \to X \in \sS$ in the metric $d_{\GHP}$ as $n \to \infty$. Then
	\begin{equation*}
		\scl(\alpha_n, \beta_n) X_n \to \scl(\alpha, \beta) X, \mbox{ in $d_{\GHP}$ as } n \to \infty.
	\end{equation*}
\end{prop}
As in the above proposition and the rest of this paper, we will always equip $\sS$ with the topology generated by $d_{\GHP}$.

\subsection{Graphs, trees and ordered trees}
\label{sec:graph-tree-prelim}
All graphs $\cG$ in this study will be simple undirected graphs. We will typically write $V(\cG)$ for the vertex set of $\cG$ and $E(\cG)$ for the corresponding edge set. We will write an edge as $e=(u,v)\in E(\cG)$ with the understanding that $(u,v)$ represents an undirected edge and is equal to $(v,u)$. As before we write $[n] = \set{1,2,\ldots, n}$. We will typically denote a connected component of $\cG$ by $\cC \subseteq \cG$.
 A connected component $\cC$, can be viewed as a metric space by imposing the usual graph distance $d_{\text{G}}$ namely
\begin{equation*}
	d_{\text{G}}(v,u) = \mbox{ number of edges on the shortest path between $v$ and $u$}, \qquad u,v\in \cC.
\end{equation*}
Recall that to construct the random graph, we started with a collection of vertex weights $\set{w_i: i\in [n]}$. Thus there are two natural measures for a connected graph $\cG$ with associated vertex weights $\vw:=\set{w_v: v\in \cG}$:
\begin{enumeratei}
	\item {\bf Counting measure:} $\mu_{\text{ct}}(A): = |A|$, for $A \subset V(\cG)$.
	\item {\bf Weighted measure:} $\mu_\vw(A) :=  \sum_{v \in A} w_v$, for $A \subset V(\cG)$. If no weights are specified then the default convention is to take $w_v \equiv 1$ for all $v$ thus resulting in $\mu_{\vw} = \mu_{\text{ct}}$.
\end{enumeratei}
 For a fixed finite connected graph $\cG$ equipped with vertex weights $\set{w_v:v\in \cG}$, one obtains a compact measured metric space $(V(\cG), d_{\text{G}}, \mu)$, where $\mu$ is either $\mu_{\text{ct}}$ or $\mu_\vw$.  We use $\cG$ for both the graph and the corresponding measured metric space.

A tree $\vt$ is a connected graph with no cycles. A \emph{rooted} tree is a pair $(\vt, r)$ where $\vt$ is a tree and $r \in V(\vt)$ is a distinguished vertex referred to as the root. All trees in the sequel will be rooted trees. Thinking of $r$ as the original progenitor of a genealogy, for vertices in the tree the notions \emph{parent, children, ancestors, siblings, generations} and \emph{heights} have their usual interpretation. An \emph{ordered} tree is a rooted tree in which an order is specified amongst the children of each vertex so that one can talk about the first child, the second child etc. Such trees will be represented as $(\vt,\mvpi)$, where $\vt$ is a rooted tree and $\mvpi$ is the corresponding order. These trees will sometimes be referred to as planar trees as they can be embedded in the plane, arranging children of each vertex from left to right in increasing value of the order.  

For $n \in \bN$, write $\bG_n$, $\bG_n^{\con}$, $\bT_n$ and $\bT_n^{\ord}$ for the collection of all graphs, connected graphs, rooted trees and ordered trees, respectively, with vertex set $[n]$. For ease of notation, we suppress $\mvpi$ in the pair $(\vt, \mvpi)$ for ordered trees and just write $\vt \in \bT_n^{\ord}$. Planar trees can be treated as connected graphs by forgetting about the root and order. Therefore all notation for graphs apply to rooted ordered trees as well. In particular, any tree $\vt$ can be viewed as a measured metric space in $\sS$.




\section{Results}
\label{sec:results}
%
%
%
%

We are now in a position to describe our main results. Section \ref{sec:res-rank-one} describes our main results for the rank-one model (and the associated Chung-Lu model and Britton-Deijfen-Martin-L\"{o}f model) in the critical regime. In Section \ref{sec:res-p-trees} we describe tail bounds on the diameter of $\vp$-trees which play a crucial role in the proof of convergence in the $l^4$ metric.

\subsection{Scaling limits for the rank-one random graph at criticality}
\label{sec:res-rank-one}
We start by stating the assumptions on the weight sequence $\vw$ used to construct the random graph $\cG_n^{\nr}(\vw,\lambda)$. Note that throughout $w_i = w_i(n)$ might depend on $n$ but we suppress this dependence.  Define
\begin{align*}
	&\sigma_k(n) := n^{-1} \sum_{i=1}^n w_i^k  \qquad \mbox{ for } k=1,2,3,\\
	&w_{\max} := \max_{i \in [n]} w_i, \qquad \mbox{ and }  w_{\min} := \min_{i \in [n]} w_i.
\end{align*}

We make the following assumptions on our weight sequence:
\begin{ass} \ 
\label{ass:wt-seq}
	\begin{enumeratea}
		\item {\bf Convergence of three moments: }  There exist constants $\sigma_k >0$ for $k=1,2,3$ such that
			\begin{align*}
				\max\set{	n^{1/3}|\sigma_1(n) - \sigma_1|,~ n^{1/3}|\sigma_2(n) - \sigma_2|~, |\sigma_3(n) - \sigma_3| } \to 0  \mbox{ as } n \to \infty.
			\end{align*}
		\item {\bf Critical regime:} $\sigma_1 = \sigma_2$.
		\item {\bf Bound on the maximum}: There exists $\eta_0 \in (0,1/6)$ such that $w_{\max}= o( n^{1/6 - \eta_0})$.
		\item {\bf Bound on the minimum}: There exists $\gamma_0 > 0$ such that $ 1/w_{\min} = o( n^{\gamma_0})$. Thus minimal weights decrease at most polynomially fast to zero.
	\end{enumeratea}
\end{ass}

For convergence in the $\cT_2$ topology, we will need the following stronger assumption on $w_{\max}$:
\begin{ass} 
\label{ass:high-moment}
There exists $\eta_1 \in (0,1/48)$ such that $w_{\max}= o( n^{1/48 - \eta_1})$.
\end{ass}




For fixed $\lambda\in \bR$, recall the rank-one random graph $\cG_n^{\nr}(\vw,\lambda)$ defined below \eqref{eqn:pij-def}. For $i\geq 1$, let $\cC_n^{\sss(i)}(\lambda)$, denote the $i$-th largest component of $\cG_n^{\nr}(\vw,\lambda)$; for the rest of the sequel to ease notation we will suppress dependence on $\lambda$ and write $\cC_n^{\sss(i)}$ instead of $\cC_n^{\sss(i)}(\lambda)$.  If the number of components is less than $i$, define $\cC_n^{\sss(i)}$ to be a degenerate measured metric space with $\diam(\cC_n^{\sss(i)}) = \mass(\cC_n^{\sss(i)}) = 0$. As described in Section \ref{sec:intro-conn-ph-trans}, for the rank-one model in the critical regime, the number of vertices in $\cC_n^{\sss(i)}$ is of order $n^{2/3}$ (we describe the precise limit result in Section \ref{sec:descp} for each $i \in [n]$). Equipping $\cC_n^{\sss(i)}(\lambda)$ with the graph distance metric and assigning weight $w_v$ to each vertex $v$, we view each of these components as measured metric spaces (see Section \ref{sec:graph-tree-prelim}).
 Our main result is about the limit of the scaled metric spaces defined as
\begin{equation*}
	\vM^n(\lambda) := ( \scl(n^{-1/3}, n^{-2/3}) \cdot \cC_n^{\sss(i)}(\lambda): i \geq 1),
\end{equation*}
namely, rescaling graph distance by $n^{-1/3}$ and each of the weights by $n^{-2/3}$.

\begin{thm} 
	\label{thm:inhom-random-graph}
	Fix $\lambda\in \bR$. Consider the rank-one random graph model $\cG_n^{\nr}(\vw,\lambda)$ with weight sequence $\vw$.
	\begin{enumeratei}
		\item Under Assumption \ref{ass:wt-seq},  as $n \to \infty$,
		\begin{equation}\label{eqn:prod-topology-convergence}
		\vM^n(\lambda) \weakc \vM(\lambda)
		\end{equation}
		where $\vM(\lambda)=(M_i(\lambda): i \geq 1)$ is an $\cM$-valued random variable and
		the convergence takes place with respect to the product topology $\cT_1$. The construction of $\vM(\lambda)$ is given in Section \ref{sec:descp}.
		\item Under the additional Assumption  \ref{ass:high-moment}, the above convergence takes place with respect to the $\cT_2$ topology as in \eqref{eqn:l4-dist-def}.	
	\end{enumeratei}
\end{thm}

\begin{rem}
	Write $\cG^{\er}(n,p)$ for the \erdos random graph with vertex set $[n]$ and connection probability $p$. The critical scaling window corresponds to $p = 1/n+\lambda/n^{4/3}$ with $\lambda\in \bR$ (\cite{erd6s1960evolution,bollobas1984evolution,luczak1990component,luczak1994structure, janson1993birth}). Write $\cC_n^{\sss(i), \er}(\lambda)$ for the $i$-th largest component of $\cG^{\er}(n, 1/n + \lambda/n^{4/3})$, and view these components as measured metric spaces via the counting measure. Define $\vM^{n, \er}(\lambda)$ as
	\[\vM^{n, \er}(\lambda) := \left( \scl(n^{-1/3}, n^{-2/3}) \cdot \cC^{\sss(i), \er}_{n}(\lambda): i\geq 1\right). \]
	Building on the analysis by Aldous \cite{aldous1997brownian} on the size and surplus of components at criticality, Addario-Berry, Broutin and Goldschmidt showed in \cite[Theorem 2]{addario2012continuum} and \cite[Section 4]{addario2013scaling} that, as $n \to \infty$,
	\begin{equation}
	\label{eqn:adb-limit}
		\vM^{n, \er}(\lambda) \weakc \vM^{\er} (\lambda) = (M_i^{\er}(\lambda): i \geq 1) \in \cM,
	\end{equation}
where $\vM^{\er}(\lambda)$ is described in detail in Section \ref{sec:descp}, and convergence is with respect to the $\cT_2$ topology. It will be shown in Lemma \ref{lem:irg-erg-scaling} that the limiting metric spaces in Theorem \ref{thm:inhom-random-graph} satisfy
$$\vM(\lambda)\stackrel{d}{=}\scl\left(\frac{\sigma_1}{\sigma_3^{2/3}}, \frac{\sigma_1}{\sigma_3^{1/3}}\right) \cdot\vM^{\tt er}(\lambda\sigma_1/\sigma_3^{2/3}).$$
This actually shows that under the assumptions of Theorem \ref{thm:inhom-random-graph}, critical rank-one inhomogeneous random graphs viewed
as metric spaces belong to the \erdos universality class.


\end{rem}

The following corollary gives a simple choice of weights which satisfy the relevant assumptions.

\begin{cor}
	Let $\set{w_i: i \in [n]}$ be \emph{i.i.d.} copies of a strictly positive random variable $W$ that satisfies
	$$\E(W)= \E(W^2), \quad \lim_{x \downarrow 0} x^{-\epsilon}\pr(W \leq x) = 0 \mbox{ for some } \epsilon > 0.$$
	Conditional on the weights $\vw = \set{w_i: i \in [n]}$ construct $\cG_n^{\nr}(\vw, \lambda)$ as above.
	\begin{enumeratei}
		\item Assume $\E W^{6+\epsilon} < \infty$. Then the maximal components in $\cG_n^{\nr}(\vw,\lambda)$ satisfy \eqref{eqn:prod-topology-convergence} where the convergence holds in the $\cT_1$ topology.
		\item Assume $\E W^{48+\epsilon} < \infty$.   Then the convergence in \eqref{eqn:prod-topology-convergence} holds in the $\cT_2$ topology.
	\end{enumeratei}
\end{cor}

Now write $D_n$ for the diameter of the graph $\cG_n^{\nr}(\vw,\lambda)$, namely the largest graph distance between two vertices in the same component in $\cG_n^{\nr}(\vw,\lambda)$. Once one is able to prove convergence in the $l^4$ metric, as in \cite{addario2012continuum}, one gets asymptotics for the diameter as well:
\begin{thm}
	\label{thm:diam}
	Assume that the weight sequence satisfies Assumptions \ref{ass:wt-seq} and \ref{ass:high-moment}. Then
	\[\frac{D_n}{n^{1/3}} \weakc \Xi_\infty\]
	where $\Xi_\infty$ is a positive random variable that has an absolutely continuous distribution.
\end{thm}

By \cite[Corollary 2.12]{janson2010asymptotic}, the Norros-Reittu random graph model is
asymptotically equivalent (in the sense of \cite{janson2010asymptotic}) to the Chung-Lu model and the Britton-Deijfen-Martin-L\"{o}f
model under Assumptions \ref{ass:wt-seq}. Hence, an immediate consequence of Theorem
\ref{thm:inhom-random-graph} and Theorem \ref{thm:diam} is the following corollary:
\begin{thm}
\label{thm:cl-bdl}
In the setup of Theorems \ref{thm:inhom-random-graph} and \ref{thm:diam},
the conclusions hold for the Chung-Lu model and Britton-Deijfen-Martin-L\"{o}f model.
\end{thm}

\subsection{Height of $\vp$-trees}
\label{sec:res-p-trees}
In this section, we define a family of random tree models called $\vp$-trees, which play a key role in the proof and the resulting probability bounds on the height of these trees are of independent interest. We refer the interested reader to \cite{pitman2001random} for a comprehensive survey including their role in linking combinatorial objects such as the Abel-Cayley-Hurwitz multinomial expansions to probability. Fix $m \geq 1$, and a probability mass function $\vp = (p_1, p_2,\ldots, p_m)$ with $p_i > 0$ for all $i\in [m]$. A $\vp$-tree is a random tree in $\bT_m$, with law as follows. For any fixed $\vt \in \bT_m$ and $v\in \vt$, write $d_v(\vt)$, for the number of children of $v$ in the tree $\vt$. Then the law of the $\vp$-tree, denoted by $\pr_{\text{tree}}$, is defined as:
\begin{equation}
\label{eqn:p-tree-def}
	\pr_{\text{tree}}(\vt) = \pr_{\text{tree}}(\vt; \vp) = \prod_{v\in [m]} p_v^{d_v(\vt)}, \qquad \vt \in \bT_m.
\end{equation}
Generating a random $\vp$-tree $\cT\sim \pr_{\text{tree}}$ and then assigning a uniform random order on the children of every vertex $v\in \cT$ gives a random element with law $\pr_{\ord}(\cdot ; \vp)$ given by
\begin{equation}
\label{eqn:ordered-p-tree-def}
	\pr_{\ord}(\vt) = \pr_{\ord}(\vt; \vp) = \prod_{v\in [m]} \frac{p_v^{d_v(\vt)}}{(d_v(\vt)) !}, \qquad \vt \in \bT_m^{\ord}.
\end{equation}
Obviously a $\vp$-tree can be constructed by first generating an ordered $\vp$-tree with the above distribution and then forgetting about the order.


Suppressing $m$ in the notation, define
\begin{equation*}
	\sigma(\vp) = \sqrt{\sum_{i=1}^m p_i^2}, \qquad p_{\max} = \max_{1 \leq i \leq m}p_i, \qquad p_{\min} = \min_{ 1 \leq i \leq m}p_i.
\end{equation*}
We will prove the following tail bound for the height of $\vp$-trees. Let $\cT \in \bT_m$ be a random $\vp$-tree with distribution as in \eqref{eqn:p-tree-def}. Let $\height(\cT)$ be the height of the tree $\cT$.

\begin{thm}[Tail bounds]
	\label{thm:ht-p-tree}
  Assume that there exist $\epsilon_0  \in (0, 1/2)$ and $r_0 \in (2, \infty)$ such that
 \begin{equation}
	 \label{eqn:803}
 	\sigma(\vp) \leq \frac{1}{2^{20}}, \qquad \frac{p_{\max}}{[\sigma(\vp)]^{3/2 + \epsilon_0}} \leq 1, \qquad \frac{[\sigma(\vp)]^{r_0}}{p_{\min}} \leq 1.
 \end{equation}
	  Then for any  $r > 0$, there exists some constant $K_{\ref{thm:ht-p-tree}} = K_{\ref{thm:ht-p-tree}}(r) > 0$, such that
\begin{equation}
	\label{eqn:ht-p-tree-2406}
	\pr\left(\height(\cT) \geq \frac{x}{\sigma(\vp)}\right) \leq \frac{K_{\ref{thm:ht-p-tree}}}{x^{r}}, \qquad \mbox{ for } x \geq 1.
\end{equation}
\end{thm}
\begin{rem} 
A natural question the reader might have are the slightly alien assumptions in \eqref{eqn:803}. We will describe the explicit connection between $\vp$-trees and critical rank-one random graphs while proving the main results but to assuage the reader, we give a high-level overview.  In Section \ref{sec:ptree-dfs-connected-irg}, we will relate a natural notion of a ``depth-first'' spanning tree of the maximal components of $\cG^{\nr}_n(\vw, \lambda)$ to tilted versions of associated $\vp$-trees. For example, the depth-first-search tree of the largest component $\cC_n^{\sss(1)}$ can be related to a $\vp$-tree with parameters
\begin{equation*}
	\vp =  (p_v: v \in \cC_n^{\sss (1)}) := \left( \frac{w_v}{\sum_{u\in \cC_n^{\sss(1)}} w_u} : v \in \cC_n^{\sss(i)} \right).
\end{equation*}
Our proofs coupled with Assumption \ref{ass:wt-seq} will imply that \eqref{eqn:803} is satisfied with high probability. 
 Another important implication of the above theorem is that it gives a uniform tail bound for all $\vp$ satisfying \eqref{eqn:803}. Since the diameter of a component is bounded by the height of its associated depth-first tree, the bound in Theorem \ref{thm:ht-p-tree} makes it possible to control the diameter for all components in $\cG^{\nr}_n(\vw, \lambda)$ uniformly and prove convergence in the $\cT_2$ topology in Theorem \ref{thm:inhom-random-graph}.
\end{rem}

\section{Description of limit objects}
\label{sec:descp}
In this section we describe the limit objects $\vM(\lambda)$ arising in Theorem \ref{thm:inhom-random-graph}, first constructed for the \erdos random graph in \cite{addario2012continuum}. We need the following three ingredients:
\begin{enumeratei}
	\item {\bf Real trees:} An abstract notion of ``tree-like'' metric spaces.
	\item {\bf Shortcuts:} A procedure describing when and where to identify points in the real tree to take into account the fact that maximal components in the critical regime may not be trees and could have non-zero surplus or complexity.
	\item {\bf Tilted Brownian excursions:} We will need Brownian excursions whose lengths are described by the limit of component sizes (appropriately rescaled) of the rank-one model as proven in \cite{SBVHJL10},  tilted in favor of excursions with ``large area''.
\end{enumeratei}

\subsection{Real trees and shortcuts}
\label{sec:real-tree-shortcut-limit-object}
A compact metric space $(X,d)$ is called a \emph{real tree} \cite{legall-book,evans-book} if between every two points there is a unique geodesic such that this path is also the only non self-intersecting path between the two points. Functions encoding excursions from zero can be used to construct such metric spaces which we now describe.

For $0 < a < b <\infty$, an \emph{excursion} on $[a,b]$ is a continuous function $h \in C([a,b], \bR)$ with $h(a)=0=h(b)$ and $h(t) > 0$ for $t \in (a,b)$. The length of such an excursion is $b-a$. For $l \in(0,\infty)$, let $\cE_l$ be the space of all excursions on the interval $[0,l]$. Given an excursion $h \in \cE_l$, one can construct a real tree as follows. Define the pseudo-metric $d_h$ on $[0,l]$:
\begin{equation}
\label{eqn:d-pseudo}
	d_h(s,t):= h(s) + h(t) - 2 \inf_{u \in [s,t]}h(u), \; \mbox{ for } s,t  \in [0,l].
\end{equation}
Define the equivalence relation $s \sim t \Leftrightarrow d_h(s,t) = 0$. Let $[0,l]/\sim$ denote the corresponding quotient space and consider the metric space $\cT(h):= ([0,l]/\sim, \bar d_h)$, where $\bar d_h$ is the metric on the equivalence classes induced by $d_h$. Then $\cT(h)$ is a real tree (\cite{legall-book,evans-book}). Let $q_h:[0,l] \to \cT(h)$ be the canonical projection and write $\mu_h$ for the push-forward of the Lebesgue measure on $[0,l]$ onto $\cT(h)$ via $q_h$. Equipped with $\mu_h$, $\cT(h)$ is now a measured compact metric space.

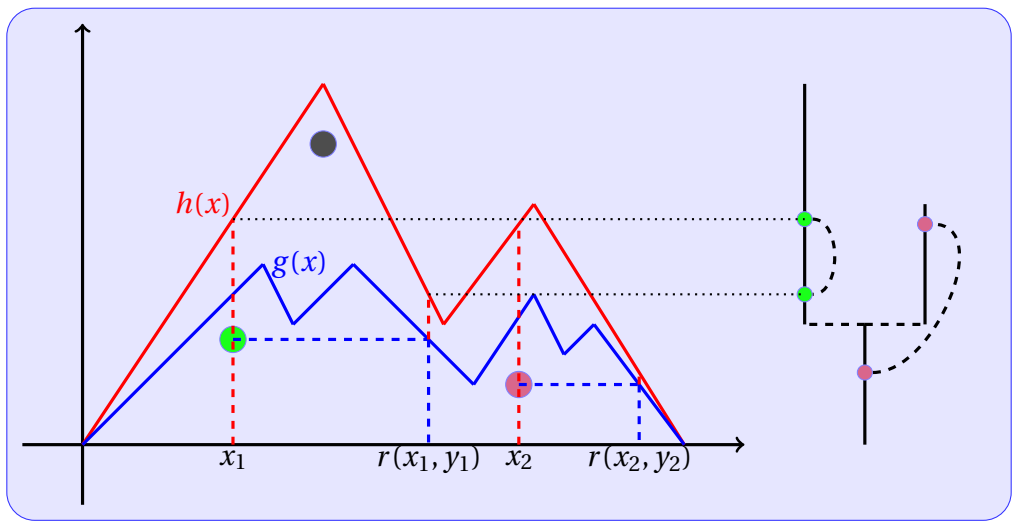
\begin{figure}[htbp]
	\centering
			
		 \begin{tikzpicture}[scale=.8,background rectangle/.style=
		    {draw=blue!80,fill=blue!10,rounded corners=2ex},
		     show background rectangle]
	
		 	\draw[black, very thick, ->] (0,-1) -- (0,7);
		 	\draw[black, very thick, ->] (-1,0) -- (11,0);
		 	\draw[red, very thick] (0,0) -- (4,6);
		 	\draw[red, very thick] (4,6) -- (6,2);
		 	\draw[red, very thick] (6,2) -- (7.5,4);
		 	\draw[red, very thick] (7.5,4) -- (10,0);
		 	\node at (2,4) [red] {$h(x)$};
	
		 	\draw[blue, very thick] (0,0) -- (3,3);
		 	\draw[blue, very thick] (3,3) -- (3.5,2);
		 	\draw[blue, very thick] (3.5,2) -- (4.5,3);
		 	\draw[blue, very thick] (4.5,3) -- (6.5,1);
		 	\draw[blue, very thick] (6.5,1) -- (7.5,2.5);
		 	\draw[blue, very thick] (7.5,2.5) -- (8,1.5);
		 	\draw[blue, very thick] (8,1.5) -- (8.5,2);
		 	\draw[blue, very thick] (8.5,2) -- (10,0);
		 	\node at (3.6,3) [blue] {$g(x)$};
	
		 	\node at (2.5,1.75) [circle,draw=blue!50,fill=green!90] {};
		 	\draw[ red, dashed, very thick] (2.5,0) -- (2.5,3.75);
		 	\node at (2.5,-.25) [] {$x_1$};
		 	\draw[ blue, dashed, very thick] (2.5,1.75) -- (5.75,1.75);
		 	\draw[ blue, dashed, very thick] (5.75,1.75) -- (5.75,0);
		 	\node at (5.75,-.25) [] {$r(x_1,y_1)$};
		 	\draw[ red, dashed, very thick] (5.75,1.75) -- (5.75,2.5);
	
		 	\node at (7.25,1) [circle,draw=blue!50,fill=purple!60] {};
		 	\draw[ red, dashed, very thick] (7.25,0) -- (7.25,3.667);
		 	\node at (7.25,-.25) [] {$x_2$};
		 	\draw[ blue, dashed, very thick] (7.25,1) -- (9.25,1);
		 	\draw[ blue, dashed, very thick] (9.25,1) -- (9.25,0);
		 	\node at (9.25,-.25) [] {$r(x_2,y_2)$};
		 	\draw[ red, dashed, very thick] (9.25,1) -- (9.25,1.2);

		 	\node at (4,5) [circle,draw=blue!50,fill=black!70] {};

	
		 	\draw[black, very thick] (13,0) -- (13,2);
		 	\draw[black, very thick, dashed] (13,2) -- (12,2);
		 	\draw[black, very thick] (12,2) -- (12,6);
		     \node (sca) at (12, 3.75) [circle,draw=blue!50,fill=green!90, inner sep=0pt,minimum size=2mm] {};	
		     \node (scb) at (12, 2.5) [circle,draw=blue!50,fill=green!90, inner sep=0pt,minimum size=2mm] {};	
		 	\draw [black, very thick, dashed] (sca)        to [out=360,in=360] (scb);
		 	\draw[black, dotted, thick] (2.5, 3.75) -- (12, 3.75);
		 	\draw[black, dotted, thick] (5.75, 2.5) -- (12, 2.5);

		 	\draw[black, very thick, dashed] (13,2) -- (14,2);
		 	\draw[black, very thick] (14,2) -- (14,4);
		 	\node (scc) at (14, 3.667) [circle,draw=blue!50,fill=purple!60, inner sep=0pt,minimum size=2mm] {};	
		 	\node (scd) at (13, 1.2) [circle,draw=blue!50,fill=purple!60, inner sep=0pt,minimum size=2mm] {};	
		 	\draw [black, very thick, dashed] (scd)        to [out=360,in=360] (scc);

		 \end{tikzpicture}

		\caption{The above shows a pictorial construction of the metric space $\cG(h,g,\cP)$ starting from two given functions $g, h$. In the left panel, $g\cap \cP$ contains two points of the Poisson process $\cP$ (while one point, colored black, of $\cP$ not contained is also represented; this has no effect on the construction). The right panel shows $\cG(h,g,\cP)$. Here pairs of points connected by dotted lines are identified. }
		
		\label{fig:ghp-construction}
		
\end{figure}

Since our limit objects will not necessarily be trees, we define a procedure that incorporates ``short cuts'' (more precisely identification of points) on a real tree. Let $h,g \in \cE_l$ be two excursions, and $\cP \subseteq \bR_+\times \bR_+$ be a countable set with
\begin{equation*}
	g \cap \cP := \set{(x,y) \in \cP: 0 \leq x \leq l, \; 0 \leq y < g(x)  }.
\end{equation*}
Using these three ingredients, construct a metric space $\cG(h,g,\cP)$ as follows. Let $\cT(h)$ be the real tree associated with $h$ and $q_h: [0,l] \to \cT(h)$ be the canonical projection.  Note that $|g \cap \cP| < \infty$ and write $g \cap \cP  = \set{(x_i, y_i): 1\leq i\leq k}$ for some $k< \infty$. For each $i \leq k$, define
\begin{equation}
	\label{eqn:r-x-y-def}
	r(x_i,y_i) := \inf\set{x: x \geq x_i, \; g(x) \leq y_i}.
\end{equation}
 For $1 \leq i \leq k$, identify the points $q_h(x_i)$ and $q_h(r(x_i,y_i))$ in $\cT(h)$. Call the resulting metric space $\cG(h,g,\cP)$. Equipping this metric space with the push forward of the measure $\mu_h$ on $\cT(h)$ makes $\cG(h,g,\cP)$ a measured compact metric space. $\cG(h,g,\cP)$ can be viewed as the metric space obtained by adding $k$ shortcuts in $\cT(h)$, with the location of the shortcuts determined by the excursion $g$ and the collection of points $\cP$. A shortcut between two points $u, v \in \cT(h)$ identifies these two points as a single point. See Figure \ref{fig:ghp-construction} for an explicit example of this construction.

In principle, one can use the same function for both $h$ and $g$ in the above construction of $\cG(h,g,\cP)$. We retain the flexibility of using different $h$ and $g$, since we will always want $\cP$ to be a \textbf{rate one} Poisson point process, and as a result in many application settings $h$ and $g$ differ by a constant multiple.

\subsection{Scaling limits for component sizes of the Norros-Reittu model}
\label{sec:com-size-nr-limit-object}
Let $\set{B(s):s\geq 0}$ be a standard Brownian motion. For $\kappa,\sigma \in (0,\infty)$ and $\lambda \in \bR$, define
\begin{equation}
\label{eqn:bm-lamb-kapp-def}
W_{\kappa,\sigma}^\lambda(s) := \kappa B(s) + \lambda s - \sigma\frac{s^2}{2}, \qquad s\geq 0.
\end{equation}
Define the reflected process
\begin{equation}
\label{eqn:reflected-pro-def}
	\bar{W}_{\kappa,\sigma}^\lambda(s) := W_{\kappa,\sigma}^\lambda(s) -\min_{0\leq u\leq s} W_{\kappa,\sigma}^\lambda(u),   \qquad s\geq 0.
\end{equation}
 Consider the metric space
\[ \ldown:= \set{\vx= (x_i:i\geq 1): x_1\geq  x_2 \geq \ldots \geq 0, \sum_{i=1}^\infty x_i^2 < \infty},\]
equipped with the natural metric inherited from $l^2$.
It was shown by Aldous in \cite{aldous1997brownian} that the excursions of $\bar{W}_{\kappa,\sigma}^\lambda$ from zero can be arranged in decreasing order of their lengths as
\begin{equation}
\label{eqn:mvxi-def}
	\mvxi_{\kappa,\sigma}^\lambda = (\xi_{\kappa,\sigma}^\lambda(i): i \geq 1),
\end{equation}
where $\xi_{\kappa,\sigma}^\lambda(i)$ denotes the length of the $i$-th largest excursion, and further $\mvxi_{\kappa,\sigma}^\lambda \in \ldown$.

Consider the critical rank-one model $\cG_n^{\nr}(\vw,\lambda)$ with connection probabilities as in \eqref{eqn:pij-def}. Let $|\cC_n^{\sss(i)}(\lambda)|$ be the size of the $i$-th largest component, for $i \geq 1$. In \cite{SBVHJL10} the following was shown about the normalized component sizes under finite third-moment assumptions on the weight sequence:
\begin{thm}[\cite{SBVHJL10}]
	\label{thm:comp-sizes-b-van}
	Fix $\lambda \in \bR$. Under Assumption \ref{ass:wt-seq} (a), (b) and (c), as $n \to \infty$,
\[\left(\frac{|\cC_n^{\sss(i)}(\lambda)|}{n^{2/3}}: i\geq 1\right)\weakc \mvxi_{\sqrt{\sigma_3/\sigma_1}, {\sigma_3}/{\sigma_1^{2}}}^\lambda,
\]
where the weak convergence is with respect to the topology generated by the $l^2$ norm.
\end{thm}
\begin{rem}
	It was assumed in \cite{SBVHJL10} that $w_{\max}=o(n^{1/3})$ and that the empirical distribution $\frac{1}{n}\sum_{i=1}^n \delta_{w_i}$ converges in distribution to some limiting distribution $F$. This second assumption was not used in the proof of Theorem \ref{thm:comp-sizes-b-van} other than to express the limit constants $\sigma_i$ in terms of moments of $W\sim F$ and thus can be removed as long as one has Assumptions \ref{ass:wt-seq} (a) and (b). We place stronger assumptions on $w_{\max}$ as in Assumption \ref{ass:wt-seq} (c). 
\end{rem}
To ease notation, for the rest of this section we shall suppress the dependence on $\sigma_1, \sigma_3$ and $\lambda$ and write the limiting component sizes as
\begin{equation}
\label{eqn:z-defn}
		\mvxi_{\sqrt{\sigma_3/\sigma_1}, {\sigma_3}/{\sigma_1^{2}}}^\lambda:= \vZ = (Z_i : i \geq 1).
\end{equation}
Using Brownian scaling $\set{a^{-1/2}B(as): s\geq 0} \stackrel{d}{=} \set{B(s): s\geq 0}$ with $a = (\kappa/\sigma)^{2/3}$ implies that, for $s \geq 0$,
	\begin{equation*}
		W_{\kappa,\sigma}^\lambda( a s ) = \frac{\kappa^{4/3}}{\sigma^{1/3}}\left[ a^{-1/2} B(as) + \frac{\lambda}{\kappa^{2/3}\sigma^{1/3}} s - \frac{s^2}{2}  \right]  \stackrel{d}{=} \frac{\kappa^{4/3}}{\sigma^{1/3}} W_{1,1}^{\lambda/\kappa^{2/3}\sigma^{1/3}}(s).
	\end{equation*}
	Thus we have $\mvxi_{\kappa,\sigma}^\lambda  \stackrel{d}{=} (\kappa/\sigma)^{2/3} \mvxi_{1,1}^{\lambda /\kappa^{2/3}\sigma^{1/3}}$. Therefore the limit object in Theorem \ref{thm:comp-sizes-b-van} satisfies
	\begin{equation}
	\label{eqn:brown-scal}
		\vZ \stackrel{d}{=} \frac{\sigma_1}{\sigma_3^{1/3}} \mvxi_{1,1}^{\lambda \sigma_1/\sigma_3^{2/3}}.
	\end{equation}
	This relation will be useful in proving Lemma \ref{lem:irg-erg-scaling} below.

\subsection{Tilted Brownian excursions}
\label{sec:tilt-brownian-exc-limit-object}
For fixed $l> 0$, recall that $\cE_l$ denotes the space of excursions of length $l$. We can treat $\cE_l$ as a subset of $C([0,\infty), [0,\infty)) $ by identifying $h \in \cE_l$ with $g \in C([0,\infty), [0,\infty))$ where  $g(s)=h(s)$ for $s \in [0,l]$ and $g(s)=0$ for $s >l$. Write $\cE := \cup_{l >0} \cE_l$ for the space of all finite-length excursions from zero and equip $\cE$ with the $L^\infty$ norm, namely, $\|h\|_\infty = \sup_{s \in [0,\infty)} |h(s)|$.

Let $\set{\ve_l(s): s \in [0,l]}$ be a standard Brownian excursion of length $l$. For $l > 0$ and $\theta >0$, define the tilted Brownian excursion $\tilde \ve_l^\theta$ as an $\cE$-valued random variable such that for all bounded continuous functions $f : \cE \to \bR$,
\begin{equation}
	\label{eqn:tilt-exc-def}
	\E[f(\tilde \ve_l^\theta)] = \frac{ \E\left[f(\ve_l) \exp\left(\theta\int_0^l \ve_l(s)ds\right) \right] }{\E\left[\exp\left(\theta\int_0^l \ve_l(s)ds\right)\right]}.
\end{equation}
Note that $\ve_l$ and $\tilde \ve_l^\theta$ are both supported on $\cE_l$. Writing $\nu_l$ and $\tilde \nu_l^\theta$ respectively for the law of $\ve_l$ and $\tilde \ve_l^\theta$ on $\cE_l$,  the Radon-Nikodym derivative is given by
\begin{equation*}
	\frac{d{\tilde \nu_l^\theta}}{d\nu_l} (h)
	= \frac{\exp\left(\theta \int_0^l h(s) ds \right)}{\int_{\cE_l} \exp \left(\theta \int_0^l g(s) ds \right) d\nu_l(dg)}, \qquad h\in \cE_l.
\end{equation*}
 We use $\ve(\cdot) = \ve_1(\cdot)$ for the standard Brownian excursion. Similarly, we write $\tilde \ve^{\theta} (\cdot) = \tilde \ve_1^{\theta}(\cdot)$ and $\tilde \ve_l(\cdot) = \tilde  \ve_l^{1}(\cdot)$.



By Brownian scaling, $\set{ \sqrt{a} \ve_l(s/a) : s \in [0, al] } \stackrel{d}{=} \set{\ve_{al}(s): s \in [0,al]}$ for $a >0$, thus $\int_0^{al} \ve_{al}(s)ds \stackrel{d}{=} a^{3/2} \int_0^l \ve_l(s) ds$. Taking $a = \theta^{2/3}$ and applying this to \eqref{eqn:tilt-exc-def} gives
\begin{align*}
	\E[f(\tilde \ve_l^\theta)] = \frac{ \E\left[f( \frac{1}{\sqrt{a}}\ve_{al}(a\cdot) ) \exp\left(\int_0^{al} \ve_{al}(s)ds\right) \right] }{\E\left[\exp\left(\int_0^{al} \ve_{al}(s)ds\right)\right]} = \E \left[ f\left(\frac{1}{\sqrt{a}}\tilde \ve_{al}(a\cdot)\right) \right].
\end{align*}
Thus the tilted excursions obey the following scaling relation:
\begin{equation}
	\label{eqn:tilt-brown-exc-scale}
	\set{\tilde \ve_l^\theta(s) : s \in [0,l]} \stackrel{d}{=} \set{\frac{1}{\theta^{1/3}}\tilde \ve_{\theta^{2/3}l}(\theta^{2/3} s) : s \in [0,l]}.
\end{equation}
Note that the scaling relation below \cite[Equation (20)]{addario2012continuum} has a small typo which says $\| \tilde \ve_l \|_\infty \stackrel{d}{=} \sqrt{l}\| \tilde \ve \|_\infty $ in our notation. Based on \eqref{eqn:tilt-brown-exc-scale}, the correct version should be  $\| \tilde \ve_l \|_\infty \stackrel{d}{=} \sqrt{l}\| \tilde \ve^{l^{3/2}} \|_\infty.$

In general, for any $l,\gamma, \theta >0$,
\begin{equation}
	\label{eqn:tilt-brown-exc-scale-general}
	\set{\tilde \ve_l^{\theta \gamma}(s) : s \in [0,l]} \stackrel{d}{=} \set{\frac{1}{\theta^{1/3}}\tilde \ve_{\theta^{2/3}l}^{\gamma}(\theta^{2/3} s) : s \in [0,l]}.
\end{equation}
We generalize the above definitions to (tilted) excursions with random lengths as follows: Fix $\theta > 0$ and let $L$ be an $\bR_+$-valued random variable. Define $\ve_{L}$ [resp. $\tilde \ve_{L}^\theta$] to be a $\cE$-valued random variable such that for every $l > 0$, we have $\pr( \ve_L \in \cdot \mid L=l ) = P( \ve_l \in \cdot )$ [resp. $\pr( \tilde \ve_L^\theta \in \cdot \mid L=l ) = P( \tilde \ve_l^\theta \in \cdot )$].


\subsection{Construction of the scaling limit}
\label{sec:construct-limit-object}
We are now in a position to describe the scaling limits $\vM(\lambda) = (M_i(\lambda): i\geq 1)$ in Theorem \ref{thm:inhom-random-graph}.
Let $\vZ$ be an $\ldown$-valued random variable as defined in \eqref{eqn:z-defn} representing limits of normalized component sizes. Conditional on $\vZ$, generate a sequence $\vh := (h_i: i\geq 1)$ of independent random excursions in $\cE$ via the prescription
\[h_i \stackrel{d}{=} \tilde \ve_{Z_i}^{\sigma_3^{1/2}/\sigma_1^{3/2}}, \qquad i\geq 1.\]
Let $\vP=(\cP_i: i \in \bN)$ be a sequence of \emph{i.i.d.} rate one Poisson point processes on $\bR_+\times \bR_+$, independent of $(\vZ, \vh)$. Define the metric spaces $M_i(\lambda) \in \sS$, $i \geq 1$ as
\begin{equation}
	\label{eqn:def-limit-mi-lambda}
	M_i(\lambda) := \cG\left(\frac{2\sigma_1^{1/2}}{\sigma_3^{1/2}} h_i, \frac{\sigma_3^{1/2}}{\sigma_1^{3/2}}h_i, \cP_i\right) \stackrel{d}{=} \cG\left(\frac{2\sigma_1^{1/2}}{\sigma_3^{1/2}} \tilde \ve_{Z_i}^{\sigma_3^{1/2}/\sigma_1^{3/2}}, \frac{\sigma_3^{1/2}}{\sigma_1^{3/2}} \tilde \ve_{Z_i}^{\sigma_3^{1/2}/\sigma_1^{3/2}}, \cP_i\right),
\end{equation}
where we recall the construction of the metric spaces $\cG(h,g,\cP)$ using the real trees encoded by excursions $h$ and shortcuts generated by the excursion $g$ and collection of points $\cP$,  as introduced in Section \ref{sec:real-tree-shortcut-limit-object}. Write $\vM(\lambda) = (M_i(\lambda): i\geq 1)$ for the sequence of random metric spaces so constructed. Then this is the asserted continuum limit of the critical components in the Norros-Reittu model in Theorem \ref{thm:inhom-random-graph}.

Now we compare this limit to the limit metric spaces for \erdos random graphs as proved in \cite{addario2012continuum}.

\begin{lemma}\label{lem:irg-erg-scaling}
	The limit objects for the rank-one model satisfy the distributional equivalence
	$$\vM(\lambda)\stackrel{d}{=}\scl\left(\frac{\sigma_1}{\sigma_3^{2/3}}, \frac{\sigma_1}{\sigma_3^{1/3}}\right) \cdot\vM^{\tt er}\left(\frac{\lambda\sigma_1}{\sigma_3^{2/3}}\right),$$
	where for any $\lambda^\prime \in \bR$, $\vM^{\tt er}(\lambda^\prime)$ denote the limit objects for the \erdos random graph $\cG_n^{\tt er}(n,1/n +\lambda^\prime/n^{4/3})$ as constructed in \cite{addario2012continuum}.
\end{lemma}

\noindent{\bf Proof:} Write
\begin{equation}
	\mvxi_{1,1}^\lambda := \mvgamma(\lambda) = (\gamma_i(\lambda): i \in \bN).
\end{equation}
In \cite{addario2012continuum} it is shown that the scaling limits for the \erdos model $\cG^{\er}(n, 1/n + \lambda/n^{4/3})$ is
\begin{equation*}
	\vM^{\er}(\lambda) = ( M_i^{\er}(\lambda) : i \geq 1), \mbox{ where } M_i^{\er}(\lambda) := \cG( 2 \tilde \ve_{\gamma_i(\lambda)}, \tilde \ve_{\gamma_i(\lambda)}, \cP_i).
\end{equation*}
In order to compare $\vM(\lambda)$ and $\vM^{\er}(\lambda)$, we again use Brownian scaling. By \eqref{eqn:brown-scal},
\begin{equation}
	\label{eqn:882}
	Z_i = \frac{\sigma_1}{\sigma_3^{1/3}} \gamma_i(\lambda \sigma_1/\sigma_3^{2/3}).
\end{equation}
By the definition of $\cG(h,g,\cP)$ in Section \ref{sec:real-tree-shortcut-limit-object} and $\scl(\alpha,\beta)$ in Section \ref{sec:topology-prelim}, for $\alpha, \beta > 0$,
\begin{equation}
	\label{eqn:887}
	\scl(\alpha, \beta) \cdot \cG(h,g,\cP) = \cG(\alpha h (\cdot/\beta), \frac{1}{\beta} g(\cdot/\beta), \cP^{\beta}),
\end{equation}
where $\cP^{\beta} := \set{(\beta x, y/\beta) : (x,y) \in \cP}$. With the these ingredients, letting $\bar \theta := \sigma_3^{1/2}/\sigma_1^{3/2}$,
\begin{align*}
	M_i(\lambda)
	=& \cG\left(\frac{2\sigma_1^{1/2}}{\sigma_3^{1/2}} \tilde \ve_{Z_i}^{\bar \theta}, \frac{\sigma_3^{1/2}}{\sigma_1^{3/2}} \tilde \ve_{Z_i}^{\bar \theta}, \cP_i\right) \\
	\stackrel{d}{=}& \cG\left(\frac{2\sigma_1^{1/2}}{\sigma_3^{1/2} \bar \theta^{1/3}} \tilde \ve_{\bar \theta^{2/3}Z_i}(\bar \theta^{2/3} \cdot), \frac{\sigma_3^{1/2}}{\sigma_1^{3/2}\bar \theta^{1/3}} \tilde \ve_{\bar \theta^{2/3}Z_i}(\bar \theta^{2/3} \cdot), \cP_i\right)\\
	=& \cG\left(\frac{2\sigma_1^{1/2}}{\sigma_3^{1/2} \bar \theta^{1/3}} \tilde \ve_{\gamma_i(\lambda \sigma_1/\sigma_3^{2/3})}(\bar \theta^{2/3} \cdot), \frac{\sigma_3^{1/2}}{\sigma_1^{3/2}\bar \theta^{1/3}} \tilde \ve_{\gamma_i(\lambda \sigma_1/\sigma_3^{2/3})}(\bar \theta^{2/3} \cdot), \cP_i\right),
\end{align*}
where the first line is by definition, the second line is due to \eqref{eqn:tilt-brown-exc-scale}, and the third line follows from \eqref{eqn:882}. To ease notation write $\gamma_i = \gamma_i(\lambda \sigma_1/\sigma_3^{2/3})$. Then taking $\bar \alpha = \sigma_3^{2/3}/\sigma_1$,
\begin{align*}
	\scl( \bar \alpha, \bar \theta^{2/3}) \cdot 	M_i(\lambda)
	=& \cG\left(\frac{2\bar \alpha \sigma_1^{1/2}}{\sigma_3^{1/2} \bar \theta^{1/3}} \tilde \ve_{\gamma_i}(\cdot), \frac{\sigma_3^{1/2}}{\sigma_1^{3/2}\bar \theta} \tilde \ve_{\gamma_i}(\cdot), \cP_i^{\bar \theta^{2/3}}\right)\\
	=& \cG\left(2 \tilde \ve_{\gamma_i}(\cdot),  \tilde \ve_{\gamma_i}(\cdot), \cP_i^{\bar \theta^{2/3}}\right)\\
	\stackrel{d}{=}& \cG\left(2 \tilde \ve_{\gamma_i}(\cdot),  \tilde \ve_{\gamma_i}(\cdot), \cP_i\right).
\end{align*}
Here the first line uses \eqref{eqn:887} and the second line follows from the definition of $\bar \alpha$ and $\bar \theta$.  The third line follows from standard properties of Poisson point processes resulting in $\cP_i^{\beta} \stackrel{d}{=} \cP_i$ for all $\beta >0$, and the independence between $\cP_i$ and $h_i$ in the construction. This completes the proof of Lemma \ref{lem:irg-erg-scaling}. \qed

%

%
%
%

%

\section{Discussion}
\label{sec:disc}
Before proceeding to the proofs, let us briefly describe the relevance of these results and their connection to the existing literature on random graphs.
\begin{enumeratea}
	\item {\bf Relevance of the results:} As described in the introduction, it turns out that the proof techniques developed to analyze the particular random graph model in this study turn out to be the key technical tools to understand the continuum limits of the metric structure of maximal components of a host of other models in the critical regime. This is developed further in \cite{SBSSXW-universal}. We will try to give a general idea why this is true.  Many random graph models can be formulated as dynamic random graph processes $\set{\cG_n(t): t\geq 0}$ where edges are added to the system in some model-dependent manner.  Further the critical regime corresponds to a (model-dependent) critical time $t_c$ while the critical scaling window corresponds to times of the form $t_c+\lambda/n^{1/3}$ for fixed $\lambda \in \bR$. Fix $\delta < 1/3$ and call the configuration at time $t_n  = t_c -n^{-\delta}$, the \emph{barely subcritical regime}. For typical applications it turns out that one needs $\delta \in (1/6, 1/5)$, see \cite{SBSSXW-universal}. It turns out that while it is hard to approximate connections (edges) formed all the way from $t=0$ in most random graph models,  connections formed between components in the time window $[t_n, t_c+\lambda/n^{1/3}]$ can be strongly coupled to a rank-one random graph model.  The general framework developed in \cite{SBSSXW-universal} consists of two main ingredients:
	\begin{enumeratei}
		\item Show that components at time $t_n$, the so-called entrance boundary, satisfy good regularity properties.
		\item Show that connections between components in $[t_n, t_c+\lambda/n^{1/3}]$ can be approximated via the rank-one model and then use techniques developed in Section \ref{sec:scale-limit-connected-irg} to analyze the metric space generated by these connections.
	\end{enumeratei}
	This general framework is then used to show that components in the critical regime for the configuration model \cite{BC78,B80,MR98}, inhomogeneous random graph models with finite type space and general kernel $\kappa$ \cite{BBSJOR07} and bounded-size rules \cite{spencer2007birth} all satisfy results analogous to Theorem \ref{thm:inhom-random-graph}.

	\item {\bf Connection to existing results:} As remarked in Section \ref{sec:mot}, the main aim of the paper was to rigorously understand conjectures in statistical physics on scaling limits of distances in the critical regime for inhomogeneous random graph models, which then predict distances for the minimal spanning tree (on the giant component) in the supercritical regime where each edge has Uniform$[0,1]$ edge weights. See \cite{braunstein2003optimal} and the references therein. This entire program has been rigorously carried out for the complete graph see \cite{addario2009critical,addario2012continuum,addario2010critical} for a sequence of results including distance scaling for the maximal components in the critical regime for the \erdos random graph, finally culminating in the scaling limit of the minimal spanning tree of the complete graph equipped with uniform edge weights. Most influential to this study is \cite{addario2012continuum} which constructs the scaling limit of these components in the critical regime. Extending these results to the context of general inhomogeneous random graph models turned out to be challenging since the homogeneous nature of the \erdos played an important role in various parts of the proof in \cite{addario2012continuum}.
	
	While there have been few rigorous results on the actual structure of components in the critical regime, if one were interested in only sizes of the maximal components, then this has witnessed significant progress over the last few years, see \cite{joseph2010component,nachmias2010critical,riordan2012phase,hatami2012scaling} for results on the configuration model, \cite{bhamidi2012augmented} for results on a class of dynamic random graph processes called the bounded-size rules and most relevant for this work \cite{SBVHJL10,turova2009diffusion} for results for the rank-one inhomogeneous random graph. See \cite{bollobas2013phase} for a recent survey.
	
	\item {\bf Importance of the assumptions:} Consider the critical rank-one model studied in this paper. To prove our main results, we needed moment Assumptions \ref{ass:wt-seq}. If one were interested in just proving results on the sizes of components (Theorem \ref{thm:comp-sizes-b-van}), finite third moment assumptions, namely Assumptions \ref{ass:wt-seq}(a),(b) and (c) replaced by $w_{\max} = o(n^{1/3})$ suffice. However to understand the actual metric structure of the component, one is lead to these stronger assumptions owing to technical conditions in \cite{aldous2004exploration}, required to show that in our setting the associated (untilted) $\vp$-trees (properly rescaled by $n^{1/3}$) have scaling limits related to the continuum random tree.  The results in \cite{aldous2004exploration} in fact assume exponential moments, however as in the remark under  \cite[Theorem 3]{aldous2004exploration}, these can be extended without much work following the same proof technique as in \cite{aldous2004exploration} to the setting of finite moment conditions.  We state this extension in Theorem \ref{thm:AMP} below. We believe that in fact these results can be extended all the way to finite third moments and are in the process of understanding how to refine these results further. Given the technical nature of the proof even with sufficient moment assumptions, we defer this to future work. In the context of infinite third moments but finite second moments of the weight sequence where one deals with rank-one models with degree exponent $\tau \in (3,4)$, Proposition \ref{prop:distrib-gone-tildg} in this paper is the starting point in \cite{bhamidi2015multiplicative} to derive completely new novel scaling limits for the metric structure of maximal components. These scaling limits are based on tilted versions of inhomogeneous continuum random trees (ICRTs) \cite{aldous2000inhomogeneous}. In the context of \cite{bhamidi2015multiplicative} there are no known results about the exploration processes (analogs of Theorem \ref{thm:AMP}). Thus starting from Proposition \ref{prop:distrib-gone-tildg}, in \cite{bhamidi2015multiplicative} a completely different approach was required; here first Gromov weak convergence was established and then a global lower mass bound derived to establish conclude the analysis.  	
\end{enumeratea}


\section{Proof preliminaries and outline}
\label{sec:proofs-comp}

%
%
We now commence on the proofs of the main results. In this section we start with an outline of the proof and describe some preliminary properties of the rank-one model.

\subsection{Outline of proof}
\label{sec:proof-idea}
 We start in Section \ref{sec:connected-comp} where we describe how the connected components of the rank-one model can be constructed in two steps. In particular this construction will imply that there are two major parts in understanding the maximal connected components:
\begin{enumeratea}
	\item {\bf Construction and asymptotics for rank-one model conditional on being connected: } In Section \ref{sec:ptree-dfs-connected-irg}, we explore
a random graph generated by the rank-one model in a randomized depth-first manner and show that the law of this  depth-first tree is that of an
\emph{ordered tilted $\vp$-tree} (Proposition \ref{prop:distrib-gone-tildg}). This will imply an alternate way of constructing
a rank-one graph conditioned on being connected: first generate an ordered tilted $\vp$-tree and then add the surplus
edges independently with appropriate probabilities. To understand asymptotics, the proof proceeds in two steps:
\begin{enumeratei}
	\item \emph{Untilted $\vp$-trees with shortcuts:} We start by studying the effect of adding surplus edges to the original $\vp$-trees in Section \ref{sec:scale-limit-connected-irg}.   Strengthening the results of \cite{aldous2004exploration},
it follows that an ordinary $\vp$-tree converges in the Gromov-Hausdorff topology (after the tree distance has been appropriately rescaled)
to a continuum random tree under some regularity conditions on the driving probability mass function $\vp$. A careful analysis culminates in the proof of Proposition \ref{prop:ghp-convg-non-tilt} describing joint convergence of the associated $\vp$-tree with shortcuts as well as an associated functional $L(\cdot)$.
\item \emph{Tilted $\vp$-trees and asymptotics for large connected components:} Provided we can show the corresponding tilt is ``nice'', this would imply that rank-one random graphs conditioned on being connected converge to a tilted continuum tree where certain pairs of points
have been identified. We achieve this with the help of Lemma \ref{lem:wts-permitted}
(which shows that the tilt converges pointwise) and Lemma \ref{lem:fexp-linfty-tail} (which yields uniform integrability
of the tilt).
\end{enumeratei}
\item {\bf Generating connected components in the rank-one model and regularity of vertex weights in the maximal components: } To study this, we start in Section \ref{sec:size-bias-explor} by describing the exploration of the graph $\cG_n^{\nr}(\vw, \lambda)$ in a size-biased random order first used in \cite{aldous1997brownian} and later used in \cite{SBVHJL10} to prove Theorem \ref{thm:comp-sizes-b-van} on the scaling of the maximal components in the critical regime. We use this exploration process to generate the connected components and establish strong regularity properties of the \emph{weights} of vertices in these maximal components (Proposition \ref{prop:weight-control}).
\end{enumeratea}

Now conditional on these regularity properties being satisfied within each maximal component, the internal structure
of such a component is simply that of a rank-one inhomogeneous random graph conditioned on being connected, namely the object analyzed in (a) above. We then combine these two parts to prove convergence of the scaled components in Section \ref{sec:finish-proof}, where we first prove convergence in the product topology $\cT_1$. Extending the convergence of the components in the $\cT_2$ topology ($l^4$ metric) amounts to proving a tail bound on the diameter of
the components. Since the depth-first tree of each component spans the component and is distributed as a tilted $\vp$-tree,
it is enough to get a tail bound on heights of $\vp$-trees. We achieve this in Section \ref{sec:proof-height}
by using techniques from \cite{camarri2000limit, aldous2004exploration}.

\subsection{Connected components of the model}
\label{sec:connected-comp}

Recall that $(\cC_n^{\sss(i)}: i \geq 1)$ are the components of $\cG_n^{\nr}(\vw,\lambda)$ ranked in decreasing order of their size. 
Fix $\cV\subset [n]$ and write $\bG_{\cV}^{\con}$ for the space of all simple connected graphs with vertex set $\cV$.
For fixed $a > 0$, and probability mass function $\vp = (p_v: v \in \cV)$, define probability distributions $\pr_{\con}(\cdot; \vp, a, \cV)$ on $\bG_{\cV}^{\con}$ as follows: Define for $i,j \in \cV$,
\begin{equation}
\label{eqn:qij-def}
	q_{ij}:= 1-\exp(-a p_i p_j).
\end{equation}
Then
\begin{equation}
	\label{eqn:pr-con-vp-a-cV-def}
	\pr_{\con}(G; \vp, a, \cV): = \frac{1}{Z(\vp,a,\cV)} \prod_{(i,j)\in E(G)} q_{ij} \prod_{(i,j)\notin E(G)} (1-q_{ij}),\qquad  \mbox{ for } G \in \bG_{\cV}^{\con},
\end{equation}
where $Z(\vp,a,\cV)$ is the normalizing constant
$$Z(\vp,a,\cV) := \sum_{G \in \bG_{\cV}^{\con}}{\prod_{(i,j)\in E(G)} q_{ij} \prod_{(i,j)\notin E(G)}(1-q_{ij})}.$$
We will refer to this distribution as the rank-one random graph \emph{conditioned on being connected}.
Now let $\cV^{\sss(i)} := V(\cC_n^{\sss(i)})$ be the vertex set of $\cC_n^{\sss(i)}$ for $i \geq 1$ and note that $(\cV^{\sss(i)}: i\geq 1)$ denotes a random finite partition of the complete vertex set $[n]$. The next proposition characterizes the distribution of the random graphs $(\cC_n^{\sss(i)}: i \geq 1)$ conditioned on the partition $(\cV^{\sss(i)}: i \geq 1)$:
\begin{prop}
	\label{prop:generate-nr-given-partition}
	For $i \geq 1$ define
	\begin{equation}
		\label{eqn:838}
		\vp^{\sss(i)} := \left( \frac{w_v}{\sum_{v \in \cV^{\sss(i)}}w_v } : v \in \cV^{\sss(i)} \right), \;\;  a^{\sss(i)}:= \left(1+ \frac{\lambda}{n^{1/3}}\right)\frac{(\sum_{v \in \cV^{\sss(i)}}w_v )^2}{l_n}.
	\end{equation}
	For each fixed $i \geq 1$, let $G_i \in  \bG_{\cV^{\sss(i)}}^{\con}$ be connected simple graphs with vertex set $\cV^{\sss(i)}$. Then
	\begin{equation*}
		\pr\left(\cC_n^{\sss(i)} = G_i, \;\; \forall i \geq 1 \mid (\cV^{\sss(i)}: i \geq 1) \right) = \prod_{i\geq 1} \pr_{\con}( G_i; \vp^{\sss(i)}, a^{\sss(i)}, \cV^{\sss(i)}).
	\end{equation*}
\end{prop}
\noindent\textbf{Proof: } Let $\set{V_i : i \in [k]}$, $k \in \bN$, be a partition of $[n]$ such that $|V_1| \geq \cdots \geq |V_k| > 0$. In this proof we will fix such a partition $\set{V_i : i \in [k]}$ and define, for $i \in [k]$, $\vp^{\sss(i)} := (p^{\sss(i)}_v: v \in V_i)$ and $a^{\sss(i)}$ as in \eqref{eqn:838}, but with $\cV^{\sss(i)}$ replaced by $V_i$. It is sufficient to show that
\begin{equation} \label{eqn:1180-0}
	\frac{\pr\left(\cC_n^{\sss(i)} = G_i, \cV^{\sss(i)} = V_i, \;\; \forall i \in [k]  \right)}{\prod_{i \in [k]} \pr_{\con}( G_i; \vp^{\sss(i)}, a^{\sss(i)}, V_i)} \mbox{ does not depend on } (G_i \in \bG_{V_i}^{\con}: i \in [k]).
\end{equation}
Indeed, by \eqref{eqn:pij-def} and the definition of $\vp^{\sss(i)}$ and $a^{\sss(i)}$, for $v,u \in [n]$,
\begin{equation}
	\label{eqn:1180-1}
	p_{vu}(\lambda)= 1-\exp\left(-\left(1+\frac{\lambda}{n^{1/3}}\right)\frac{w_v w_u}{l_n}\right) = 1 - \exp \left( -a^{\sss(i)} p_v^{\sss(i)}p_u^{\sss(i)}\right).	
\end{equation}
We suppress the dependency on $\lambda$ and write $p_{vu} = p_{vu}(\lambda)$. Denote by $[n]_2$ (resp. $E_i$, for $i \in [k]$) the collection of all possible edges on a graph with vertex set $[n]$ (resp. $V_i$, for $i\in [k]$). Then
\begin{align}
	&\pr\left(\cC_n^{\sss(i)} = G_i, \cV^{\sss(i)} = V_i, \;\; \forall i \in [k]  \right) \nonumber \\
	=& \prod_{(v,u) \in E(G_1)\cup \cdots \cup E(G_k)} p_{vu} \prod_{(v,u) \in [n]_2 \setminus (E(G_1)\cup \cdots \cup E(G_k))} (1-p_{vu}) \nonumber \\
	=& \prod_{(v,u) \in [n]_2 \setminus (E_1 \cup \cdots \cup E_k)} (1-p_{vu}) \; \; \prod_{i \in [k]} \left( \prod_{(v,u) \in E(G_i)} p_{vu}\prod_{(v,u) \in E_i \setminus E(G_i)}(1-p_{vu}) \right), \label{eqn:1180-2}
\end{align}
where the last line uses the following relation:
\begin{equation*}
	[n]_2 \setminus (E(G_1)\cup \cdots \cup E(G_k)) = \left[[n]_2 \setminus (E_1 \cup \cdots \cup E_k)\right] \cup \left[\cup_{ i \in [k]}(E_i \setminus E(G_i))\right].
\end{equation*}
On the other hand, by \eqref{eqn:pr-con-vp-a-cV-def} and \eqref{eqn:1180-1},
\begin{equation*}
	\prod_{i \in [k]} \pr_{\con}( G_i; \vp^{\sss(i)}, a^{\sss(i)}, V_i) = \prod_{i \in [k]} \left( \prod_{(v,u) \in E(G_i)} p_{vu}\prod_{(v,u) \in E_i \setminus E(G_i)}(1-p_{vu}) \right)  \prod_{i \in [k]} \frac{1}{Z(\vp^{\sss(i)},a^{\sss(i)}, V_i)}.
\end{equation*}
Comparing the above display and \eqref{eqn:1180-2}, we complete the proof of \eqref{eqn:1180-0}. \qed

\ \\
Under the assumptions on $\vw$, it will be shown in Section \ref{sec:finish-proof} that for each fixed $i\geq 1$,  $a^{\sss(i)} = O_P(n^{1/3})$. Now  the above proposition says that the random graph $\cG_n^{\nr}(\vw,\lambda)$ can be generated in two stages:
\begin{enumeratei}
	\item In the first stage generate the partition of the vertices
into different components, i.e. $(\cV^{\sss(i)}: i \geq 1)$. This can be achieved by generating the full graph of $\cG_n^{\nr}(\vw,\lambda)$, recording the partition of the vertices, and then erasing all the edges.
\item In the second stage, given the partition, we generate the internal structure of each component following the law of $\pr_{\con}(\cdot ; \vp^{\sss(i)}, a^{\sss(i)}, \cV^{\sss(i)})$, independently across different components.
\end{enumeratei}
	
	The plan of the next section is to study the technically more challenging question of the rank-one random graph conditioned on being connected. 

\section{Scaling limit of rank-one graphs conditioned on being connected}
\label{sec:scale-limit-connected-irg}


Recall Proposition \ref{prop:generate-nr-given-partition} where conditional on the partition of the vertices into their respective connected components resulted in the probability distributions $\pr_{\con}(\cdot; \vp^{\sss(i)}, a^{\sss(i)}, \cV^{\sss(i)})$ on the space of connected graphs with vertex set $\cV^{\sss(i)}$ where the parameters $\vp^{\sss(i)}$ and $a^{\sss(i)}$ are as defined in \eqref{eqn:838} with respect to the weights of the vertices in the component $\cV^{\sss(i)}$. To simplify presentation, write $a$, $\vp$ and $\cV$ for the generic setup as above.   This section heralds the entry of $\vp$-trees as defined in Section \ref{sec:res-p-trees}.  We will give an alternate construction of the distribution $\pr_{\con}(\cdot; \vp, a, \cV)$ using random trees whose distribution is given via an appropriate tilt of the original $\vp$-tree distribution, where for component $\cC_n^{\sss(i)}$ the corresponding driving $\vp$-tree model uses the parameters in \eqref{eqn:838}. We will start by describing two assumptions on these driving parameters. These assumptions will be elaborated on in Section \ref{sec:size-bias-explor}  for the maximal components.

Fix $m \geq 1$ and  let $\vp^{\sss(m)} = (p_i^{\sss (m)}: i \in [m])$ be a probability mass function. We will suppress $m$ in the notation and just write $\vp = (p_i: i \in [m])$. Recall the definitions
\[\sigma(\vp) := \sqrt{\sum_{i\in [m]}p_i^2}, \qquad p_{\max} := \max_{i \in [m]} p_i,\qquad p_{\min} := \min_{ i \in [m]} p_i. \]

\begin{ass}
	\label{ass:aldous-AMP}
There exists $\epsilon > 0$ and $r > 0$ such that, as $m \to \infty$,
	\begin{equation*}
		\sigma(\vp) \to 0, \qquad \frac{p_{\max}}{[\sigma(\vp)]^{3/2+\epsilon}} \to 0, \qquad \frac{[\sigma(\vp)]^r}{p_{\min}} \to 0.
	\end{equation*}
\end{ass}

Write $\bG_m^{\con}$ for the collection of all connected graphs with vertex set $[m]$. Let $\set{a(m):m\geq 1}$ be a sequence of positive real numbers. We will use $a=a(m)$ and $\vp$ to construct a probability measure $\pr_{\con}$ on $\bG_m^{\con}$ as in \eqref{eqn:pr-con-vp-a-cV-def} namely
\begin{equation}
\label{eqn:rg-dist-connected}
	\pr_{\con}(G): = \pr_{\con}(G; \vp, a, [m]), \qquad G \in \bG_m^{\con}.
\end{equation}
 Let $\cG_m$ be a $\bG_m^{\con}$-valued random variable with distribution $\pr_{\con}$. Thus $\cG_m$ has the same distribution as a rank-one random graph with vertex set $[m]$ and connection probabilities $q_{ij} = 1-\exp(-ap_ip_j)$, {conditioned on being connected}.
We will think of $\cG_m \in \bG_m^{\con}$ as a measured metric space as described in Section \ref{sec:graph-tree-prelim} using the graph distance as the metric, and assigning mass $p_i$ to vertex $i \in [m]$. The main result of this section shows that under some regularity conditions on $\vp$ and $a$ as $m\to\infty$, the metric space $\cG_m$ with graph distance rescaled by $\sigma(\vp)$ converges to a measured (random) metric space with distribution as described in Section \ref{sec:real-tree-shortcut-limit-object}. In addition to Assumption \ref{ass:aldous-AMP}, we need the following assumption on $a(m)$:
\begin{ass}
	\label{ass:additional-connected}
	There exists constant  $\bar \gamma \in (0,\infty)$ such that
	\begin{equation}
		\label{eqn:ass-asymptotic-am}
	\lim_{m \to \infty}a \sigma(\vp)  = \bar \gamma.
	\end{equation}
\end{ass}
The main aim of this section is to prove the following result:
\begin{thm}
	\label{thm:conv-condition-on-connected}
	 Let $\cG_m$ be a $\bG_m^{\con}$-valued random variable with law $\pr_{\con}$. Under Assumptions \ref{ass:aldous-AMP} and \ref{ass:additional-connected}, as $m \to \infty$,
	\begin{equation*}
		\scl\left(\sigma(\vp), 1\right) \cdot \cG_m \weakc \cG( 2 \tilde \ve^{\bar \gamma}, \bar \gamma \tilde \ve^{\bar \gamma}, \cP),
	\end{equation*}
	where $~\tilde \ve^{\bar \gamma}$ is the tilted Brownian excursion as defined in \eqref{eqn:tilt-exc-def}, $\cP$ is a rate one Poisson point process on $\bR_+^2$ independent of $\tilde \ve^{\bar \gamma}$, and $\cG( 2 \tilde \ve^{\bar \gamma}, \bar \gamma \tilde \ve^{\bar \gamma}, \cP)$ is the random compact measured metric space constructed in Section \ref{sec:real-tree-shortcut-limit-object}.
\end{thm}

The rest of this section is organized as follows. We start in Section \ref{sec:ptree-dfs-connected-irg} where we will give an alternative construction of the law $\pr_{\con}$ from an ordered $\vp$-tree by tilting this distribution appropriately. In Section \ref{sec:conv-height-cont-connected-irg}, we study the scaling limit of a random connected graph without applying the tilt. Finally in Section \ref{sec:connected-tightness-tilt}, we prove the tightness of the tilt and complete the proof of Theorem \ref{thm:conv-condition-on-connected}.

\subsection{Distribution of connected components and tilted $\vp$-trees}
\label{sec:ptree-dfs-connected-irg}

 Recall that $\bT_m^{\ord}$ denotes the space of ordered (planar) trees on $m$ vertices. We start by introducing the following \textbf{randomized Depth First Search (rDFS)} procedure, which takes a graph $G \in \bG_m^{\con}$ as the input and outputs a random ordered tree in $\bT_m^{\ord}$, denote by $\Gamma_{\vp}(G)$, as its output. Given $G \in \bG_m^{\con}$, the rDFS consists of two stages:

\textbf{I. Selection of a root:} Pick $v(1) \in [m]$ at random from the distribution $\vp$. The vertex $v(1)$ is the starting point of the rDFS algorithm on the graph $G$ and also the root of the ensuing tree $\Gamma_{\vp}(G)$.

\textbf{II. Depth-First-Search:} At each step $1\leq i \leq m$, we will keep track of three types of vertices:
\begin{enumeratea}
	\item The set of already explored vertices, $\cO(i)$.
	\item The set of active vertices $\cA(i)$. We view $\cA$ as a vertical \emph{stack} with $~\cA(i)$ denoting the state of the stack at the \emph{end} of the $i$-th step.
	\item The set of unexplored vertices $\cU(i) := [m]\setminus (\cA(i) \cup \cO(i))$.
\end{enumeratea}
Initialize with $\cA(0) = \set{v(1)}$, $\cO(0) = \emptyset$. At step $i \geq 1$, let $v(i)$ denote the vertex on \emph{top} of the stack $\cA(i-1)$ and let
\[\cD(i):= \set{ u \in \cU(i-1): (v(i),u) \in E(G) },\]
namely $\cD(i)$ is the set of unexplored neighbors of $v(i)$. Let $d_{v(i)} = |\cD(i)|$ and suppose $\cD(i) = \set{ u(j) : 1\leq j\leq d_{v(i)}}$. Then update the stack $\cA(\cdot)$ in the following manner:
\begin{enumeratei}
	\item \textbf{Randomization:} Generate $\mvpi = \mvpi(i)$ a uniform random permutation on $[d(i)]$.
	\item Delete $v(i)$ from $\cA (i-1)$.
	\item Arrange the vertices $\cD(i)$ on top of the stack $\cA(i-1)$ using the order $\mvpi$ generated in (i).
\end{enumeratei}
Define $\cA(i)$ to be the state of the stack $\cA$ after the above operations. As sets, $\cA(i) = \cA(i-1) \cup \cD(i) \setminus \set{v(i)}$. Define $\cO(i):= \cO(i-1) \cup \set{v(i)}$ and $\cU(i) := \cU(i-1) \setminus \cD(i)$.

Note that in the above rDFS algorithm, we have $|\cO(i)|=i$ for $i \in [m]$. Thus after $m$ steps we complete the exploration of all vertices in $G$. At the end of the procedure we are left with a rooted random tree $\Gamma_{\vp}(G) \in \bT_m^{\ord}$ with $v(1)$ as the root and with edge set $E(\Gamma_{\vp}(G)) : = \set{ (v(i),u): i \in [m], \; u \in \cD(i)}$. Carrying the order $\set{\mvpi(i): i\in [m]}$ used to order the vertices at each stage of the procedure then makes the resulting tree an ordered tree that we explore in a depth first manner,  resulting in the order $(v(1), \ldots, v(m))$.
This completes the construction of $\Gamma_{\vp}(G) \in \bT_m^{\ord}$. Note that for fixed $G$, $\Gamma_{\vp}(G)$ is a $\bT_m^{\ord}$-valued random variable.

The rDFS algorithm incorporates randomization in two places: First, the root is chosen randomly using the probability mass function $\vp$; second, the children (unexplored vertices) of each vertex are explored in uniform random order. Given an ordered tree $\vt \in \bT_m^{\ord}$, one can run a depth first search on $\vt$ in a deterministic manner starting from the root of the tree $\vt$ and exploring the children using the associated order of the tree. Let $(\cO(i), \cA(i), \cU(i), \cD(i) : i\in [m])$ be the corresponding sets of vertices obtained from this deterministic depth first search of the tree $\vt$. Write $\sP(\vt,\mvpi)$ for the set of edges $\set{u,v}$ such that there exists $1\leq i\leq m-1$ such that $u,v\in \cA(i)$, namely both are active but have not yet been explored. Using terminology from \cite{addario2012continuum} call this collection of edges, the set of \emph{permitted edges}. By definition,
\begin{equation}
	\sP(\vt) := \set{(v(i), j): i \in [m], \;\; j \in \cA(i-1) \setminus \set{v(i)}}. \label{eqn:permit-edges-characterization}
\end{equation}
Write $[m]_2$ for the collection of all possible edges on a graph with vertex set $[m]$ and recall that $E(\vt)$ denotes the edge set of $\vt$. Call the remaining edges i.e.,
\[\sF(\vt) := [m]_2\setminus (\sP(\vt) \cup E(\vt)), \]
the set of \emph{forbidden} edges. 

For a fixed planar tree $\vt \in \bT_m^{\ord}$, define the subset of simple connected graphs $\bG(\vt) \subset \bG_m^{\con}$ as
\begin{equation}
	\label{eqn:def-gtpi}
	\bG (\vt) := \set{ G \in \bG_m^{\con} : E(\vt) \subset E(G) \subset E(\vt) \cup \sP(\vt)}.
\end{equation}

For fixed $G \in \bG_m^{\con}$, let $\nu^{\dfs}(\cdot ; G)$ be the probability distribution of $\Gamma_{\vp}(G)$ on $\bT_m^{\ord}$. When $G \notin \bG (\vt)$, by \cite[Lemma 7]{addario2012continuum}, we have  $\nu^{\dfs}(\vt; G) = 0$. This also explains the terms ``\emph{permitted}'' and ``\emph{forbidden}''. Indeed, if $\Gamma_\vp(G) = \vt$, then edges in $\sF(\vt)$ are forbidden in the sense they cannot be present in $G$, while the only other edges that $G$ can contain in addition to those in $\vt$ belong to the collection of permitted edges $\cP(\vt)$. When $G \in \bG(\vt)$, by construction,
\begin{equation}
	\label{eqn:nu-dfs}
	\nu^{\dfs} ( \vt ; G) = p_{r(\vt)} \prod_{i \in [m]} \frac{1}{d_i(\vt)!},
\end{equation}
where $r(\vt)$ denotes the root of $\vt$ and $d_i(\vt)$ denotes the number of children of $i$ in $\vt$.

Recall the law of an ordered $\vp$-tree, denoted by $\pr_{\ord}(\cdot)$, as defined in \eqref{eqn:ordered-p-tree-def}. Define the function  $L : \bT_m^{\ord} \to \bR_+$ by
\begin{equation}
\label{eqn:ltpi-def}
	\displaystyle L(\vt):= \prod_{(k,\ell)\in E(\vt)} \left[\frac{\exp(a p_k p_{\ell})- 1}{ap_k p_{\ell}} \right] \exp\left(\sum_{(k,\ell) \in \sP(\vt)} a p_k p_{\ell}\right), \qquad \vt \in \bT_m^{\ord}.
\end{equation}
We use $L(\cdot)$ to tilt the distribution of the $\vp$-tree results in the distribution
\begin{equation}
	\label{eqn:tilt-ord-dist-def}
	\tilde \pr_{\ord}( \vt) := \pr_{\ord}(\vt) \cdot \frac{L(\vt)}{\E_{\ord}[ L(\cdot)]}, \qquad \vt \in \bT_m^{\ord},
\end{equation}
where $\E_{\ord}[ L(\cdot)]$ denotes the expectation of $L(\cdot)$ with respect to the law $\pr_{\ord}$.

Now note that given a fixed planar tree $\vt \in \bT_m^{\ord}$, one can construct a connected random graph by adding each possible permitted edge $(i,j)\in \sP(\vt)$ independently with probability $q_{ij}= 1-\exp(-a p_i p_j)$. Write $\nu^{\per}(\cdot ; \vt)$ for the probability distribution of this random graph, where ``per'' stands for ``permitted edges''. Obviously by construction, the support of $\nu^{\per}(\cdot ; \vt)$ is the set $\bG(\vt)$ as defined in \eqref{eqn:def-gtpi} and has the explicit form,
\begin{equation}
	\label{eqn:prob-cg}
	\nu^{\per} ( G ; \vt) := \ind{\set{ G \in  \bG (\vt)}} \prod_{(i,j) \in \sP(\vt) \cap E(G)} q_{ij} \prod_{(i,j) \in \sP(\vt) \setminus E(G)} (1 - q_{ij}).
\end{equation}
The main result of this section is the following proposition. In words what this result says is the following: one can sample a connected random graph $\cG_m\sim \pr_{\con}$ with distribution in \eqref{eqn:rg-dist-connected}, in the following two-step procedure:
\begin{enumeratea}
	\item Generate a random planar tree $\tilde{\cT}$ using the tilted $\vp$-tree distribution $\tilde{\pr}_{\ord}(\cdot)$ given in \eqref{eqn:tilt-ord-dist-def} via the tilt $L(\cdot)$.
	\item Conditional on $\tilde{\cT}$, add each of the permitted edges $(i,j)\in\sP(\tilde{\cT})$ independently with the appropriate probability $q_{ij}$.
\end{enumeratea}

\begin{prop}
	\label{prop:distrib-gone-tildg}
	For all $G \in \bG_m^{\con}$ and $\vt \in \bT_m^{\ord}$,
	\begin{equation}
		\label{eqn:1484}
		\pr_{\con}(G) \nu^{\dfs}(\vt ; G) = \tilde \pr_{\ord}(\vt) \nu^{\per}(G;\vt).
	\end{equation}
	In particular,  $\pr_{\con}(G) = \sum_{\vt \in \bT_m^{\ord}} \tilde \pr_{\ord}(\vt) \nu^{\per}(G;\vt)$.
\end{prop}

\noindent\textbf{Proof:}. From the definition of $\nu^{\dfs}(\vt ; G)$ and $\nu^{\per}(G;\vt)$, the left hand side and right hand side of \eqref{eqn:1484} are non zero if and only if $G \in \bG(\vt)$. When $G \in \bG(\vt)$, using \eqref{eqn:rg-dist-connected} and \eqref{eqn:nu-dfs} for the left hand side gives
\begin{align}
	\label{eqn:1494}
	\pr_{\con}(G) \nu^{\dfs}(\vt ; G) = \frac{1}{Z(\vp) }{\prod_{(i,j)\in E(G)} (1 - e^{-ap_ip_j}) \prod_{(i,j)\notin E(G)} e^{-ap_i p_j}} \times p_{r(\vt)} \prod_{i \in [m]} \frac{1}{d_i(\vt)!}.
\end{align}
To ease notation write $d_i(\vt) = d_i$ for the number of children of vertex $i$ in $\vt$. Let us now simplify the right hand side of \eqref{eqn:1484}. Using \eqref{eqn:ltpi-def} and the fact that for a fixed tree $\vt$, $\sF(\vt), \sP(\vt)$ and $E(\vt)$ form a partition of all possible edges (denoted by $[m]_2$) on the vertex set $[m]$ gives

\begin{align*}
	L(\vt)
	=& \prod_{(i,j) \in E(\vt)}\left[\frac{e^{a p_i p_j}- 1}{ap_ip_j} \right] \prod_{(i,j) \in \sP(\vt)} e^{a p_i p_j}\\
	=& \frac{p_{r(\vt)}}{a^{m-1}}\prod_{i \in [m]} \frac{1}{p_i^{d_i+1}} \prod_{(i,j) \in E(\vt)}(e^{a p_i p_j}- 1) \prod_{(i,j) \in \sP(\vt)} e^{a p_i p_j}\\
	=& \frac{p_{r(\vt)}}{a^{m-1}}\prod_{i \in [m]} \frac{1}{p_i^{d_i+1}} \prod_{(i,j) \in [m]_2} e^{ap_ip_j} \prod_{(i,j) \in E(\vt)}(1- e^{-a p_i p_j}) \prod_{(i,j) \in \sF(\vt)} e^{-a p_i p_j}.
\end{align*}
Using the above display, \eqref{eqn:ordered-p-tree-def} and \eqref{eqn:prob-cg},
\begin{align*}
	\tilde \pr_{\ord}(\vt) \nu^{\per}(G;\vt)
	=& \frac{1}{\E_{\ord}[ L(\cdot)]} \prod_{i \in [m]} \frac{p_i^{d_i}}{d_i!} \times \prod_{(i,j) \in \sP(\vt) \cap E(G)} (1 - e^{-ap_ip_j}) \prod_{(i,j) \in \sP(\vt) \setminus E(G)} e^{-ap_ip_j} \times L(\vt)\\
	=& \frac{a^{-(m-1)}{\prod\nolimits_{(i,j) \in [m]_2}} e^{ap_ip_j}}{\E_{\ord}[ L(\cdot)] \prod_{i \in [m]} p_i} \times p_{r(\vt)} \prod_{i \in [m]} \frac{1}{d_i(\vt)!} \times  \prod_{(i,j)\in E(G)} (1 - e^{-ap_ip_j}) \prod_{(i,j)\notin E(G)} e^{-ap_i p_j},
\end{align*}
where the last display is obtained by using $E(G) = E(\vt) \cup (\sP(\vt) \cap E(G))$ and $E(G)^c = \sF(\vt) \cup (\sP(\vt) \setminus E(G))$. Comparing the above expression with \eqref{eqn:1494},
\begin{equation*}
	\frac{\pr_{\con}(G) \nu^{\dfs}(\vt ; G)}{\tilde \pr_{\ord}(\vt) \nu^{\per}(G;\vt)} = f(\vp, a, m),
\end{equation*}
where $f(\vp,a,m)$ is a constant, independent of $\vt$ or $G$. Since both the left and the right hand sides are probability distributions, $f(\vp,a,m) \equiv 1$. This completes the proof. \qed


\subsection{Convergence of untilted graphs}
\label{sec:conv-height-cont-connected-irg}
Using Proposition \ref{prop:distrib-gone-tildg}, define the probability distribution $\tilde \nu^{\jt}(\cdot,\cdot)$ on $\bT_m^{\ord} \times \bG_m^{\con}$ via the prescription,
\begin{equation}
	\label{eqn:nu-tilde-joint-def}
	\tilde \nu^{\jt}(\vt, G) :=  \pr_{\con}(G) \nu^{\dfs}(\vt ; G) = \tilde \pr_{\ord}(\vt) \nu^{\per}(G;\vt), \mbox{ for } \vt \in \bT_m^{\ord}, G \in \bG_m^{\con}.
\end{equation}
This is the main object of interest. Let us first study the simpler object which does not incorporate the tilt. More precisely define the probability distribution $\nu^{\jt}(\cdot, \cdot)$  on $\bT_m^{\ord} \times \bG_m^{\con}$ as follows:
\begin{equation}
	\label{eqn:nu-joint-def}
	 \nu^{\jt}(\vt, G) :=  \pr_{\ord}(\vt) \nu^{\per}(G;\vt), \mbox{ for } \vt \in \bT_m^{\ord}, G \in \bG_m^{\con}.
\end{equation}
In this section, we will study the limit behavior of $\nu^{\jt}$ and $L(\vt)$ under $\nu^{\jt}$. Write $(\cT^{\vp},\cG^{\vp})\sim \nu^{\jt}$ for the $\bT_m^{\ord} \times \bG_m^{\con}$-valued random variable with distribution $\nu^{\jt}$. The main aim of this section is the following result for the untilted object. The next section studies the tilted version.

\begin{prop}
	\label{prop:ghp-convg-non-tilt}
	
	Let $(\cT^{\vp},\cG^{\vp})$ be $\bT_m^{\ord} \times \bG_m^{\con}$-valued random variable with distribution $\nu^{\jt}$ viewed as measured metric spaces using the vertex weights $\vp$. Then under Assumptions \ref{ass:aldous-AMP} and \ref{ass:additional-connected}, as $m\to\infty$,
	\begin{equation*}
		\left(\scl\left( \sigma(\vp), 1 \right) \cG^{\vp}, L(\cT^{\vp})\right) \weakc \left( \cG( 2 \ve,\bar \gamma \ve, \cP),  \exp\left( \bar \gamma \int_0^1 \ve(s) ds \right)\right).
	\end{equation*}
\end{prop}

Before diving into the proof, we start by giving an explicit construction of $(\cT^{\vp}, \cG^{\vp})$ from $(\vX, \cP)$, where $\vX = (X_i: i \in [m])$ are \emph{i.i.d.} Uniform$[0,1]$ r.v.s and $\cP$ is a rate one Poisson point process on $\bR_+^2$, independent of $\vX$.  The construction is based on \cite{aldous2004exploration} which starts by setting up a map $\psi_{\vp}: (0,1)^m \to \bT_m^{\ord} $ as follows. 
Fix a collection of distinct points $\vx = (x_i: i \in [m])\in (0,1)^m$. Define
\begin{equation}
\label{eqn:fp-defn}
	F^{\vp}(u) := - u + \sum_{i=1}^m p_i \ind{\set{x_i\leq u}}, \qquad u\in [0,1].
\end{equation}

\begin{figure}
\centering
\begin{minipage}{.5\textwidth}
  \centering
  \includegraphics[trim=2.5cm 2.5cm 2.2cm 3cm, clip=true, angle=0, scale=0.45]{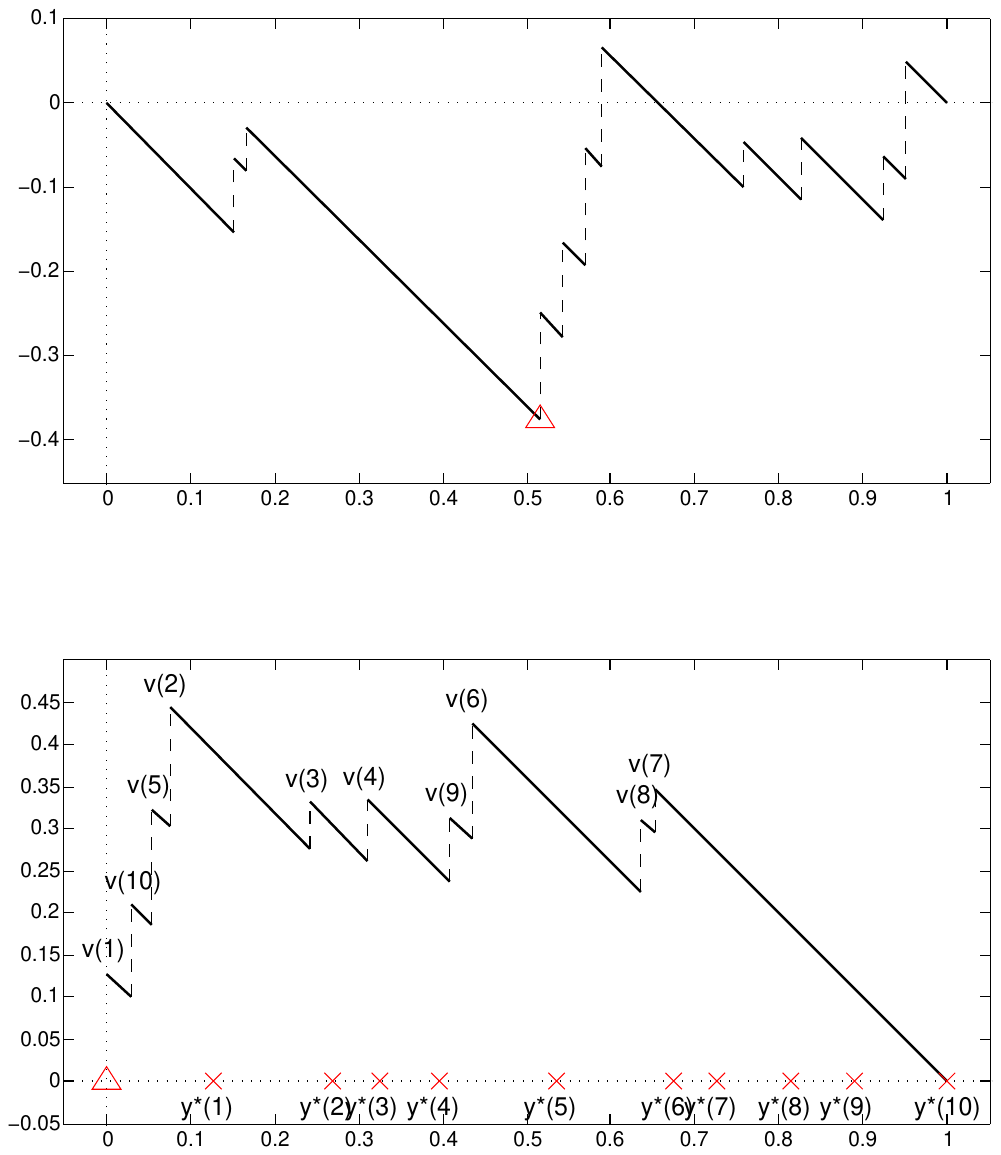}
  \captionof{figure}{The functions $F^{\vp}$ on the top and the corresponding function $F^{\exec,\vp}$ for a specific choice of $m=10$ points and a pmf $\vp$.}
  \label{fig:test1}
\end{minipage}%
\begin{minipage}{.5\textwidth}
  \centering
  \includegraphics[trim=2.5cm 0cm 1cm 0cm, clip=true, angle=0, scale=0.4]{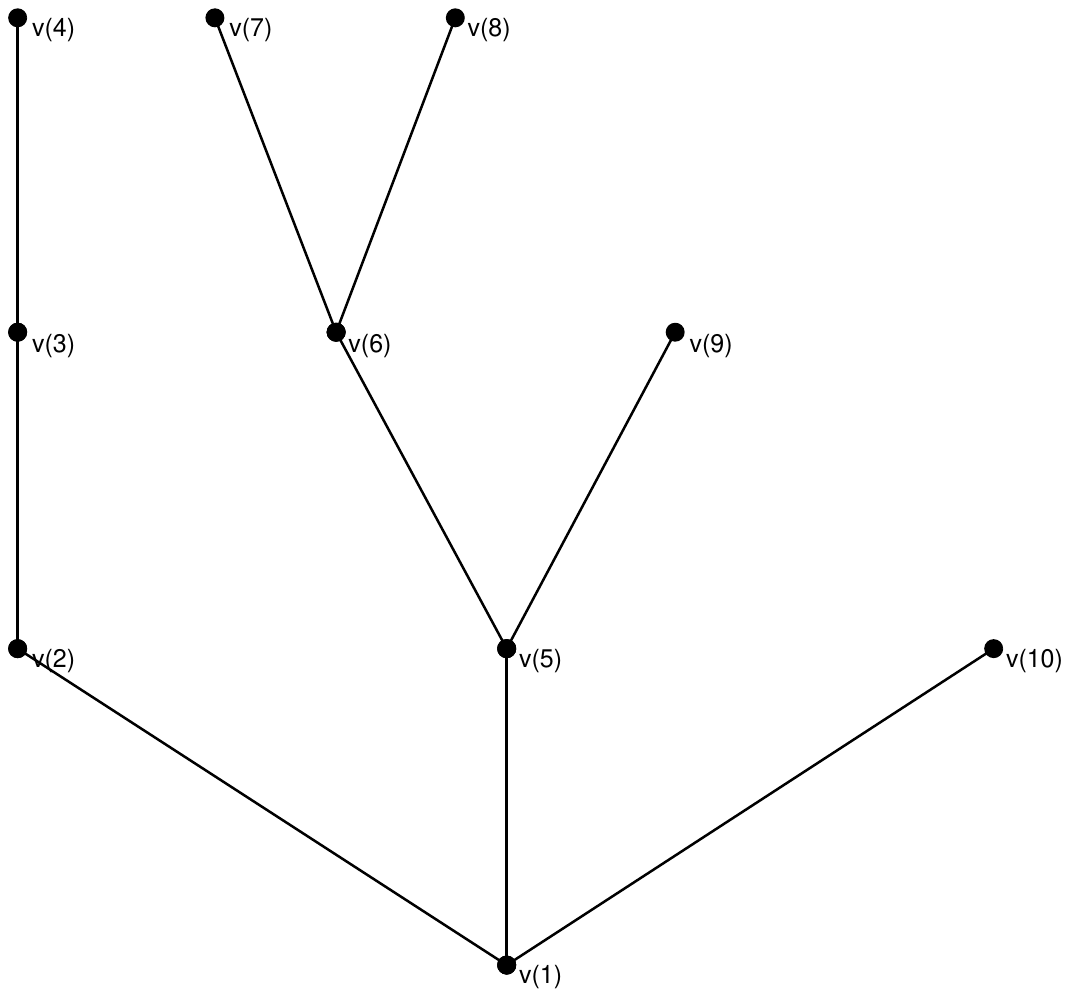}
  \captionof{figure}{The tree obtained from the depth-first construction}
  \label{fig:test2}
\end{minipage}
\end{figure}

Assume that there exists a unique point $v^* \in [m] $ such that $F^{\vp}(x_{v^*}-) = \min_{u\in [0,1]} F^{\vp}(u)$. Set $v^*$ to be the root of the tree $\psi_{\vp}(\vx)$. Define $y_i := x_i - x_{v^*}$ for $i \in [m]$, and
\begin{equation*}
\label{eqn:fexc-def}
	F^{\exec,\vp}(u):= F^{\vp}( x_{v^*} + u \mbox{ mod } 1) - F^{\vp}(x_{v*}-), \qquad 0\leq u < 1.
\end{equation*}
 Then $F^{\exec,\vp}(1-) = 0$ and $F^{\exec,\vp}(u) > 0$ for $u \in [0,1)$. Extend the definition of $F^{\exec,\vp}$ to $u \in [0,1]$ by defining $F^{\exec,\vp}(1) = 0$. We will use $F^{\exec,\vp}$ to construct a depth-first-search of an ordered tree whose exploration in this depth first manner is encoded by the function $F^{\exec,\vp}$. This in turn defines the tree $\psi_{\vp}(\vx)$. As before, in this construction we will carry along a set of explored vertices $\cO(i)$, active vertices $\cA(i)$ and unexplored vertices $\cU(i) = [m]\setminus (\cA(i)\cup \cO(i))$, for $0\leq i \leq m$. In particular we will view $\cA(i)$ as the state of a vertical stack $\cA$ after the $i$-th step in the depth-first-search.

Initialize with $\cO(0) = \emptyset$, $\cA(0) = \set{v^*}$, $\cU(0) = [m] \setminus \set{v(1)}$, and define $y^*(0) = 0$. At step $i \in [m]$, let $v(i)$ be the value that is on the top of the stack $\cA(i-1)$ and define $y^*(i) := y^*(i-1)+p_{v(i)}$. Define $\cD(i) := \set{ i \in [m] : y^*(i-1) < y_i < y^*(i) }$. Suppose $\cD(i) = \set{ u(j) : 1\leq j \leq k}$ where we have ordered these vertices in the sequence that they are found in this interval namely
\[y^*(i-1) < y_{u(1)} <... < y_{u(k)} < y^*(i).\]
Update the stack $\cA(i-1)$ as follows:
\begin{enumeratei}
	\item Delete $v(i)$ from $\cA(i-1)$.
	\item Push $u(j)$, $1\leq j\leq k$, to the top of $\cA(i-1)$ sequentially (so that $u(k)$ will be on the {\bf top} of the stack at the end).
\end{enumeratei}
Let $\cA(i)$ be the state of the stack after the above operations. Our approach here is not exactly the same as the one in \cite{aldous2004exploration}, where the vertices are pushed to the stack in the reverse order. However as remarked in \cite{aldous2004exploration} this does not effect the resulting distribution of the tree. 
Update $\cO(i) := \cO(i-1) \cup \set{v(i)}$ and $\cU(i) := \cU(i-1)\setminus \cD(i) $.

The tree $\psi_{\vp}(\vx) \in \bT_m^{\ord}$ is constructed by putting the edges $\set{ (v(i),u): i \in [m], u \in \cD(i)}$ and using the order prescribed in the above exploration to make the tree an ordered tree. The fact that this procedure actually produces a tree is proved in \cite{aldous2004exploration}.  So far we have given the construction of a deterministic tree $\psi_{\vp}(\vx)$ using $\vx \in (0,1)^m$. Using the collection of uniform random variables $\vX$ then results in $\psi_{\vp}(\vX)$ being a random ordered tree in $\bT_m^{\ord}$.
It is further shown in \cite{aldous2004exploration} that $\psi_{\vp}(\vX)$ has the same distribution as an ordered $\vp$-tree, i.e., $\psi_{\vp}(\vX)$ has the law $\pr_{\ord}$ in \eqref{eqn:ordered-p-tree-def}.

We use the same notation to denote the various constructs in the above construction when replacing $\vx$ with $\vX$, so that notation such as $\cA(i)$, $\cD(i)$ and $y^*(i)$ now correspond to random objects. Define
\begin{equation}
	\label{eqn:hp-height-def}
	H^{\vp}(u) := \mbox{ height of $v(i)$ in $\psi_{\vp}(\vX)$}, \qquad u\in (y^*(i-1), y^*(i)], i \in [m].
\end{equation}
 Extend $H^{\vp}(u)$ to $u=0$ continuously. $F^{\exec,\vp}$ in \eqref{eqn:fexc-def} and $H^{\vp}$ are random elements in $D([0,1],\bR)$.

The following theorem was proved \cite[Proposition 3 and Theorem 3]{aldous2004exploration} under a set of assumptions stronger than Assumption \ref{ass:aldous-AMP}. As remarked in \cite{aldous2004exploration}, their assumptions can be relaxed to just moment assumptions. See Remark \ref{rem:amp} below.  It turns out that Assumption \ref{ass:aldous-AMP} is a sufficient condition for the same result.
\begin{thm}
\label{thm:AMP}
Under Assumptions \ref{ass:aldous-AMP}, as $m\to\infty$,
\begin{equation}
\label{eqn:convg-fexec-hp}
	\left(\frac{F^{\exec,\vp}(\cdot)}{\sigma(\vp)},{\sigma(\vp)}H^{\vp}(\cdot)\right) \weakc (\ve, 2\ve),
\end{equation}
where $\ve$ is a standard Brownian excursion.
\end{thm}

Lemma \ref{lem:AMP} stated below is an ingredient in the proof of Theorem \ref{thm:AMP}. The proof of this lemma is similar to that of \cite[Proposition 4]{aldous2004exploration}. We relax the assumption about exponential moments used in \cite[Proposition 4]{aldous2004exploration} to the bound on $p_{\max}$ as in Assumption \ref{ass:aldous-AMP}, and the price is a stronger assumption on $p_{\min}$. The proof of Lemma \ref{lem:AMP} is outlined in Appendix \ref{Appendix} briefly.

\begin{lem}\label{lem:AMP}
Under Assumption \ref{ass:aldous-AMP}, as $m \to \infty$,
\begin{align*}
\sup_{u\in[0,1]}\left|\frac{1}{2}\sigma(\vp)H^{\vp}(u)
-\frac{1}{\sigma(\vp)}F^{\exec, \vp}(u)\right|\probc 0.
\end{align*}
\end{lem}
\noindent\textbf{Proof of Theorem \ref{thm:AMP}: }
By \cite[Equation (19)]{aldous2004exploration}, under the assumptions $\lim_{m \to \infty} \sigma(\vp) = 0$ and $ \lim_{m \to \infty} {p_{\max}}/{\sigma(\vp)} = 0$, as $m \to \infty$,
\begin{equation}
	\label{eqn:1494-1010}
	\frac{1}{\sigma(\vp)} F^{\exec, \vp}(\cdot) \weakc \ve(\cdot).
\end{equation}
The proof of Theorem \ref{thm:AMP} is completed by combining \eqref{eqn:1494-1010} and Lemma \ref{lem:AMP}. \qed

\begin{rem}\label{rem:amp}
The proof of \cite[Proposition 4]{aldous2004exploration} uses large deviation inequalities; this is where the assumption on exponential moments is used. The use of large deviation inequalities makes the proof simpler. However, as observed by the authors of \cite{aldous2004exploration} (see the remark after the statement of \cite[Theorem 3]{aldous2004exploration}), it is possible to prove this result simply by using Markov's inequality and the Burkholder-Davis-Gundy inequality (instead of large deviation bounds) and it turns out that Assumption \ref{ass:aldous-AMP} is sufficient for this purpose. We will use similar techniques in the proof of Theorem \ref{thm:ht-p-tree}, so to avoid repetition, we will only provide an outline of a proof of Lemma \ref{lem:AMP} in Appendix \ref{Appendix}.
\end{rem}
\begin{rem}\label{rem:amp-fix}
During a conversation with Gr\'egory Miermont, we discovered a mistake in the proof of \cite[Proposition 4]{aldous2004exploration}. Indeed, the claim in \cite[Lemma 11]{aldous2004exploration} that $X^\ast$ has the same distribution as $X$ defined in \cite[Lemma 10]{aldous2004exploration} is incorrect because of the following reason: The distribution of the ordered $\vp$-tree obtained from $F^{\exec,\vp}$ via the depth-first construction is given by \eqref{eqn:ordered-p-tree-def}. This is equivalent to starting from an unordered $\vp$-tree having distribution \eqref{eqn:p-tree-def}, and then ordering the children of each vertex using a uniform permutation. Thus, if $y$ is a child of some vertex on the ancestral line of another vertex $v\in[m]$, then the relative position of $y$ with respect to the ancestral line of $v$ cannot be decided using \emph{i.i.d.} Bernoulli random variables.

However, this problem can be fixed by imitating the proof of Proposition \ref{prop:qv-weird-bound}. Thus, there is an easy fix for the mistake in the proof of \cite[Proposition 4]{aldous2004exploration} under Assumption \ref{ass:aldous-AMP}.
\end{rem}

Next, we will construct a random graph $\psi_{\vp}^{G}(\vX) \in \bG_m^{\con}$ such that $(\psi_{\vp}^{G}(\vX), \psi_{\vp}(\vX)) \stackrel{d}{=} (\cG^{\vp}, \cT^{\vp} )$ as defined in Theorem \ref{prop:ghp-convg-non-tilt}. For $i \in [m]$, let $\cS(i) := \cA(i-1) \setminus \set{v(i)}$. Define the function $A_m(\cdot)$ on $[0,1]$ via
\begin{equation}
\label{eqn:ant-def}
	A_m(u):= \sum_{v\in \cS(i)}  p_v, \qquad \mbox{ for } u \in (y^*(i-1), y^*(i)], i \in [m].
\end{equation}
Define $\bar A_m(u) := a A_m (u)$, $u \in [0,1]$, where $a$ is the scaling constant in the definition of the edge probabilities $q_{ij}$. Recall that $\cP$ is a rate one Poisson point process on $\bR_+^2$, independent of $\vX$.
Let $\bar A_m \cap \cP := \set{(x,y) \in \cP: 0 \leq x \leq 1,  y \leq \bar A_m(x) }$. For each point $(x,y) \in \bar A_m \cap \cP$, define
\begin{equation}
\label{eqn:rn-x-y-def}
	r_m(x,y) = \inf\set{x' \geq x: \bar A_m(x')  < y}.
\end{equation}

Conditioned on $\psi_{\vp}(\vX)$, the graph $\psi_{\vp}^{G}(\vX)$ is constructed as follows: Suppose $\bar A_m \cap \cP  = \set{(x_l,y_l) : l \in [k] }$. Then for each $l \in [k]$ define $i_l \in [m]$ to be such that $y^*(i_l-1) < x_l < y^*(i_l) $, and define $j_l \in [m]$ to be such that $y^*(j_l) = r_m(x_l,y_l)$. Let $\psi_{\vp}^{G}(\vX)$ be the graph obtained by adding edges $(v(i_l),v(j_l))$, $l \in [k]$, to $\psi_{\vp}(\vX)$. There is a small probability that multiple edges are placed between two vertices if there are multiple points in $\cP$ that are very close to each other. In that case, let $\psi_{\vp}^{G}(\vX)$ be the simple graph obtained by replacing all multi-edges with simple edges.

The key observation is that for every edge in $\sP(\psi_{\vp}(\vX))$ of the form,
\[(v(i),v(j)) \in \sP(\psi_{\vp}(\vX)) \text{ such that } v(j) \in \cA(i-1)\setminus \set{v(i)}, \]
we can find a unique corresponding rectangle in $\bR_+^2$ below the path $\bar A(\cdot)$:
\begin{equation*}
	R(i,j) := \set{(x,y) \in \bR_+^2: y^*(i-1) \leq x < y^*(i), \; \bar A_m(y^*(j)) < y \leq \bar A_m(y^*(j)-)}.
\end{equation*}
Note that these rectangles have the following properties:
\begin{enumeratea}
	\item They consist of a partition of $\set{(x,y) \in \bR_+^2: 0 \leq x < 1, \; 0 < y \leq \bar A_m(x)}$.
	\item $R(i,j)$ has width $p_{v(i)}$ and height $a p_{v(j)}$.
	\item $(v(i),v(j))$ is an edge in $\psi_{\vp}^{G}(\vX)$ if and only if $R(i,j) \cap \cP \neq \emptyset$.
\end{enumeratea}
Based on the above observation, since $\cP$ is a Poisson point process, for $(v(i),v(j)) \in \sP(\psi_{\vp}(\vX))$,
\begin{align}
\pr(\mbox{$(v(i),v(j))$ is added to }  \psi_{\vp}(\vX) |  \psi_{\vp}(\vX) )
&= 1-\exp\left(-a p_{v(i)} p_{v(j)}\right). \label{eqn:prob-add-edges}
\end{align}
Further $\pr( \psi_{\vp}^{G}(\vX) = G |\psi_{\vp}(\vX) = \vt) = \nu^{\per}(G; \vt)$ and thus
\begin{equation}
	\label{eqn:1574}
	(\psi_{\vp}^{G}(\vX), \psi_{\vp}(\vX)) \stackrel{d}{=} (\cG^{\vp}, \cT^{\vp}).
\end{equation}

\noindent {\bf Proof of Proposition \ref{prop:ghp-convg-non-tilt}:} Using \eqref{eqn:convg-fexec-hp} and the Skorohod embedding, we can construct $\set{F^{\exec,\vp}, H^{\vp}: m \in \bN}$ on a common probability space $\Omega_1$ such that
\begin{equation*}
	\left(\frac{F^{\exec,\vp}(\cdot)}{\sigma(\vp)},{\sigma(\vp)}H^{\vp}(\cdot)\right) \convas (\ve, 2\ve),
\end{equation*}
where the convergence is with respect to the product of the Skorohod topology on $D([0,1],\bR) \times D([0,1],\bR)$.
Let $\cP$ be a rate one Poisson point process on $\bR_+^2$, independent of $\set{F^{\exec,\vp}, H^{\vp}: m \geq 1}$ and the almost sure limit $\ve$.  
 By \eqref{eqn:1574}, we can write
\begin{equation*}
	 (\cG^{\vp}, \cT^{\vp}) := (\psi_{\vp}^{G}(\vX), \psi_{\vp}(\vX)).
\end{equation*}
We start with a preliminary lemma analyzing asymptotics for $A_m(\cdot)$ in \eqref{eqn:ant-def}.

\begin{lem}
	\label{lem:fexc-ant}
	As $m\to\infty$,
	\[\sup_{t\in [0,1]} \left|\frac{F^{\exec,\vp}(t) - A_m(t)}{\sigma(\vp)}\right|\convas 0.\]
\end{lem}
\noindent {\bf Proof:} By the definition of $F^{\exec,\vp}$,
\begin{equation}
	\label{eqn:sum-of-wts-active}	
	F^{\exec,\vp}(y^*(i)) = \sum_{v \in \cA(i) } p_v, \;\; \mbox{ for } i \in [m].
\end{equation}
Recall that $\cS(i) = \cA(i-1)\setminus \set{v(i)}$.  By \eqref{eqn:ant-def},
\begin{equation}
	\label{eqn:1534}
	A_m(t) = \sum_{v \in \cS(i)} p_v = \sum_{v \in \cA(i-1)}p_v - p_{v(i)}, \;\; \mbox{ for } t \in (y^*(i-1), y^*(i)].
\end{equation}
Thus
\begin{align*}
	\sup_{t \in (y^*(i-1), y^*(i)]}|A_m(t)-F^{\exec,\vp}(t)| \leq& | A_m(y^*(i)) - F^{\exec, \vp}(y^*(i-1)) | \\
	& \qquad + 	\sup_{t \in (y^*(i-1), y^*(i)]}|F^{\exec,\vp}(t) - F^{\exec, \vp}(y^*(i-1))|\\
	=& p_{v(i)} + 	\sup_{t \in (y^*(i-1), y^*(i)]}|F^{\exec,\vp}(t) - F^{\exec, \vp}(y^*(i-1))|
\end{align*}
Denoting $\Delta_m(\delta) := \sup_{0 \leq s < t \leq 1, |s-t| \leq \delta} |F^{\exec,\vp}(s)- F^{\exec,\vp}(t)|$, then
\begin{equation*}
	\sup_{t\in [0,1]} \left|\frac{F^{\exec,\vp}(t) - A_m(t)}{\sigma(\vp)}\right| \leq \frac{p_{\max}}{\sigma(\vp)} + \frac{\Delta_m(p_{\max})}{\sigma(\vp)}.
\end{equation*}
By Assumption \ref{ass:additional-connected}, $p_{\max}/\sigma(\vp) \to 0$ and $p_{\max} \to 0$ as $m \to \infty$. In addition, since $\sup_{t \in [0,1]}|F^{\exec,\vp}(t)/\sigma(\vp) - \ve(t)| \to 0$ and $\ve(\cdot)$ is continuous on $[0,1]$, we have ${\Delta_m(p_{\max})}/{\sigma(\vp)} \to 0$ as $m \to \infty$ as well. The proof of Lemma \ref{lem:fexc-ant} is completed. \qed
The proof of the above lemma also implies the following result, which we state here for later use:
\begin{lem} \label{lem:1555}
	For all $m \geq 1$,
	$$\|A_m\|_\infty \leq \|F^{\exec,\vp}\|_\infty.$$
\end{lem}
\noindent\textbf{Proof: } Note that $A_m(\cdot)$ is piecewise constant and $F^{\exec,\vp}(\cdot)$ is piecewise linear. By \eqref{eqn:sum-of-wts-active} and \eqref{eqn:1534},
\begin{equation*}
	\|A_m\|_\infty = \sup_{i \in [m]}\left(\sum_{v \in \cA(i-1)}p_v - p_{v(i)}\right) \leq \sup_{i \in [m]}\left(\sum_{v \in \cA(i-1)}p_v \right) = \|F^{\exec,\vp}\|_\infty.
\end{equation*}
This completes the proof of Lemma \ref{lem:1555}. \qed\\

By Lemma \ref{lem:fexc-ant} and the construction of the point process $\cP$, since $a \sigma(\vp) \to \bar \gamma$,
\begin{equation}
	\label{eqn:1555}
	\left(\frac{1}{\sigma(\vp)} F^{\exec,\vp} ,{\sigma(\vp)} H^{\vp}, \bar A_m, \cP \right) \convas (\ve, 2\ve, \bar\gamma \ve, \cP).
\end{equation}
Note that ${\sigma(\vp)} H^{\vp}$ encodes the information of distances with regards to the underlying $\vp$-tree, while $(\bar A_m, \cP)$ encodes the information of the shortcuts or surplus edges added to the tree. We will combine these two pieces of information to prove the convergence of the untilted graphs in Lemma \ref{lem:1672}. In addition, we will use the convergence of $\bar A_m$ again in Lemma \ref{lem:wts-permitted} to prove the convergence of $L(\cT^{\vp})$. Now let us continue with the proof of Proposition \ref{prop:ghp-convg-non-tilt}.

By \eqref{eqn:1555}, there exists $k \in \bN_0$ such that for all $m$ large enough
\begin{equation*}
	\bar A_m \cap \cP = \set{(x_l,y_l) : l =1,2,...,k}.
\end{equation*}
Recall from Section \ref{sec:real-tree-shortcut-limit-object} that given any excursion $h$ we can construct a real tree $\cT(h)$. Let $(v(i_l),v(j_l))$ be as defined below \eqref{eqn:rn-x-y-def}, $r(x_l,y_l)$ be as defined in \eqref{eqn:r-x-y-def} by replacing $g$ with $\bar \gamma \ve$, and $q_{2\ve}$ be the canonical map $[0,1] \to \cT(2\ve)$.
Then $\cG^{\vp}$ and $\cG(2\ve,\bar \gamma \ve, \cP)$ are obtained from identifying the pairs $(v(i_l),v(j_l))$ and $(q_{2\ve}(x_l), q_{2\ve}(r(x_l,y_l)))$ respectively, for $1\leq l \leq k$. Denote
\begin{equation*}
	\cG_m^{\vp} :=  \scl\left( \sigma(\vp), 1  \right) \cdot \cG^{\vp} \;\;\mbox{ and } \;\;\cT_m^{\vp} :=  \scl\left( \sigma(\vp), 1  \right) \cdot \cT^{\vp}.
\end{equation*}
In order to complete the proof of Proposition \ref{prop:ghp-convg-non-tilt}, we will rely on the following two lemmas:
\begin{lem}
	\label{lem:1672}
	$\cG_m^{\vp} \convas  \cG( 2 \ve,\bar \gamma \ve, \cP)$, as $ m \to \infty$.
\end{lem}
\textbf{Proof:} By \cite[Lemma 4.2]{addario2013scaling}, we need to construct, for each $m \in \bN$,  a correspondence $C_m$ between $\cT_m^{\vp}$ and $\cT(2 \ve)$ and a measure $\xi_m$ on the space $\cT_m^{\vp} \times \cT(2\ve)$ such that
\begin{enumeratei}
	\item $(v(i_l), q_{2\ve}(x_l)) \in C_m$ and $(v(j_l), q_{2\ve}(r(x_l,y_l))) \in C_m$, for $l =1,2,...,k$.
	\item $\xi_m(C_m^c) \convas 0$ as $m \to \infty$.
	\item $D(\xi_m) \convas 0$ as $m \to \infty$, where $D(\xi_m)$ is the discrepancy defined in \eqref{eqn:def-discrepancy}.
	\item $\dis(C_m) \convas 0$ as $m \to \infty$, where $\dis(C_m)$ is the distortion defined in \eqref{eqn:def-distortion}.
\end{enumeratei}
Note that here $k$, $C_m$ and $\xi_m$ are all random objects. Once the above conditions are verified, by \cite[Lemma 4.2]{addario2013scaling},
\begin{align*}
	d_{\GHP}(\cG_m^{\vp}, \cG(2\ve, \bar \gamma \ve, \cP) ) \leq (k+1) \max\left( \frac{1}{2}\dis(C_m), D(\xi_m), \xi_m(C_m^c) \right) \convas 0,
\end{align*}
as $m \to \infty$ and therefore Lemma \ref{lem:1672} is proved.

Now we describe the construction of $C_m$ and $\xi_m$. Define
\begin{equation}
	\epsilon_m := 2 \sup_{l=1,2,...,k}|r_m(x_l,y_l) - r(x_l,y_l)|.
\end{equation}
By definition of $r(x,y)$,
\begin{equation*}
	\ve(x) > \ve(r(x_l,y_l)) \mbox{ for } x \in [x_l, r(x_l, y_l)), l=1,2,...,k.
\end{equation*}
Further by the property of Brownian excursions, for each $\delta > 0$, there exists $x \in [r(x_l,x_l),r(x_l,x_l) + \delta)$ such that $\ve(x) < \ve(r(x_l,y_l))$. Since $\sup_{t\in [0,1]}|\bar A_m(t) - \bar \gamma \ve(t)| \convas 0$, then
\begin{equation*}
	|r_m(x_l,y_l) - r(x_l,y_l)| \convas 0 \mbox{ as } m \to \infty, \mbox{ for } l =1,2,...,k.
\end{equation*}
Thus $\epsilon_m \convas 0$ as $m \to \infty$.

Define the correspondence $C_m$ as
\begin{equation*}
	C_m := \set{ (v(i), q_{2\ve}(x)): i \in [m], x\in [ 0 \vee (y^*(i - 1) - \epsilon_m), 1 \wedge (y^*(i) + \epsilon_m)] }.
\end{equation*}
By the definition of $\epsilon_m$, the condition (i) is automatically satisfied. Define the measure $\xi_m$ as
\begin{equation}
	\xi_m( \set{v(i)} \times A) := \mbox{Leb}\left(q_{2\ve}^{-1}(A) \cap [y^*( i-1), y^*( i)]\right),
\end{equation}
for $i \in [m]$, $A \subset [0,1]$ measurable. Since the map $i \mapsto v(i)$ is 1-1, $C_m$ and $\xi_m$ above are well defined. It is easy to check that $\xi_m(C_m) \equiv 1$ and $D(\xi_m) \equiv 0$, thus the conditions (ii) and (iii) are also satisfied. We only need to check the condition (iv). Let $(v(i_1), q_{2\ve}(u_1))$ and $(v(i_2), q_{2\ve}(u_2))$ be two elements in $C_m$. Denote $d_1$ and $d_2$ for the metric on $\cT_m^{\vp}$ and $\cT(2\ve)$ respectively. Observe that if either one is an ancestor of the other, we have
$$ {d_1(v(i_1),v(i_2)) }/{\sigma(\vp)} = H^{\vp}(y^*( i_1)) + H^{\vp}(y^*( i_2)) - 2\inf_{t \in [y^*( i_1), y^*( i_2)]} H^{\vp}(t), $$
 otherwise
$$ {d_1(v(i_1),v(i_2)) }/{\sigma(\vp)} = H^{\vp}(y^*( i_1)) + H^{\vp}(y^*( i_2)) - 2\inf_{t \in [y^*( i_1), y^*( i_2)]} H^{\vp}(t) + 2.$$
Thus
\begin{align*}
	&|d_1(v(i_1), v(i_2)) - d_2(q_{2\ve}(u_1), q_{2\ve}(u_2))| \\
	\leq& \left| \sigma(\vp)H^{\vp}(y^*(\bar i_1)) + \sigma(\vp)H^{\vp}(y^*(\bar i_2)) - 2\sigma(\vp)\inf_{t \in [y^*(\bar i_1), y^*(\bar i_2)]} H^{\vp}(t) \right. \\
	& \left. - 2 \ve(u_1)- 2\ve(u_2) + 4 \inf_{t \in [u_1,u_2]}\ve(t) \right| + 2 \sigma(\vp)\\
	\leq& 4 \sup_{t \in [0,1]}\left| \sigma(\vp)H^{\vp}(t) - 2 \ve(t) \right| + 8\Delta_{\ve}(\epsilon_m) + 2 \sigma(\vp),
\end{align*}
where $\Delta_{\ve}(\delta) = \sup_{0 \leq s < t \leq 1, |s-t| < \delta} | \ve(s)- \ve(t)|$, for $\delta > 0$. Since the above bound holds for all $(v(i_1), q_{2\ve}(u_1))$ and $(v(i_2), q_{2\ve}(u_2))$ in the correspondence $C_m$, therefore it is also an upper bound on the distortion $\dis(C_m)$ (see the defintion in \eqref{eqn:def-distortion}). In addition, the above expression convergence to zero almost surely as $m \to \infty$. Thus Condition (iv) is verified.
The proof of Lemma \ref{lem:1672} is completed.\qed

The last lemma that we need to complete the proof of Proposition \ref{prop:ghp-convg-non-tilt} is the following:
\begin{lem}\label{lem:wts-permitted}
 As $m \to \infty$,
\begin{equation*}
	L(\cT^{\vp}) = \left[\prod_{(i,j)\in E(\cT^{\vp})} \frac{\exp(a p_i p_j)- 1}{ap_ip_j} \right] \exp\left(\sum_{(i,j) \in \sP(\cT^{\vp})} a p_i p_j\right) \convas \exp\left(\bar \gamma \int_0^1 \ve(s)ds  \right).
\end{equation*}
\end{lem}
\noindent\textbf{Proof:} By the basic inequalities $1 \leq (e^x-1)/x \leq e^x$ for $x >0$, we have for $\vt \in \bT_m^{\ord}$,
\begin{equation*}
	 1 \leq \prod_{(i,j)\in E(\vt)} \frac{\exp(a p_i p_j)- 1}{ap_ip_j} \leq \exp\left( a\sum_{(i,j)\in E(\vt)} p_i p_j  \right) \leq \exp( a p_{\max}),
\end{equation*}
where the last inequality follows using the fact that $\vt$ is a tree, thus for each $(i,j) \in E(\vt)$ such that $i$ is the parent of $j$ we have $p_ip_j \leq p_{\max} p_j$; further by definition of $\vp$ we have $\sum_{j\in \vt} p_j\leq 1$. By Assumption \ref{ass:additional-connected}, we have $ap_{\max} \to 0$, thus the above display goes to one as $m \to \infty$. Then notice that
\begin{align*}
	\sum_{(i,j) \in \sP(\cT^{\vp})} a p_i p_j = a \sum_{i \in [m]} \sum_{ j \in \cS(i)} p_i p_j= \int_0^1 \bar A_m(s) ds \to \bar \gamma \int_0^1 \ve(s)ds,
\end{align*}
as $m \to \infty$, where the last convergence follows since $\bar A_m \convas \bar \gamma \ve$. The proof of Lemma \ref{lem:wts-permitted} is thus completed. \qed\\

\noindent\textbf{Completing the proof of Proposition \ref{prop:ghp-convg-non-tilt}:} The proof follows from Lemmas \ref{lem:1672} and \ref{lem:wts-permitted}. \qed

\subsection{Uniform integrability of the tilt}\label{sec:connected-tightness-tilt}

The final ingredient needed in proving Theorem \ref{thm:conv-condition-on-connected} is the tightness of $L(\cT^{\vp})$ where $L(\cdot)$ is the tilt as in \eqref{eqn:ltpi-def} and $\cT^{\vp}$ is a random $\vp$-tree with distribution $\pr_{\ord}$ as in \eqref{eqn:ordered-p-tree-def}. We start with a concentration inequality on $\|F^{\exec,\vp}\|_\infty$ that allows us to control the tilt appearing on the right hand side of \eqref{eqn:ltpi-def}. A key step is a concentration inequality for partial sums when sampling without replacement, a problem studied in a slightly different setting in \cite{serfling1974probability}.

\begin{lem}
	\label{lem:fexp-linfty-tail}
	 Recall that $\sigma(\vp) = \sqrt{\sum_{i=1}^m p_i^2}$ and $p_{\max} = \max_{i \in [m]} p_i$. Assume that
	\begin{equation}
		\label{eqn:1743}
		 4 p_{\max} \leq x \leq \frac{16\sigma^2(\vp)}{p_{\max}}.
	\end{equation}
	Then
	\begin{equation*}
		\pr ( \|F^{\exec, \vp}\|_\infty > x ) \leq 12 \exp\left( - \frac{x^2}{1024(\sigma(\vp))^2}\right).
	\end{equation*}
\end{lem}

%
%
\noindent \textbf{Proof:} Write $X_1, \ldots X_m$ for the \emph{i.i.d.} Uniform$[0,1]$ random variables used to construct $F^{\vp}$ which is then used to construct $F^{\exec,\vp}$ from \eqref{eqn:fexc-def}. Let $X_{\sss(1)}< X_{\sss(2)}< \cdots X_{\sss(m)}$ be the corresponding order statistics and let $\pi$ denote the corresponding permutation of $[m]$ namely $X_{\sss(i)} = X_{\pi(i)}$. Obviously $\pi$ is a uniform random permutation. Now by definition
\begin{equation}
\label{eqn:fp-exec-max-min-bd}
	\|F^{\exec,\vp}\|_{\infty} \leq \sup_{t \in [0,1]} F^{\vp}(t) + \left|\inf_{t \in [0,1]}F^{\vp}(t)\right|.
\end{equation}
Let us analyze the first term. Define $\vartheta_i:=-X_{\sss(i)}+ \sum_{j=1}^i p_{\pi(j)} $, namely the value $F^{\vp}(\cdot)$ at each location with a positive jump. Since $\sup_{t\in [0,1]} F^{\vp}(t)  \leq \sup_{i \in [m]}|\vartheta_i|$, we consider
\begin{align}
	\pr\left( \sup_{i \in [m]}|\vartheta_i| \geq \frac{x}{2}\right)
	 \leq& \pr\left(\sup_{i \in [m]}\left|-X_{\sss(i)}+\frac{i}{m}\right|\geq \frac{x}{4}\right)+ \pr\left(\sup_{i \in [m]}\left|\sum_{j=1}^i p_{\pi(j)} - \frac{i}{m}\right|\geq \frac{x}{4}\right)\notag\\
	:=& T_1+ T_2 \label{eqn:t1-t2}
\end{align}
Let $F_m(u) := n^{-1}\sum_{i=1}^m \ind{\set{X_i \leq u}}$, $u \in [0,1]$, denote the empirical distribution function of $(X_i: 1\leq i\leq n)$ so that $F_m(X_{\sss(i)}) = i/m$. Thus by the Dvoretzky–Kiefer–Wolfowitz (DKW) inequality \cite{massart1990}
\begin{equation}
\label{eqn:t1-bd}
	T_1 =\pr\left(\sup_{i \in [m]}\left|F_m(X_{\sss(i)}) -X_{\sss(i)}\right|\geq \frac{x}{4}\right) = \pr\left(\sup_{u\in [0,1]}|F_m(u) - u|\geq \frac{x}{4}\right) \leq 2\exp(-mx^2/8).
\end{equation}
We now analyze $T_2$. Since $\vp$ is a probability distribution, for any $m/2\leq k \leq m-1$, $|\sum_{j=1}^k p_{\pi(j)} - k/m| = |\sum_{j=k+1}^m p_{\pi(j)} - (m-k)/m| $. Without loss of generality, assume $m$ is even. Define
\[\displaystyle p(m,x):=\pr\left(\sup_{k \in [m/2]}\left|\sum_{j=1}^k p_{\pi(j)} - \frac{k}{m}\right|\geq \frac{x}{4}\right).\]
Now
\begin{align}
	T_2 &\leq p(m,x)+ \pr\left(\sup_{m/2\leq k\leq m-1}\left|\sum_{j=1}^k p_{\pi(j)} - \frac{k}{m}\right|\geq \frac{x}{4}\right)\notag\\
	&\leq p(m,x) + \pr\left(\sup_{m/2\leq k\leq m-1}\left|\sum_{j=k+1}^m p_{\pi(j)} - \frac{m-k}{m}\right|\geq \frac{x}{4}\right)\notag\\
	&=p(m,x) + \pr\left(\sup_{k^\prime \in [m/2]}\left|\sum_{l=1}^{k^\prime}p_{\pi(m-l+1)}-\frac{k^\prime}{m} \right|\geq \frac{x}{4}\right)=2 p(m,x), \label{eqn:t1-2pmx}
\end{align}
where the last line follows by noting that the permutation $\pi^\prime$ defined via $\pi(l) = \pi(m-l+1) $ is again a uniform permutation on $[m]$.  We are now left with bounding $p(m,x)$. Assume that we generate $\pi$ by sequentially drawing without replacement from $[m]$. For $k\geq 1$, let $\cF_k$ denote the $\sigma$-field generated by $(\pi(1),\ldots, \pi(k))$.  Writing $S_0 = 0$ and $S_k := \sum_{j=1}^k p_{\pi(j)}$, $k \in [m]$, it is easy to check that $\set{Y_k: k=0,1,...,m-1}$ defined by the following is an $\cF_k$-martingale:
\begin{equation*}
	Y_k := \frac{S_k - k/m}{m-k}, \mbox{ for } k=0, 1,..., m-1.
\end{equation*}
Note that $\sup_{i \in [m/2]}|S_i - i/m| \leq m \sup_{i \in [m/2]} |Y_i|$, thus
\begin{align}
	\label{eqn:skvn-bound}
	p(m,x) \leq \pr\left(\sup_{k \in [m/2]}|Y_k| \geq \frac{x}{4m} \right),
\end{align}
For $h > 0$, since $\exp(hx) >0$ is convex in $x$, then $\exp(h Y_k)$ is a sub-martingale. Hence,
\begin{equation}
\label{eqn:ix-bd}
	\pr\left(\sup_{k \in [m/2]}Y_k \geq \frac{x}{4m} \right) \leq \exp( - \frac{hx}{4m})\E\left[ \exp(hY_{m/2}) \right].
\end{equation}
By a similar bound on $\pr\left(\inf_{k \in [m/2]}Y_k \leq -\frac{x}{4m} \right)$ and the fact $Y_{m/2} \stackrel{d}{=} - Y_{m/2}$, following \eqref{eqn:skvn-bound},
\begin{equation}
	\label{eqn:1805}
	p(m,x) \leq 2  \exp( - \frac{hx}{4m})\E\left[ \exp(hY_{m/2}) \right] = 2  \exp( - \frac{hx}{4m}) \E \left[ \exp\left(\frac{2h}{m}S_{m/2} - \frac{h}{m}\right)\right].
\end{equation}
Now we use the standard technique of bounding the moment generating function of $S_{m/2}$ by repeatedly conditioning on the previous time steps. Note that for $0 < \delta < 1/p_{\max}$ and $k \in [m/2]$ ,
\begin{align}
	\E[ \exp(\delta p_{\pi(k+1)}) \mid \cF_k]
	=& \frac{1}{m-k} \sum_{j \notin \set{ v(i): i \in [k]}} \exp(\delta p_{j}) \nonumber \\
	\leq& \frac{1}{m-k} \sum_{j \notin \set{ v(i): i \in [k]}} ( 1 + \delta p_j + \delta^2 p_j^2) \nonumber\\
	\leq& 1 + \frac{\delta}{m-k}(1 - \sum_{j \in [k]}p_{\pi(j)}) + \frac{2\delta^2 \sigma^2(\vp) }{m} \nonumber\\
	\leq& \exp\left( \frac{\delta}{m-k}(1 - S_k) + \frac{2\delta^2 \sigma^2(\vp) }{m} \right), \label{eqn:1812}
\end{align}
where the second line uses the fact that $e^x < 1 + x + x^2$ for $x \in [0,1]$ and the third line uses the fact $\sum_{j \in [m]} p_{\pi(j)} = 1$ and $k \leq m/2$. Using \eqref{eqn:1812} repeatedly in evaluating $\E[\exp(\delta S_k)]$ for $k \leq m/2$,
\begin{align*}
	\E[\exp(\delta S_k)]
	=& \E\left[ \exp(\delta S_{k-1}) \E\left[ \exp(\delta p_{\pi(k)}) \mid \cF_{k-1} \right] \right]\\
	\leq&  \E\left[ \exp(\delta S_{k-1}) \exp\left(\frac{\delta}{m-(k-1)}(1 - S_{k-1}) + \frac{2\delta^2 \sigma^2(\vp) }{m} \right)  \right]\\
	=&  \E\left[  \exp\left( \frac{m-k}{m-k+1} \delta S_{k-1} \right) \right]  \exp\left(\frac{2\delta^2 \sigma^2(\vp) }{m}\right) \exp\left(\frac{\delta}{m-k+1}\right)\\
	\leq& \E\left[  \exp\left( \frac{m-k}{m-k+1} \cdot \frac{m-k+1}{m-k+2} \delta S_{k-2} \right) \right]  \exp\left( 2 \cdot \frac{2\delta^2 \sigma^2(\vp) }{m}\right)\\
	 &\times \exp\left( \frac{\delta(m-k)}{(m-k)(m-k+1)} + \frac{\delta(m-k)}{(m-k+1)(m-k+2)} \right).
\end{align*}
Proceeding inductively, we have
\begin{align*}
	\E[\exp(\delta S_k)]
	\leq& \E\left[ \frac{m-k}{m} \delta S_0\right] \exp\left( k \cdot \frac{2\delta^2 \sigma^2(\vp) }{m} + {(m-k)\delta} \cdot \sum_{j=0}^{k-1} \frac{1}{(m-k+j)(m-k+j+1)} \right)\\
	=&  \exp\left( k \cdot \frac{2\delta^2 \sigma^2(\vp) }{m} +  {(m-k)\delta}\cdot\frac{k}{m ( m-k)} \right).
\end{align*}
Note that in the $l$-th iteration of applying \eqref{eqn:1812}, $\delta$ is replaced by $\delta(m-k)/(m-k+l-1)$, which is less than $\delta$. Therefore, by assuming $\delta < 1/p_{\max}$, all iterative use of \eqref{eqn:1812} are valid. Taking $k=m/2$ in the above inequality, we have
\begin{equation*}
	\E[\exp(\delta S_{m/2})] \leq \exp( \delta^2 \sigma^2(\vp) + \delta/2).
\end{equation*}
Using the above bound with $\delta = 2h/m$ in \eqref{eqn:1805}, we have
\begin{equation}
	\label{eqn:1842}
	p(m,x) \leq 2 \exp\left(-\frac{hx}{4m} + \frac{4h^2\sigma^2(\vp)}{m^2} + \frac{h}{m} - \frac{h}{m}\right) \leq 2\exp\left( - \frac{x^2 }{256 \sigma^2(\vp)}\right),
\end{equation}
where the last inequality follows from taking $h=mx/32\sigma^2(\vp)$. By our choice of $\delta$ and $h$, the restriction $\delta < 1/p_{\max}$ reduces to the upper bound in the assumption \eqref{eqn:1743}.

Now combining \eqref{eqn:t1-bd}, \eqref{eqn:t1-2pmx} and \eqref{eqn:1743},
\begin{align}
 \pr\left(\sup_{t\in [0,1]} F^{\vp}(t) \geq \frac{x}{2}\right) \leq	\pr\left(\sup_{i \in [m]}|\vartheta_i| \geq \frac{x}{2}\right) \leq& 2 \exp\left(-\frac{mx^2}{8}\right) + 4\exp\left( - \frac{x^2 }{256 \sigma^2(\vp)}\right) \nonumber \\
 \leq&  6\exp\left( - \frac{x^2 }{256 \sigma^2(\vp)}\right),	\label{eqn:1894}
\end{align}
where the last bound uses the fact that,
\begin{equation}
	\label{eqn:m-sigma-p}
	m \sigma^2(\vp) = m \cdot \left(\sum_{i \in [m]} p_i ^2\right) \geq \left( \sum_{i \in [m]} p_i\right)^2 = 1.
\end{equation}
This tackles the first term in \eqref{eqn:fp-exec-max-min-bd}. To deal with the second term, define $\vartheta_i^\prime = -X_{\sss(i)}+ \sum_{j=1}^{i-1}p_{\pi(j)}$ so that for $\vartheta_i$ as defined after \eqref{eqn:fp-exec-max-min-bd}, $\vartheta_i = \vartheta_i^\prime + p_{\pi(i)}$. Then
\[\left|\inf_{t\in [0,1]}F^{\vp}(t)\right| = \sup_{i \in [m]}|\vartheta_i^\prime|\leq\sup_{i \in [m]} |\vartheta_i|+ p_{\max}. \]
By the assumption $p_{\max} < x/4$ and \eqref{eqn:1894},
\begin{align*}
\pr\left(|\inf_{t\in [0,1]} F^{\vp}(t)|\geq \frac{x}{2}\right)
&\leq \pr\left(\sup_{i \in [m]}\left|\vartheta_i\right|\geq \frac{x}{4}\right) \leq 6\exp\left( - \frac{x^2 }{1024 \sigma^2(\vp)}\right).
\end{align*}
This together with \eqref{eqn:fp-exec-max-min-bd} and \eqref{eqn:1894} completes the proof of Lemma \ref{lem:fexp-linfty-tail}.\qed

The next result uses Lemma \ref{lem:fexp-linfty-tail} to obtain bounds on the moment generating function of $\|F^{\exec,\vp} \|_\infty$.

\begin{cor}
	\label{cor:exp-moments-fexe-infty}
	For any $B > 0$ satisfying
	\begin{equation}\label{eqn:1770}
		\frac{p_{\max}}{[\sigma(\vp)]^{3/2}} \leq \sqrt{1/(8B)},
	\end{equation}
    there exists $K_{\ref{cor:exp-moments-fexe-infty}} = K_{\ref{cor:exp-moments-fexe-infty}}(B)$ such that
	\begin{equation*}
		 \E \left[ \exp\left( \frac{B \|F^{\exec,\vp}\|_\infty}{\sigma(\vp)}\right) \right] \leq K_{\ref{cor:exp-moments-fexe-infty}}.
	\end{equation*}
\end{cor}
\noindent\textbf{Proof:}  By the trivial bound $\|F^{\exec,\vp}\|_\infty \leq 1$,
\begin{align*}
	\E \left[ \exp\left( \frac{B \|F^{\exec,\vp}\|_\infty}{\sigma(\vp)}\right) \right] =&  \int_0^{\infty} B\exp(B y) \pr\left(  \frac{\|F^{\exec,\vp}\|_\infty}{\sigma(\vp)}  \geq y  \right) dy \\
	=& \int_0^{1/\sigma(\vp)} B\exp(B y) \pr\left(  \|F^{\exec,\vp}\|_\infty  \geq y \sigma(\vp) \right) dy.
	\end{align*}
	Decomposing the integral over the intervals $[0,4p_{\max}/\sigma(\vp)]$, $[4p_{\max}/\sigma(\vp),16\sigma(\vp)/p_{\max}]$ and $[16\sigma(\vp)/p_{\max},1/\sigma(\vp)]$, applying Lemma \ref{lem:fexp-linfty-tail} to the second interval gives
	\begin{align*}
	\E \left[ \exp\left( \frac{B \|F^{\exec,\vp}\|_\infty}{\sigma(\vp)}\right) \right]
	\leq& \frac{4Bp_{\max}}{\sigma(\vp)}\exp\left(\frac{4p_{\max}}{\sigma(\vp)}\right) + \int_{4p_{\max}/\sigma(\vp)}^{16\sigma(\vp)/p_{\max}} B \exp\left(B y - \frac{y^2}{1024}\right) dy \nonumber \\
	& + \pr\left(  \|F^{\exec,\vp}\|_\infty > 16 \sigma^2(\vp)/p_{\max} \right) \cdot \exp\left(\frac{B}{\sigma(\vp)}\right),  \nonumber\\
	:=& \cB_1 + \cB_2 + \cB_3. 
\end{align*}
For the first two terms above, using \eqref{eqn:1770} and $\sigma(\vp) \leq 1$,
\begin{equation*}
	\cB_1 \leq \frac{4B}{\sqrt{8B}}\exp\left( \frac{4}{\sqrt{8B}}\right), \qquad \cB_2 \leq \int_0^\infty B \exp\left( By - \frac{y^2}{1024}\right)dy.
\end{equation*}
For $\cB_3$, using \eqref{eqn:1770} and Lemma \ref{lem:fexp-linfty-tail},
\begin{align*}
	\cB_3 \leq 12 \exp\left(- \frac{\sigma^2(\vp)}{4p_{\max}^2} + \frac{B}{\sigma(\vp)}\right) \leq  12 \exp\left(- \frac{B}{\sigma(\vp)} \right) \leq 12 e^{-B}. 
\end{align*}
The proof of Corollary \ref{cor:exp-moments-fexe-infty} is completed. \qed

\begin{cor}
	\label{cor:tightness-l-uniform-integrable}
	Assume that $\gamma>0$, $B_1> 0$, and $ B_2 \in (0, 1/\sqrt{8\gamma B_1}]$ satisfy
	\begin{equation*}
		a \sigma(\vp) \leq B_1, \qquad  \frac{p_{\max}}{[\sigma(\vp)]^{3/2}} \leq B_2.
	\end{equation*}
	Let $\cT^{\vp}$ be a $\bT_m^{\ord}$-valued random variable with distribution $\pr_{\ord}$, and $L(\cdot )$ be as defined in \eqref{eqn:ltpi-def}. Then there exists a constant $K_{\ref{cor:tightness-l-uniform-integrable}}= K_{\ref{cor:tightness-l-uniform-integrable}}(\gamma,B_1,B_2)>0$ such that
	\begin{equation*}
		 \E[ L(\cT^{\vp})^\gamma ] < K_{\ref{cor:tightness-l-uniform-integrable}}.
	\end{equation*}
	In particular, when $ p_{\max}/[\sigma(\vp)]^{3/2} \to 0$ and $a \sigma(\vp) \to \bar \gamma$ as $m \to \infty$, the sequence $\set{L(\cT^{\vp}): m \geq 1}$ is uniformly integrable.
\end{cor}
\textbf{Proof:} Recall $F^{\exec, \vp}$ from \eqref{eqn:fp-defn} and $\bar A_m$ from below \eqref{eqn:ant-def}. Let $\vX =(X_i : i \in [m])$ be the \emph{i.i.d.} Uniform$[0,1]$ random variables used in the definition of $F^{\exec, \vp}$ and $\bar A_m$. Define $\cT^{\vp} = \psi_{\vp}(\vX)$ thus $\cT^{\vp}$ has the law $\pr_{\ord}$. We have
\begin{equation*}
 L(\cT^{\vp}) \leq \exp(a  p_{\max}) \exp\left( \int_0^1 \bar A_m(s) ds \right) \leq \exp(B_1B_2) \exp\left( \frac{B_1 \|F^{\exec,\vp}\|_\infty}{\sigma(\vp)}\right),
\end{equation*}
where the last inequality uses the fact $\|A_m\|_\infty \leq \|F^{\exec,\vp}\|_\infty$ (see Lemma \ref{lem:1555}).
 Then the corollary directly follows from Corollary \ref{cor:exp-moments-fexe-infty}, and we have $K_{\ref{cor:tightness-l-uniform-integrable}} = e^{\gamma B_1 B_2} K_{\ref{cor:exp-moments-fexe-infty}}(\gamma B_1)$. Taking $\gamma > 1$ we have the uniform integrability of $L(\cT^{\vp})$. \qed

Now we are ready to give the proof of Theorem \ref{thm:conv-condition-on-connected}:

\textbf{Proof of Theorem \ref{thm:conv-condition-on-connected}:} Recall that we view the connected random graph $\cG_m$ as a compact measured metric space using the graph distance as the metric and where each vertex $i\in [m]$ is assigned mass $p_i$.   To ease notation, write $\scl_m$ for the scaling operator
\begin{equation*}
	\scl_m =  \scl\left( \sigma(\vp), 1 \right),
\end{equation*}
namely we will scale the metric by $\sigma(\vp)$ leaving the measure as is.

Let $(\cT^{\vp}, \cG^{\vp})$ have the law $\nu^{\jt}$ as in \eqref{eqn:nu-joint-def}, and let $(\tilde \cT^{\vp}, \tilde \cG^{\vp})$ have the law $\tilde \nu^{\jt}$ as in \eqref{eqn:nu-tilde-joint-def}.
We want to show that for any bounded continuous function $f(\cdot)$ on $\sS$,
\begin{equation*}
	\E [f(\scl_m \cdot \tilde \cG^{\vp})] \to \E [ f( \cG( 2 \tilde \ve^{\bar \gamma}, \bar \gamma \tilde \ve^{\bar \gamma}, \cP)) ],\;\; \mbox{ as } m \to \infty.
\end{equation*}
Define $g_f(\vt)$ for $\vt \in \bT_m^{\ord}$ as
\begin{equation*}
	g_f(\vt) := \sum_{G \in \bG_m^{\con}} f( \scl_m  G) \nu^{\per}(G;\vt).
\end{equation*}
By the definition of $\nu^{\jt}$ and $\tilde \nu^{\jt}$, we have $\E[ f( \scl_m  \cG^{\vp}) \mid \cT^{\vp}] = g_f(\cT^{\vp})$ and $\E[f( \scl_m \cdot \tilde\cG^{\vp}) \mid \tilde \cT^{\vp}] = g_f(\tilde \cT^{\vp})$. Then by \eqref{eqn:tilt-ord-dist-def}, we have
\begin{equation}
	\label{eqn:1891}
	\E[f( \scl_m \cdot \tilde \cG^{\vp})] = \E[g_f(\tilde \cT^{\vp})] = \frac{\E[ g_f(\cT^{\vp})L(\cT^{\vp})] }{\E[L(\cT^{\vp})]} = \frac{\E[ f(\scl_m \cG^{\vp}) L(\cT^{\vp})] }{\E[L(\cT^{\vp})]}.
\end{equation}
By Proposition \ref{prop:ghp-convg-non-tilt} we have the joint convergence
\begin{align*}
	L(\cT^{\vp}) &\weakc \exp\left( \bar \gamma \int_0^1 \ve(s) ds \right), \\
	f(\scl_m \cG^{\vp}) L(\cT^{\vp}) &\weakc f( \cG( 2  \ve, \bar \gamma \ve, \cP)) \exp\left( \bar \gamma \int_0^1 \ve(s) ds \right).
\end{align*}
By \eqref{eqn:1891}, the above convergence and the uniform integrability of $L(\cT^{\vp})$ (Lemma \ref{cor:tightness-l-uniform-integrable}),
\begin{align*}
	\lim_{m \to \infty} \E[f( \scl_m \cdot \cG^{\vp})]
	= \frac{\E \left[ f( \cG( 2  \ve, \bar \gamma \ve, \cP)) \exp\left( \bar \gamma \int_0^1 \ve(s) ds \right) \right]}{\E \left[ \exp\left( \bar \gamma \int_0^1 \ve(s) ds \right) \right]}
	= \E \left[ f( \cG( 2  \tilde \ve^{\bar \gamma}, \bar \gamma \tilde \ve^{\bar \gamma}, \cP)) \right],
\end{align*}
where $\tilde \ve^{\bar \gamma}$ is the tilted Brownian excursion defined in \eqref{eqn:tilt-exc-def}. The proof of Theorem \ref{thm:conv-condition-on-connected} is completed. \qed

\section{Size-biased reordering and component exploration}
\label{sec:size-bias-explor}

Recall the definition of the excursion lengths $\vZ = (Z_i : i \geq 1)$ as in \eqref{eqn:z-defn}. By Theorem \ref{thm:comp-sizes-b-van},
$$ \left( \frac{|\cC_n^{\sss(i)}|}{n^{2/3}} :  i \geq 1\right) \weakc \vZ, $$
as $n \to \infty$ in the $\ldown$ topology therefore also in the product topology. The next proposition gives more asymptotic properties for the weights of vertices in each component.
\begin{prop}
	\label{prop:weight-control}
	Recall that $\cC_n^{\sss(i)}$ is the $i$-th largest component of $\cG_n^{\nr}(\vw, \lambda)$ for $i \geq 1$. Assume that the conditions in Assumption \ref{ass:wt-seq} (a) and (b) hold. Further, assume that $w_{\max} = o(n^{1/3})$. Then, for fixed $i \geq 1$,
		\begin{equation}
			\label{eqn:as-convg-comp-wts}
			 \left( \frac{|\cC_n^{\sss(i)}|}{n^{2/3}}, \frac{\sum_{v \in \cC_n^{\sss(i)}} w_v}{n^{2/3}}, \frac{\sum_{v \in \cC_n^{\sss(i)}} w_v^2}{n^{2/3}} \right)\weakc \left(Z_i, Z_i, \frac{\sigma_3 Z_i}{\sigma_1}\right)\text{ as }n\to\infty.
		\end{equation}
\end{prop}
\noindent{\bf Proof:}  We start by recalling some of the ideas in the proof of the convergence of component sizes, namely Theorem \ref{thm:comp-sizes-b-van} proved in \cite{SBVHJL10}. Recall that given a set $[n]$ and an associated weight sequence $\set{w_v: v\in [n]}$ with $w_v > 0$, a size-biased reordering is a random reordering of $[n]$ as $(v(1), v(2), \ldots, v(n))$ using the weight sequence where
 \begin{equation}
 \label{eqn:vone-choice}
 	\pr(v(1) = j )\propto w_j,\qquad j\in [n],
 \end{equation}
  and having selected $\set{v(1), \ldots, v(j-1)}$,  $v(j)$ is selected from $[n]\setminus \set{v(i): 1\leq i\leq j-1}$ with probability proportional to the corresponding weights $w_v,~ v\in [n]\setminus \set{v(i): 1\leq i\leq j-1} $.

Now we describe the construction.  We simultaneously construct the graph $\cG_n^{\nr}(\vw,\lambda)$ and explore it in a breadth-first manner.
Let $\set{\xi_{uv}: u,v\in [n],  u \neq v}$ be a collection of independent exponential random variables with rate
\begin{equation}
\label{eqn:ruv-def}
	r_{uv}:= \left(1+\frac{\lambda}{n^{1/3}}\right)\frac{w_v}{l_n}.
\end{equation}
We will use the above randomization to construct the graph simultaneously with an exploration process as follows. At each stage there will be a collection of active vertices $\cA(\cdot)$, a collection of explored vertices $\cO(\cdot)$ and a collection of unexplored vertices $\cU(\cdot)$.

\begin{enumeratea}
	\item {\bf Initialization:} Start by selecting the first vertex $v(1)\in [n]$ using \eqref{eqn:vone-choice} and let $\cA(0) = \set{v(1)}$. Further set $\cO(0) = \emptyset$ and $\cU(0) = [n]\setminus \set{v(1)}$.
	\item {\bf Recursion:}  For $i \geq 0$, given $\cA(i)$, $\cO(i)$ and $\cU(i)$, we construct $\cA(i+1), \cO(i+1)$ and $\cU(i+1)$ as follows. If $\cA(i) \neq \emptyset$, then we must have $\cA(i)  = \set{v(i), \cdots, v(i + |\cA(i)|-1)}$.
Arrange $\set{\xi_{v(i)u} : u\in \cU(i)}$ in increasing order $\xi_{v(i)v^\prime(1)} < \xi_{v(i)v^{\prime}(2)}< \cdots$.
Now define $\cN_i := \set{u \in \cU(i): \xi_{v(i) u} < w_{v(i)}}$ and let $c(i) := |\cN_i|$. List these vertices in $\cN_i$ as $v(i+|\cA(i)|):= v'(1),\hdots,v(i+|\cA(i)|+c(i)-1) := v'(c(i))$.
Then update $\cA(i+1):=\cA(i) \setminus \set{v(i)} \cup \cN_i$, $\cU(i+1) := \cU(i)\setminus \cN_i$ and $\cO(i+1) := \cO(i) \cup \set{v(i)}$.   If $\cA(i) = \emptyset$, then we select a new vertex $v \in \cU(i)$ with probability proportional to its weight $w_v$, and then define $\cN_i$ and $c(i)$ similarly. Now we list the vertices in $\cN_i$ as $v(i+1)=v'(1), \hdots, v(i+c(i)) = v'(c(i))$, and update $\cA(i+1) := \cN_i$, $\cU(i+1) := \cU(i)\setminus (\cN_i \cup \set{v(i)})$ and $\cO(i+1) := \cO(i) \cup \set{v(i)}$.

\end{enumeratea}

  This exploration process results in an ordering of the vertex set $[n]$ as $(v(1), v(2), \ldots, v(n))$. Consider the walk associated with the process
\begin{equation}
\label{eqn:brw-walk-def}
S_n(0) = 0, \qquad S_n(i) = S_n(i-1)+ c(i)-1. 	
\end{equation}
The construction satisfies
\begin{enumeratei}
	\item The ordering $(v(1), v(2), \ldots, v(n))$ has the same distribution as the size-biased re-ordering of the vertex set $[n]$ using the vertex weight sequence $\vw$.
	\item The walk $\set{S_n(i):i\geq 0}$ encodes the sizes of components (see \cite{aldous1997brownian}) in the following sense.  Write $T_{-k} = \min\set{i: S_n(i) = -k}$. The number of vertices in the first component explored by the walk (not necessarily the largest component) is given by $|\tilde{\cC}_1| = T_{-1}$, the size of the second component explored by the walk is given by $|\tilde{\cC}_1|:= T_{-2} - T_{-1}$ and so on and further for any $j\geq 1$ 	
	\begin{equation}
	\label{eqn:sn-walk-comp}
		S_n(T_{-j}) = -j, \qquad S_n(i) > -j \;\;\; \mbox{ for }\qquad T_{-(j-1)} < i < T_{-j}.
	\end{equation}
\end{enumeratei}
Thus excursions beyond past minima encode sizes of components in the order seen by the walk. By \cite{aldous1997brownian}, Theorem \ref{thm:comp-sizes-b-van} was proven in \cite{SBVHJL10} by showing that
\begin{equation}
\label{eqn:zn-to-wk}
\set{\frac{1}{n^{1/3}} S_n(sn^{2/3}): s\geq 0} \weakc \set{W_{\sqrt{\frac{\sigma_3}{\sigma_1}}, \frac{\sigma_3}{\sigma_1^2}}^\lambda(s): s\geq 0}.
\end{equation}
where $W^\lambda_{\kappa, \sigma}(\cdot)$ is the inhomogeneous Brownian motion as in \eqref{eqn:bm-lamb-kapp-def} and convergence is in the Skorohod metric $D(\bR_+, \bR)$. By the techniques in \cite{aldous1997brownian}, excursions beyond past minima of $S_n(\cdot)$ arranged in decreasing order converge to excursion beyond past minima of $W^\lambda_{\kappa, \sigma}$.

 Proposition \ref{prop:weight-control} describes asymptotics for the sum of the type $\sum_{v \in \cC_n^{\sss(i)}}u_v$, where $u_v := w_v$ or $u_v := w_v^2$. Note that by the above construction, the vertices within each component consist of a consecutive sequence in the above size-biased re-ordering using the weights $\vw$. The next lemma studies partial sum of this type with general weight sequences $\vw$ and $\vu$:



\begin{lem}
\label{lem:size-biased-sums}
Let $\vw = \vw^{\sss(n)} = \set{w_i^{\sss(n)} > 0: i \in [n]}$ be a set of weights, and $\vu = \vu^{\sss(n)} = \set{u_i^{\sss(n)}: i \in [n]}$ be a non-negative function on $[n]$. Let $m=m(n)\leq n$ be a increasing sequence of integers. Let $\set{v(i) \in [n] : i \in [n]}$ be a size-biased random reordering of the vertices based on the weight sequence $\vw$. Denote $w_{(i)} := w_{v(i)}$ and $u_{(i)} = u_{v(i)}$ for $i \in [n]$. Let $w_{\max}:= \max_{i \in [n]}w_i$ and $u_{\max} := \max_{i\in[n]}u_i$. Write $c_n := \sum_{i \in [n]} w_i u_i/ \sum_{i \in [n]}w_i$ for the weighted average of $\vu$. Assume that
\begin{align*}
	\lim_n\frac{m w_{\max}}{\sum_{i \in [n]} w_i} = 0 \;\; \mbox{ and } \;\; \lim_n\frac{u_{\max}}{m c_n} = 0.
\end{align*}
Define $Y(t) :=  (\sum_{i = 1 }^{\lfloor mt \rfloor}u_{(i)})/mc_n $, for $t \in [0,\infty)$, with  $u_{(i)} := 0$ for $i >n$. Then
\begin{equation*}
	 \sup_{ t \in [0,1]}  \left| Y(t) - t \right| \probc 0, \mbox{ as } n\to \infty.
\end{equation*}
\end{lem}
\begin{rem}
	The above lemma says that the average of the first $m$ values of $u_{(i)}$ is approximately $c_n$. The proof is a generalization of \cite[Lemma 2.3]{SBVHJL10}, which deals with the case when $u_i\equiv w_i^2$ and $m = n^{2/3}$.
\end{rem}
\noindent{\bf Proof of Lemma \ref{lem:size-biased-sums}:} This follows via the introduction of an extra randomization trick developed in \cite{aldous1997brownian} and also used in \cite{SBVHJL10}. We will give a full proof here. Define
\begin{equation*}
	\tau_k := \sum_{i \in [n]} w_i^k \mbox{ for } k =1,2.
\end{equation*}
For $i \in [n]$, let $\zeta_i \sim$ Exp($m w_i/\tau_1$) be independent exponential random variables. Define the process $\set{N(t): t \in [0,\infty)}$ as
\begin{equation*}
	N(t)  := \sum_{i \in [n]} \ind{\set{\zeta_i \leq t}},\;\; \mbox{ for } \;\;t \in [0,\infty).
\end{equation*}
Define the process $\set{ \tilde Y(t): t \in [0,\infty)}$ as
\begin{equation*}
	\tilde Y(t)  := \frac{1}{m c_n}\sum_{i \in [n]} u_{i} \ind{\set{\zeta_i \leq t}}, \;\;\mbox{ for } \;\; t \in [0,\infty).
\end{equation*}
Note that by the construction, we have $\set{Y(N(t)/m): t \geq 0 }\stackrel{d}{=}\set{\tilde Y(t): t \geq 0}$. Therefore when $\epsilon < 1$, on the event $\set{ |N(t)/m - t| < \epsilon, \forall t \in [0,2] }$ we have $N(2)/m > 1$ and thus
\begin{align*}
	  \sup_{ t \in [0,1]}  \left| Y(t) - t \right|
	\leq& \sup_{ t \in [0,2]}  \left| Y\left(\frac{N(t)}{m}\right) - \frac{N(t)}{m} \right|\\
	\leq& \sup_{ t \in [0,2]}  \left| Y\left(\frac{N(t)}{m}\right) - t \right| + \sup_{ t \in [0,2]}  \left| \frac{N(t)}{m} - t \right|\\
	\leq& \sup_{ t \in [0,2]}  \left| Y\left(\frac{N(t)}{m}\right) - t \right| + \epsilon.
\end{align*}
Thus
\begin{align}
	\pr\left( \sup_{ t \in [0,1]}  \left| Y(t) - t \right| > 2\epsilon \right) \leq \pr\left( \sup_{ t \in [0,2]}  \left| \tilde Y (t) - t \right| > \epsilon \right) + \pr\left( \sup_{ t \in [0,2]}  \left| \frac{N(t)}{m} - t \right| > \epsilon \right). \label{eqn:1437}
\end{align}
Then we bound the first term on the right hand side of \eqref{eqn:1437}. Define the filtration $\cF_t := \sigma(\set{\zeta_i \leq t} : i \in [n])$ for $t \geq 0$. Then for $t > s >0$,
\begin{align*}
	\E[\tilde Y(t) \mid \cF_s]
	=& \frac{1}{mc_n} \sum_{i \in [n]}\left[ u_i\ind{\set{\zeta_i \leq s}} + u_i \ind{\set{\zeta_i > s}} (1 - \exp({(t-s)mw_i}/{\tau_1} )\right]\\
	\leq& \tilde Y(s) + \frac{1}{mc_n} \sum_{i \in [n]} \frac{(t-s)mw_iu_i}{\tau_1}\\
	=& \tilde Y(s)  + (t-s).
\end{align*}
Therefore, by a supermartingale inequality \cite[Lemma 2.54.5]{rogers2000diffusions} for the supermartingale $\set{\tilde Y(t)- t : t \in [0,\infty)}$,
\begin{equation}
\label{eqn:super-mart-ineq}
	\pr(\sup_{ t \in [0,2]}|\tilde Y(t) - t|>\epsilon) \leq  \frac{9}{\epsilon} \left(|\E(\tilde Y(2) - 2)|+ \sqrt{\var(\tilde Y(2))}\right).
\end{equation}
Using the fact $x - x^2/2 \leq 1-e^{-x} \leq x $, it is easy to see that
\begin{align*}
	|\E[\tilde Y(2) - 2]| =& \left| \frac{1}{mc_n}\sum_{i \in [n]} u_i (1 - e^{-{2mw_i}/{\tau_1}})  - \frac{1}{mc_n} \sum_{i \in [n]} u_i \cdot \frac{ 2w_i m}{\tau_1} \right| \\
	 \leq& \frac{1}{m c_n} \sum_{i\in[n]} u_i \frac{4m^2 w_i^2 }{2\tau_1^2} = \frac{2m \sum_{i \in [n]}w_i^2 u_i}{\tau_1 \sum_{i\in[n]}w_iu_i}.
\end{align*}
For the variance
\begin{align*}
	\var(\tilde Y(2))
	=& \frac{1}{m^2 c_n^2} \sum_{i \in [n]} \left[u_i^2 (1 - \exp( - 2mw_i/\tau_1)) \exp( - 2mw_i/\tau_1)\right]\\
	\leq& \frac{1}{m^2 c_n^2} \sum_{i \in [n]} \frac{2mw_iu_i^2}{\tau_1} = \frac{2 \tau_1 (\sum_{i\in[n]}w_iu_i^2)}{m (\sum_{i\in[n]}w_iu_i)^2}.
\end{align*}
Similar bound holds for $\pr\left( \sup_{ t \in [0,2]}  \left| \frac{N(t)}{m} - t \right| > \epsilon \right)$ by plugging in $f(x)$ with the special choice of the function $f(x) \equiv 1$. Thus from \eqref{eqn:1437}
\begin{align*}
	&\pr\left( \sup_{ t \in [0,1]}  \left| Y(t) - t \right| > 2\epsilon \right)  \\
	\leq& \frac{9}{\epsilon}\left(\frac{2m \sum_{i \in [n]}w_i^2 u_i}{\tau_1 \sum_{i\in[n]}w_iu_i} + \sqrt{\frac{2 \tau_1 (\sum_{i\in[n]}w_iu_i^2)}{m (\sum_{i\in[n]}w_iu_i)^2}} + \frac{2m\tau_2}{\tau_1^2} + \sqrt{\frac{2}{m}} \right)\\
	\leq& \frac{9}{\epsilon}\left(\frac{2m w_{\max}}{\tau_1} + \sqrt{\frac{2 \tau_1 u_{\max}}{m \sum_{i\in[n]}w_iu_i}} + \frac{2m\tau_2}{\tau_1^2} + \sqrt{\frac{2}{m}} \right)
\end{align*}
The first two terms in the above display go to zero as $n\to\infty$, by assumptions in the lemma. Since $w_{\max} \geq \tau_2/\tau_1$ and $u_{\max} \geq c_n$, the remaining two terms also go to zero. This completes the proof of Lemma \ref{lem:size-biased-sums}. \qed\\

\noindent {\bf Completing the proof of Proposition \ref{prop:weight-control}:} Fix $i\geq 1$ and let $L(n,i)$ denote the time when we start exploring the $i$-th largest component in the above size-biased construction of $\cG_n(\lambda, \vw)$ and let $R(n,i)$ be the time when we complete the exploration of the $i$-th largest component so that $|\cC_n^{\sss(i)}|= R(n,i) - L(n,i)$. Let $L(\infty,i)$ denote the time of the start of the $i$-th largest excursion from zero of $\bar{W}_{\kappa, \sigma}^\lambda(\cdot)$ and $R(\infty,i)$ denote the end of this excursion where $\kappa, \sigma$ are as in Theorem \ref{thm:comp-sizes-b-van}. Thus the limiting component sizes are given by $Z_i = R(\infty, i) - L(\infty,i)$ and further by \cite{SBVHJL10,aldous1997brownian} $L(n,i)/n^{2/3} \weakc L(\infty,i)$ and $R(n,i)/n^{2/3}\weakc R(\infty,i)$.

Let $f : \bR_+ \to \bR_+$ be a monotone non-decreasing function. We will apply Lemma \ref{lem:size-biased-sums} to the case when $u_i \equiv f(w_i)$, $i \in [n]$, and $m = Tn^{2/3}$ for some large constant $T >0$. Define $c_n(f): ={\sum_{i \in [n]} w_i f(w_i)}/{\sum_{i \in [n]} w_i}$. Notice that
\begin{equation*}
	\frac{1}{n^{2/3}} \sum_{v\in \cC_n^{\sss(i)}} f(w_v) = \frac{1}{n^{2/3}} \sum_{j = 1}^{R(n,i)} f(w_{v(j)}) - \frac{1}{n^{2/3}} \sum_{j = 1}^{L(n,i)} f(w_{v(j)}). 	
\end{equation*}
Thus for any fixed $T> 0$ and $\eps> 0$
\begin{align*}
	&\pr\biggr( \frac{1}{n^{2/3}} \biggr| \sum_{v \in \cC_n^{\sss(i)} } f(w_v) - c_n(f)|\cC_n^{\sss(i)}| \biggr| > \eps\biggr)\\
	 \leq& \pr\left({R(n,i)} > Tn^{2/3}\right) + \pr\left(\sup_{u\leq T}\left|\frac{\sum_{i=1}^{n^{2/3} u} f(w_{v(i)})}{Tn^{2/3} c_n(f)} -  u\right| > \frac{\eps}{2T c_n(f)}\right).
\end{align*}
Take $f(x) = f_k(x) := x^k$ for $k=1,2$. By Assumption \ref{ass:wt-seq} (a) we have $\lim_{n\to \infty} c_n(f_k) = \sigma_{k+1}/\sigma_1$. The assumptions in Lemma \ref{lem:size-biased-sums} reduce to $w_{\max} = o(n^{1/3})$. Thus we can apply Lemma \ref{lem:size-biased-sums} to the second term in the above inequality. Letting $n \to  \infty$ first and then $T\to\infty$ in the above bound, for $k=1,2$,
\begin{equation}
	\label{eqn:1982}
	\left|\frac{\sum_{v\in \cC_n^{\sss(i)}} w_v^{k}}{n^{2/3}}  - c_n(f_k)\frac{|\cC_n^{\sss(i)}|}{n^{2/3}}\right|\probc 0, \mbox{ as } n \to \infty.
\end{equation}
Notice that for $k=1,2$, by Theorem \ref{thm:comp-sizes-b-van} we have
\begin{equation*}
	c_n(f_k)\frac{|\cC_n^{\sss(i)}|}{n^{2/3}} \weakc \frac{\sigma_{k+1}}{\sigma_1} Z_i \mbox{ as } n \to \infty.
\end{equation*}
Combining the above convergence, \eqref{eqn:1982} and the assumption $\sigma_2 = \sigma_1$, completes the proof of Proposition \ref{prop:weight-control}. \qed

\section{Completing the proof of Theorem \ref{thm:inhom-random-graph}}
\label{sec:finish-proof}
We shall now combine the various ingredients of the last two sections to complete the proof of Theorem \ref{thm:inhom-random-graph}. We start with the proof of convergence in the product topology.

\subsection{Convergence in the product topology}
 We work under Assumption \ref{ass:wt-seq} in this section. Due to the conditional independence of the components given the partition $(\cV^{\sss(i)}: i \geq 1)$, as described in Proposition \ref{prop:generate-nr-given-partition}, we can work with each maximal component separately. To ease notation let us work with the largest component  $\cC_n^{\sss(1)}(\lambda)$. By the Skorohod embedding, without loss of generality, we will work on a  probability space in which the convergence in Proposition \ref{prop:weight-control} holds almost surely namely,
\begin{equation}
	\label{eqn:2016}
	\left( \frac{|\cC_n^{\sss(1)}|}{n^{2/3}}, \frac{\sum_{v \in \cC_n^{\sss(1)}} w_v}{n^{2/3}}, \frac{\sum_{v \in \cC_n^{\sss(1)}} w_v^2}{n^{2/3}} \right)\convas \left(Z_i, Z_i, \frac{\sigma_3 Z_i}{\sigma_1}\right).
\end{equation}

Recall the definition of $M_i(\lambda)$ in \eqref{eqn:def-limit-mi-lambda}. Thus we need to prove
\begin{equation}
	\label{eqn:1034}
	\scl\left(\frac{1}{n^{1/3}}, \frac{1}{n^{2/3}}\right) \cdot \cC_n^{\sss(1)}(\lambda)  \weakc \cG\left(\frac{2\sigma_1^{1/2}}{\sigma_3^{1/2}} \tilde \ve_{Z_1}^{\sigma_3^{1/2}/\sigma_1^{3/2}}, \frac{\sigma_3^{1/2}}{\sigma_1^{3/2}} \tilde \ve_{Z_1}^{\sigma_3^{1/2}/\sigma_1^{3/2}}, \cP_1\right).
\end{equation}
By Proposition \ref{prop:generate-nr-given-partition}, conditional on the vertices in $\cC_n^{\sss (1)}(\lambda)$, the random graph $\cC_n^{\sss(1)}(\lambda)$ has the same distribution as a connected rank-one random graph as in \eqref{eqn:pr-con-vp-a-cV-def} using
\[\vp =\left( \frac{w_v}{\sum_{u\in \cC_n^{\sss(1)}} w_u}: v\in \cC_n^{\sss(1)} \right), \qquad a= \left(1+\frac{\lambda}{n^{1/3}}\right) \frac{(\sum_{v\in \cC_n^{\sss(1)}}w_v)^2}{l_n}. \]

\noindent\textbf{Proof of Theorem \ref{thm:inhom-random-graph} (i): } Our aim is to use Theorem \ref{thm:conv-condition-on-connected}. Let us first verify Assumptions \ref{ass:aldous-AMP} and \ref{ass:additional-connected}.  Notice that the relevant quantities are
\begin{align*}
	\sigma(\vp) = \frac{ \sqrt{  \sum_{v \in \cC_n^{\sss(1)}} w_v^2 }}{\sum_{v \in \cC_n^{\sss(1)}} w_v}, \qquad p_{\max} \leq \frac{w_{\max}}{\sum_{v \in \cC_n^{\sss(1)}} w_v},\qquad \mbox{ and } p_{\min} \geq \frac{w_{\min}}{\sum_{v \in \cC_n^{\sss(1)}} w_v}.
\end{align*}
By \eqref{eqn:2016} we have $\sigma(\vp)  = \Theta(n^{-1/3})$. Therefore, Assumption \ref{ass:aldous-AMP} can be verified with any $\epsilon \in (0, 3 \eta_0)$ and $r \in (2 + 3\gamma_0, \infty)$.
Assumption \ref{ass:additional-connected} is a consequence of \eqref{eqn:2016}:
\begin{equation*}
\lim_{n \to \infty} a \sigma(\vp)= \lim_{n \to \infty}\frac{(\sum_{v \in \cC_n^{\sss(1)}} w_v)^2}{l_n} \cdot \frac{ \sqrt{\sum_{v \in \cC_n^{\sss(1)}} w_v^2}} {\sum_{v \in \cC_n^{\sss(1)}} w_v} = \frac{\sigma_3^{1/2}}{\sigma_1^{3/2}} Z_1^{3/2} := \bar \gamma_1.
\end{equation*}
Thus Assumptions \ref{ass:aldous-AMP} and \ref{ass:additional-connected} are satisfied.
Now applying Theorem \ref{thm:conv-condition-on-connected}
\begin{equation}
	\label{eqn:1042}
	\scl\left(\sigma(\vp), \frac{1}{\sum_{v \in \cC_n^{\sss(1)}}w_v }\right) \cdot \cC_n^{\sss(1)}(\lambda) \weakc \cG( 2 \tilde \ve^{\bar \gamma_1}, \bar \gamma_1 \tilde \ve^{\bar \gamma_1}, \cP_1).
\end{equation}
By replacing $(l,\gamma,\theta)$ in the Brownian scaling \eqref{eqn:tilt-brown-exc-scale-general} with $(1, l^{3/2}, \gamma l^{-3/2})$, for all $\gamma >0$ and $l >0$,
\begin{equation}
	\label{eqn:1049}
	\set{ \tilde \ve^{\gamma}(s): s \in [0,1] } \stackrel{d}{=} \set{\frac{1}{l^{1/2}} \tilde \ve_l^{\gamma/l^{3/2}}(l s): s \in [0,1] }.
\end{equation}
By comparing the two scaling operators in \eqref{eqn:1034} and \eqref{eqn:1042},
\begin{equation*}
	\scl\left(\frac{1}{n^{1/3}}, \frac{1}{n^{2/3}}\right) = \scl\left( \frac{\sum_{v \in \cC_n^{\sss(1)}} w_v}{n^{1/3} \sqrt{\sum_{v \in \cC_n^{\sss(1)}} w_v^2}}, \frac{\sum_{v \in \cC_n^{\sss(1)}} w_v}{n^{2/3}}  \right) \cdot \scl\left(\sigma(\vp), \frac{1}{\sum_{v \in \cC_n^{\sss(1)}}w_v }\right).
\end{equation*}
Therefore by Proposition \ref{prop:scale-operator} and the convergence in \eqref{eqn:2016},
\begin{align*}
	\lim_{n\to\infty}\scl\left(\frac{1}{n^{1/3}}, \frac{1}{n^{2/3}}\right) \cdot \cC_n^{\sss(1)}(\lambda)
	\weakc &\scl\left( \sqrt{\frac{Z_1 \sigma_1}{\sigma_3}}, Z_1 \right) \cG( 2 \tilde \ve^{\bar \gamma_1}, \bar \gamma_1 \tilde \ve^{\bar \gamma_1}, \cP_1).
	\end{align*}
	Further note that the limit metric space on the right satisfies
	\begin{align*}
	\scl\left( \sqrt{\frac{Z_1 \sigma_1}{\sigma_3}}, Z_1 \right) \cG( 2 \tilde \ve^{\bar \gamma_1}, \bar \gamma_1 \tilde \ve^{\bar \gamma_1}, \cP_1) \stackrel{d}{=} &\scl\left( \sqrt{\frac{Z_1 \sigma_1}{\sigma_3}}, Z_1 \right) \cG\left( \frac{2}{Z_1^{1/2}} \tilde \ve_{Z_1}^{\sigma_3^{1/2}/\sigma_1^{3/2}}(Z_1 \cdot),  \frac{\bar \gamma_1}{Z_1^{1/2}}  \ve_{Z_1}^{\sigma_3^{1/2}/\sigma_1^{3/2}}(Z_1 \cdot), \cP_1\right)\\
	\stackrel{d}{=}& \cG\left( \frac{2\sigma_1^{1/2}}{\sigma_3^{1/2}} \tilde \ve_{Z_1}^{\sigma_3^{1/2}/\sigma_1^{3/2}},  \frac{\bar \gamma_1}{Z_1^{3/2}}  \ve_{Z_1}^{\sigma_3^{1/2}/\sigma_1^{3/2}}, \cP_1\right).
\end{align*}
Here the first line uses the scaling in \eqref{eqn:1049} with $l=Z_1$ and $\gamma = \bar \gamma_1$, and the second line uses the scaling in \eqref{eqn:887} and the scale invariance of $\cP_1$. Collecting the terms in the last display gives \eqref{eqn:1034}. The proof of Theorem \ref{thm:inhom-random-graph} (i) is completed. \qed

\subsection{Convergence in the $\cT_2$ topology}
\label{sec:conv-l4-topo-proof}

We will now strengthen the convergence in \eqref{eqn:prod-topology-convergence}
to convergence in the $\cT_2$ topology. Since $\lambda \in \bR$ is fixed, we will subsequently drop it from our notation.
We consider the Norros-Reittu model $\cG_n^{\nr}(\vw,\lambda)$ in this section. We first
require some notation. As usual, let $\cC_n^{\sss(i)}$ be the $i$-th largest component.
Denote the number of vertices in $\cC_n^{\sss(i)}$ by $|\cC_n^{\sss(i)}|$.
For $v\in[n]$, let $\cC_n(v)$ denote the component that contains $v$. For $k=1,2$ and $i\geq 1$, let
\begin{align}\label{eqn:def-X-n(v;k)}
X_n(v;k) :=\sum_{j\in\cC_n(v)}w_j^k\text{ and } X_{n, i}(k):=X_n(v;k)\text{ for any }v\in\cC_n^{\sss(i)}.
\end{align}
For $i \geq 1$ define
	\[\vp^{\sss(i)}=\left(w_j/X_{n,i}(1):\ j\in\cC_n^{\sss(i)}\right).\]
Note that $\vp^{\sss(i)}$ is the $\vp$ for $\cC_n^{\sss(i)}$, for $i \geq 1$. Let $\cF_{\partition}=\sigma\left(\set{w_v:\ v\in\cC_n^{\sss(i)}}_{i\geq 1}\right)$ be the $\sigma$-field
generated by the partition of weights into different components. Note that $X_{n,i}(k)$ is measurable with respect to $\cF_{\partition}$.

While proving convergence in the $\cT_2$ topology, the plan is to treat small components and large components differently. More precisely, we will use trivial bounds on the diameter and total mass of components $\cC_n^{\sss(i)}$ with $|\cC_n^{\sss(i)}| < n^{\alpha_0}$ for a suitably chosen $\alpha_0$, while for components with $|\cC_n^{\sss(i)}| \geq n^{\alpha_0}$, the following two lemmas will provide the necessary bounds:

\begin{lemma}\label{lem:comparability-of-functionals}
Let $X_{n}(v;k)$ be as above and set $\alpha_0=1/12-2\eta_1$ where $\eta_1$ is as in Assumption \ref{ass:high-moment}. Recall the definition of $\sigma_k$ from Section \ref{sec:res-rank-one}.
Then the following hold under Assumptions \ref{ass:wt-seq} and \ref{ass:high-moment}.
For any $r > 0$, there exists constants $n_0 > 0$ and $K_{\ref{lem:comparability-of-functionals}}=K_{\ref{lem:comparability-of-functionals}}(r)>0$ such that
\begin{enumeratea}
\item For all $v \in [n]$, $n^{\alpha_0} \leq m \leq n^{47/48}$, $k = 1,2$ and $n > n_0$,
\begin{align}\label{eqn:comparability-higher}
  \pr \left( X_n(v;k) \geq  \frac{32 \sigma_{k+1} m}{ \sigma_1} \text{ and } |\cC_n(v)| \leq m\right)\leq \frac{K_{\ref{lem:comparability-of-functionals}}}{n^{r}}.
\end{align}
\item For all $v \in [n]$, $n^{\alpha_0} \leq m \leq n^{45/48}$, $k = 1,2$ and $n > n_0$
\begin{align}\label{eqn:comparability-lower}
\pr \left(X_n(v;k) \leq \frac{\sigma_{k+1} m}{ 16\sigma_1} \text{ and }|\cC_n(v)|\geq m\right)\leq \frac{K_{\ref{lem:comparability-of-functionals}}}{n^{r}}.
\end{align}
\end{enumeratea}
\end{lemma}

Let $\underline{A} := 1/16 $ and $\overline{A} := 32\sigma_3/\sigma_1$. Define the event
\begin{equation}
	\label{eqn:def-en-2127-new}
	E_n(\alpha_0):=\set{\text{for } k = 1,2 \text{ and } v\in[n], |\cC_n(v)|\geq n^{\alpha_0} \text{ implies } \underline{A}|\cC_n(v)| \leq X_n(v; k)\leq \overline{A}|\cC_n(v)|}.
\end{equation}
\begin{lem}\label{lem:tail-bound-on-diameter-of-componenets}
Assume that Assumptions \ref{ass:wt-seq} and \ref{ass:high-moment} hold.
As in Lemma \ref{lem:comparability-of-functionals}, let $ \alpha_0 = 1/12-2\eta_1$. Then there exists constants $K_{\ref{lem:tail-bound-on-diameter-of-componenets}} > 0$ and $n_0 > 0$ such that for all $n \geq n_0$ and $\eta \in (0, 2\sigma_3/\sigma_1^{1/3})$,
\begin{align*}
\ind\set{E_n(\alpha_0)} \ind{\set{n^{\alpha_0} \leq |\cC_n^{\sss(i)}| \leq \eta n^{2/3} }} \E\left[(\diam(\cC_n^{\sss(i)}))^4 \mid \cF_{\partition} \right] \leq \frac{K_{\ref{lem:tail-bound-on-diameter-of-componenets}}}{[\sigma(\vp^{\sss(i)})]^4} \mbox{ for all } i \geq 1.
\end{align*}
\end{lem}
\noindent{\bf Proof of Lemma \ref{lem:comparability-of-functionals}:}
We break up the proof into two parts.

\noindent\textbf{Proof of \eqref{eqn:comparability-higher}:} For each $v\in[n]$, define the random permutation $\pi_v$ as follows:
$\pi_v(1)=v$ and $(\pi_v(2),\hdots,\pi_v(n))$ is a size-biased permutation of $[n]\backslash \{v\}$ where the size of $j$ is $w_j$.
Then $(\pi_v(i): i \geq 1)$ has the same law as the sequence of vertices of the random graph $\cG_n^{\nr}(\vw,\lambda)$ appear in a size-biased order during a breadth-first search starting from the vertex $v$. For ease of notation, we fix $v$ and write $\bar w_i := w_{\pi_v(i)}$ in the rest of the proof.

Hence,
\begin{equation*}
	Q_v := \pr\left( X_n(v;k) \geq  \frac{32 \sigma_{k+1} m}{ \sigma_1}  \text{ and } |\cC_n(v)| \leq m\right) \leq \pr \left( \sum_{i=1}^m \bar w_i^k \geq  \frac{32 \sigma_{k+1} m}{ \sigma_1}\right).
\end{equation*}
By Assumption \ref{ass:wt-seq} (a) there exists $n_1>0$ such that when $n \geq n_1$,
\begin{equation}\label{eqn:2718-912}
	\frac{\sigma_k}{2} <  \frac{\sum_{i=1}^n w_i^k}{n} < 2\sigma_k, \mbox{ for } k=1,2.
\end{equation}
We only give the proof when $k=2$, and the case when $k=1$ is similar. Let $\cF^v_j=\sigma\{\pi_v(i):\ 1\leq i \leq j\}$, for $1 \leq j \leq m$. Note that, for $2 \leq j \leq m$ and $n \geq n_1$, we have
\begin{align} \label{eqn:2624}
\E\left(\bar w_j^2\bigg|\cF^v_{j-1}\right)
=\frac{\sum_{i=j}^{n}\bar w_i^{3}}{\sum_{i=j}^{n} \bar w_i} \leq \frac{\sum_{i=1}^{n} w_i^{3}}{\sum_{i=1}^n w_i  - m w_{\max}}
\leq \frac{2 \sigma_{3} n }{\sigma_1 n /2 - m w_{\max}}.
\end{align}
By Assumption \ref{ass:high-moment}, there exists $n_2 > 0$ such that when $n \geq n_2$,
\begin{equation*}
	m w_{\max} \leq n^{47/48} w_{\max} < \frac{\sigma_1 n}{4}, \mbox{ and } w_{\max}^2 < \frac{8\sigma_{3} n^{1/12-2\eta_1}}{\sigma_1} \leq \frac{8\sigma_{3} m}{\sigma_1}.
\end{equation*}
By the above bound and \eqref{eqn:2624},
\begin{equation*}
	\bar w_1^2 + \sum_{j=2}^{m} \E [ \bar w_j^2 \mid \cF^v_{j-1} ] \leq w_{\max}^2 +  m \cdot \frac{8\sigma_{3}}{\sigma_1} < \frac{16\sigma_{3} m}{\sigma_1}.
\end{equation*}
Thus writing $\Delta_j :=\bar w_j^2- \E\left[\bar w_j^2\bigg|\cF^v_{j-1}\right]$, by the Burkholder-Davis-Gundy inequality (see e.g. \cite[Section 11.3]{chow2012probability}), we have for any integer $r' >0$,
\begin{align*}\label{eqn:sum-of-martingale-difference}
Q_v \leq  \pr\left(|\sum_{j=2}^{m}\Delta_j| \geq  \frac{16 \sigma_{3} m}{ \sigma_1}\right) \leq  C_1(r') \left(\frac{\sigma_1}{16\sigma_{3} m}\right)^{2r'}  \E\left[\left(\sum_{j=2}^m\Delta_j^2\right)^{r'}\right],
\end{align*}
where $C_1(r')$ is the constant in the Burkholder-Davis-Gundy inequality.
Notice that $|\Delta_j| \leq w_{\max}^2$. Further by Assumption \ref{ass:high-moment}, there exists $n_3 >0$ such that when $n > n_3$, $w_{\max} < n^{1/48-\eta_1}$. Therefore
\begin{equation}\label{eqn:2740-912}
	Q_v \leq C_1(r') \left(\frac{\sigma_1}{16\sigma_{3}}\right)^{2r'}  \frac{ w_{\max}^{4r'} }{m^{r'}}
	\leq  C_1(r')\left(\frac{\sigma_1}{16\sigma_{3}}\right)^{2r'}  \left(\frac{ n^{1/12 - 4 \eta_1} }{n^{1/12-2\eta_1}} \right)^{r'}
	= C_1(r') \left(\frac{\sigma_1}{16\sigma_{3}}\right)^{2r'}  \frac{1}{n^{2\eta_1 r'}}.
\end{equation}
We conclude the proof of \eqref{eqn:comparability-higher} by setting $n_0 = \max\set{n_1,n_2,n_3}$ and $r' = \lfloor r/2\eta_1 \rfloor + 1$.

We now turn to

\noindent\textbf{Proof of \eqref{eqn:comparability-lower}:} The idea is similar to the proof of
\eqref{eqn:comparability-higher}. Using the same notation,
\begin{equation*}
	\pr\left( X_n(v;k) \leq  \frac{\sigma_{k+1} m}{ 16\sigma_1}  \text{ and } |\cC_n(v)| \geq m\right) \leq \pr \left( \sum_{i=1}^m \bar w_i^k \leq  \frac{\sigma_{k+1} m}{ 16\sigma_1}\right).
\end{equation*}
Then when $n$ is large such that \eqref{eqn:2718-912} is true and
\begin{equation*}
	m w_{\max}^3 \leq n^{45/48}  w_{\max}^3 \leq \frac{\sigma_3 n}{4},
\end{equation*}
we have
\begin{equation*}
	\E\left(\bar w_j^2\bigg|\cF^v_{j-1}\right)
	\geq \frac{\sum_{i=1}^{n} w_i^{3} - m w_{\max}^3}{\sum_{i=1}^n w_i}
	\geq \frac{\sigma_{3}}{8 \sigma_1 }.
\end{equation*}
Then we can use a bound similar to \eqref{eqn:2740-912} to complete the proof.

In fact, the constant $K_{\ref{lem:comparability-of-functionals}}(r):=16^{2\lfloor r/2\eta_1 \rfloor + 2}C_1(\lfloor r/2\eta_1 \rfloor + 1) $ works for proving both \eqref{eqn:comparability-higher} and \eqref{eqn:comparability-lower} with $k=1,2$. The proof of Lemma \ref{lem:comparability-of-functionals} is completed. \qed

\medskip

\noindent{\bf Proof of Lemma \ref{lem:tail-bound-on-diameter-of-componenets}:}
In this proof, the constant $\gamma_0$ comes from Assumption \ref{ass:wt-seq}, and the constant $\eta_1$ comes from Assumption \ref{ass:high-moment}. For convenience, we will write $(\bar w_1,\hdots,\bar w_{m})$
for $(w_j: j\in\cC_n^{\sss(i)})$, where $m = m^{\sss(i)} := |\cC_n^{\sss(i)}|$. Let $p^{\sss(i)}_j=\bar w_j/X_{n,i}(1)$ for $1\leq j\leq m$
and let $a^{\sss(i)}$ be as in the statement of Proposition \ref{prop:generate-nr-given-partition}. Define $p_{\max}^{\sss(i)} := \max_{j \in [m]} p_j^{\sss(i)}$ and $p_{\min}^{\sss(i)} := \min_{j \in [m]} p_j^{\sss(i)}$. Further, let $L^{\sss(i)}(\vt)$ be as in \eqref{eqn:ltpi-def} with $a^{\sss(i)}, p_k^{\sss(i)}, p_{\ell}^{\sss(i)}$
replacing $a, p_k, p_{\ell}$ respectively.

Note that, $L(\vt)\geq 1$ for any $\vt\in\bT_m^{\ord}$. Define $\pr^{\sss(i)}(\cdot) := \pr_{\ord}(\cdot; \vp^{\sss(i)})$ where the latter is defined in \eqref{eqn:ordered-p-tree-def}. Thus, it follows from Proposition \ref{prop:generate-nr-given-partition} and Proposition \ref{prop:distrib-gone-tildg} that
\begin{align}
\E\left[(\diam(\cC_n^{\sss(i)}))^4 \mid \cF_{\partition} \right]
\leq \frac{\int {\height^4(\vt) L(\vt)} d\pr^{\sss(i)}(\vt)}{\int L(\vt) d\pr^{\sss(i)}(\vt)}
\leq \int \height^4(\vt) L(\vt)  d\pr^{\sss(i)}(\vt) \nonumber\\
\leq \left(\int {\height^8(\vt)} d\pr^{\sss(i)}(\vt)\right)^{1/2} \left(\int {L^2(\vt)} d\pr^{\sss(i)}(\vt)\right)^{1/2}.  \label{eqn:491}
\end{align}
Define $r_0 := 2\gamma_0/\alpha_0 + 2$ and $\epsilon_0 := 6 \eta_1$.  Define the events
\begin{equation}
	\label{eqn:def-h-2185}
	H_n^{\sss(i)} := \set{ \sigma(\vp^{\sss(i)}) \leq \frac{1}{2^{10}}, \frac{p_{\max}^{\sss(i)}}{[\sigma(\vp^{\sss(i)})]^{3/2+\epsilon_0}} \leq 1, \frac{[\sigma(\vp^{\sss(i)})]^{r_0}}{p_{\min}^{\sss(i)}} \leq 1, a^{\sss(i)} \sigma(\vp^{\sss(i)}) \leq \frac{1}{16}  }.
\end{equation}
Then restricted to $H_n^{\sss(i)}$, applying Theorem \ref{thm:ht-p-tree} with $r = 9$, we have
\begin{align}
	\label{eqn:2188}
	[\sigma(\vp^{\sss(i)})]^8 \int {\height^8(\vt)} d\pr^{\sss(i)}(\vt)
	\leq  \int_0^\infty 8 x^7 \pr^{\sss(i)} \left(\sigma(\vp^{\sss(i)})\height(\vt) \geq x \right) dx \leq 8 + \int_1^\infty 8x^7 \cdot \frac{K_{\ref{thm:ht-p-tree}}(9)}{x^9} dx,
\end{align}
Restricted to $H_n^{\sss(i)}$, applying Corollary \ref{cor:tightness-l-uniform-integrable} with $B_1 = 1/16$, $B_2 = 1$, and $\gamma = 2$,
\begin{equation}
	\label{eqn:2194}
	\int {L^2(\vt)} d\pr^{\sss(i)}(\vt) \leq K_{\ref{cor:tightness-l-uniform-integrable}}\left(2,\frac{1}{16},1\right).
\end{equation}
By \eqref{eqn:491}, \eqref{eqn:2188} and \eqref{eqn:2194}, the proof is completed once we show the following: there exist $n_0$ such that for all $n \geq n_0$,
\begin{equation}
	\label{eqn:claim-2201}
	E_n(\alpha_0) \cap \set{ n^{\alpha_0} \leq |\cC_n^{\sss(i)}| \leq \eta n^{2/3} } \subset H_n^{\sss(i)} \mbox{ for all } i \geq 1.
\end{equation}
On $E_n(\alpha_0)$, there exist absolute constants $C_1, C_2, C_3, C_4 > 0$ such that for all $i \geq 1$ and $n \geq 1$,
\begin{align*}
	&\frac{C_1}{\sqrt{m}} \leq \sigma(\vp^{\sss(i)}) = \frac{\sqrt{X_{n,i}(2)}}{X_{n,i}(1)} \leq \frac{C_2}{\sqrt{m}},   \\
	&{p_{\max}} \leq \frac{C_3 w_{\max}}{m}, \;\; p_{\min} \geq \frac{C_4 w_{\min}}{m}.
\end{align*}
The following calculation will be restricted to $E_n(\alpha_0) \cap \set{ n^{\alpha_0} \leq |\cC_n^{\sss(i)}| \leq \eta n^{2/3} }$. Note that
\begin{align}
	\frac{p_{\max}^{\sss(i)}}{[\sigma(\vp^{\sss(i)})]^{3/2+\epsilon_0}} \leq& \frac{C_3 w_{\max} /m} {(C_1/\sqrt{m})^{3/2 + \epsilon_0}} = \frac{C_2}{C_1^{3/2+\epsilon_0}} \cdot \frac{w_{\max}}{m^{1/4 - \epsilon_0/2}} \nonumber\\
	\leq& \frac{C_2}{C_1^{3/2+\epsilon_0}} \cdot \frac{n^{1/48 - \eta_1}}{n^{(1/12-2\eta_1)(1/4 - 3\eta_1)}} \leq \frac{C_2}{C_1^{3/2+\epsilon_0}} \cdot \frac{1}{n^{\eta_1/4}}, \label{eqn:2722}
\end{align}
Similarly
\begin{equation} \label{eqn:2725}
	\frac{[\sigma(\vp^{\sss(i)})]^{r_0}} {p_{\min}^{\sss(i)}} \leq \frac{(C_2/\sqrt{m})^{r_0}}{C_4 w_{\min}/m} = \frac{C_2^{r_0}}{C_4} \cdot \frac{1}{w_{\min}m^{r_0/2-1}} \leq \frac{C_2^{r_0}}{C_4} \cdot \frac{1}{w_{\min}n^{\gamma_0}}
\end{equation}
By \eqref{eqn:2722}, \eqref{eqn:2725} and Assumption \ref{ass:wt-seq} (d), there exists $n_1$ such that when $n \geq n_1$ the first three conditions in \eqref{eqn:def-h-2185} hold uniformly for all $i \geq 1$. Now we only need to verify the last condition in \eqref{eqn:def-h-2185}. Let $n_2$ be such that when $n \geq n_2$, $|\lambda|/n^{1/3} < 1$ and $l_n > n \sigma_1/2$, then when $n \geq n_2$,
\begin{equation*}
	a^{\sss(i)}\sigma(\vp^{\sss(i)}) = \left(1 + \frac{\lambda}{n^{1/3}}\right)\frac{(X_{n,i}(1))^2}{l_n} \cdot \frac{\sqrt{X_{n,i}(2)}}{X_{n,i}(1)} \leq \frac{4(\overline{A})^{3/2}}{\sigma_1} \cdot \frac{(m)^{3/2}}{n} \leq \frac{4(\overline{A})^{3/2}}{\sigma_1} \cdot \eta^{3/2} \leq \frac{1}{16},
\end{equation*}
where the last inequality uses $\overline{A} = 32\sigma_3/\sigma_1$ and $\eta < 2 \sigma_3/\sigma_1^{1/3}$. Therefore, when $n \geq n_0 := \max{\set{n_1, n_2}}$, the claim \eqref{eqn:claim-2201} is true. The proof of Lemma \ref{lem:tail-bound-on-diameter-of-componenets} is completed. \qed

\medskip

Before we start to prove Theorem \ref{thm:inhom-random-graph} (ii), we make a comment on the exponent $1/48 - \eta_1$ in Assumption \ref{ass:high-moment}, which is by no mean optimal. There are two key steps in our proof that contribute to this exponent: First, in \eqref{eqn:why-high-moment-1}, we use a very rough bound to control the diameters of all components of size less than $n^{\alpha_0}$, where the method only works when $\alpha_0 < 1/12$. Second, in the proof of Lemma \ref{lem:comparability-of-functionals}, in order to make \eqref{eqn:2740-912} work and obtain a tail bound on $X_n(v;2)$ for a component of size $|\cC_n(v)| \geq n^{\alpha_0}$, we introduce another factor of 4.

\noindent\textbf{Proof of Theorem \ref{thm:inhom-random-graph} (ii):} By choosing trivial $C$ and $\pi$ in the definition of $d_{\GHP}$ in \eqref{eqn:dghp}, for any $X_1$ and $X_2 \in \cS$,
\begin{equation*}
	d_{\GHP}(X_1, X_2) \leq \diam(X_1) + \diam(X_2) +  \mass(X_1) + \mass(X_2).
\end{equation*}
Therefore
\begin{align*}
&d_{\GHP}\left(\scl\left(\frac{1}{n^{1/3}}, \frac{1}{n^{2/3}}\right) \cC_n^{\sss(i)}, M_i(\lambda)\right) \\
&\hskip40pt\leq \frac{\diam\left(\cC_n^{\sss(i)}\right)}{n^{1/3}}
+\diam\left(M_i(\lambda)\right)
+ \frac{\mass\left(\cC_n^{\sss(i)}\right)}{n^{2/3}}
+\mass(M_i(\lambda)).
\end{align*}
Since we have already proved the convergence in the product topology, to prove convergence in the $\cT_2$ topology, it is enough to show that for any $\epsilon>0$
\begin{align}\label{eqn:501}
\limsup_{N\to\infty}\limsup_{n\to\infty} \left[
\pr\left(\sum_{i\geq N}\frac{\diam^4(\cC_n^{\sss(i)})}{n^{4/3}} > \epsilon \right)
+ \pr \left( \sum_{i\geq N}\frac{X_{n,i}(1)^4}{n^{8/3}}> \epsilon\right) \right]=0.
\end{align}

First we consider the first term in \eqref{eqn:501}. Notice that
\begin{equation}
	\label{eqn:why-high-moment-1}
	\frac{1}{n^{4/3}} \sum_{i \geq 1} \ind{\set{|\cC_n^{\sss(i)}| < n^{\alpha_0}}} \diam^4(\cC_n^{\sss(i)}) \leq \frac{n \cdot n^{4\alpha_0}}{n^{4/3}} = \frac{1}{n^{8\eta_1}}.
\end{equation}
So it is enough to focus on the components with size at least $n^{\alpha_0}$. Recall $E_n(\alpha_0)$ from \eqref{eqn:def-en-2127-new}.
We will first show
\begin{equation}\label{eqn:501A}
\pr\left(E_n(\alpha_0)^c\right)\to 0\text{ as }n\to\infty.
\end{equation}
By Lemma \ref{lem:comparability-of-functionals} (b), there exists $n_1$ such that when $n > n_1$, for $ k = 1, 2$,
\begin{align*}
&\pr\left(\exists v\in[n]\text{ with }|\cC_n(v)|\geq n^{\alpha_0} \text{ and }X_n(v; k)\geq \frac{32\sigma_{k+1}}{\sigma_1}|\cC_n(v)|\right)\\
\leq& \pr\left( |\cC_n^{\sss(1)}| > n^{3/4} \right) +  \sum_{v \in [n]} \sum_{m=n^{\alpha_0}}^{n^{3/4}}\pr\left(
|\cC_n(v)| = m \text{ and }
X_n(v; k) \geq \frac{32\sigma_{k+1}}{\sigma_1}|\cC_n(v)|\right)\\
\leq&  \pr\left( |\cC_n^{\sss(1)}| > n^{3/4} \right) + n \cdot n^{3/4} \cdot \frac{K_{\ref{lem:comparability-of-functionals}}}{n^2} =o(1),
\end{align*}
where the third line is a consequence of \eqref{eqn:comparability-higher} with $r =2$. By a similar argument
and an application of \eqref{eqn:comparability-lower}, we can show that
\[\pr\left(\exists v\in[n]\text{ with }|\cC_n(v)|\geq n^{\alpha_0}
\text{ and }X_n(v; k)\leq \frac{\sigma_{k+1}}{16\sigma_1}|\cC_n(v)|\right)=o(1).\]

By Cauchy-Schwarz inequality, $\sigma_3(n)\sigma_1(n)\geq\sigma_2(n)^2$. Letting $n\to\infty$ and using Assumption \ref{ass:wt-seq}(a) and (b), we get $\sigma_3 \geq \sigma_2 = \sigma_1$. Using this in the above equation yields \eqref{eqn:501A}.

Fix $\eta\in(0, 2\sigma_3/ \sigma_1^{1/3})$ (the upper bound of $\eta$ is due to Lemma \ref{lem:tail-bound-on-diameter-of-componenets}) and by Theorem \ref{thm:comp-sizes-b-van} we can find  $N_{\eta}$ such that, for all $n \geq 1$,
$$\pr\left(\sum_{i\geq N_{\eta}}n^{-4/3}|\cC_n^{\sss(i)}|^2>\eta\right)\leq\eta.$$
Define $$G_n(\alpha_0, \eta) := E_n(\alpha_0)\cap\set{\sum_{i\geq N_{\eta}}n^{-4/3}|\cC_n^{\sss(i)}|^2\leq\eta}.$$ Let $\sum_1$ denote the sum over all components $\cC_n^{\sss(i)}$ for which $i\geq N_{\eta}$ and $|\cC_n^{\sss(i)}|\geq n^{\alpha_0}$.
Then
\begin{align}\label{eqn:502}
&\pr\left(\sum\nolimits_1 n^{-4/3}\diam^4(\cC_n^{\sss(i)})>\eps\right)\\
\leq&\E\left[\ind\set{G_n(\alpha_0,\eta)}
\pr\left(\sum\nolimits_1 n^{-4/3}\diam^4(\cC_n^{\sss(i)})>\eps
\bigg|\cF_{\partition}\right)\right]+\pr(E_n(\alpha_0)^c)+\eta.\notag\\
\leq& \frac{1}{\eps n^{4/3}}\E\left[\ind\set{G_n(\alpha_0,\eta)}
\sum\nolimits_1 \E\left(\diam^4(\cC_n^{\sss(i)})\bigg|\cF_{\partition}\right)\right] + \pr(E_n(\alpha_0)^c)+\eta. \notag
\end{align}
By Lemma \ref{lem:tail-bound-on-diameter-of-componenets}, there exists $n_2$ such that for $n\geq n_2$,
\begin{align}\label{eqn:504}
& \frac{1}{\eps n^{4/3}}\E\left[\ind\set{G_n(\alpha_0,\eta)}
\sum\nolimits_1 \E\left(\diam^4(\cC_n^{\sss(i)})\bigg|\cF_{\partition}\right)\right] \\
\leq& \frac{1}{\eps n^{4/3}}\E\left[ \ind\set{G_n(\alpha_0,\eta)}\sum\nolimits_1 \frac{K_{\ref{lem:tail-bound-on-diameter-of-componenets}}}{[\sigma(\vp^{\sss(i)})]^4} \right]
\leq \frac{K_{\ref{lem:tail-bound-on-diameter-of-componenets}}}{\eps n^{4/3}}\E\left[  \ind\set{G_n(\alpha_0,\eta)} \sum\nolimits_1 |\cC_n^{\sss(i)}|^2 \right]
\leq \frac{\eta K_{\ref{lem:tail-bound-on-diameter-of-componenets}}}{\eps},\nonumber
\end{align}
where the last line uses the fact $[\sigma(\vp^{\sss(i)})]^2 |\cC_n^{\sss(i)}| \geq 1$ (see \eqref{eqn:m-sigma-p}) and the definition of $G_n(\alpha_0, \eta)$.

Combining \eqref{eqn:501A}, \eqref{eqn:502}, and \eqref{eqn:504}, we arrive at
\begin{align*}
\limsup_{n}\pr\left(n^{-4/3}\sum\nolimits_1 \diam^4(\cC_n^{\sss(i)})>\eps\right)
\leq \eta+ \frac{\eta K_{\ref{lem:tail-bound-on-diameter-of-componenets}}}{\eps}
\end{align*}
Since $\eta$ can be arbitrarily small, we conclude that
\begin{align}\label{eqn:507}
\limsup_{N\to\infty}\limsup_{n\to\infty}\pr\left(n^{-4/3}\sum_{i\geq N}\diam^4(\cC_n^{\sss(i)})>\eps\right)=0.
\end{align}

Next, we consider the second term in \eqref{eqn:501}. For components with $|\cC_n^{\sss(i)}| < n^{\alpha_0}$, on the event
\begin{equation}
	\label{eqn:2304-925}
	\set{ \mbox{for all } v \in [n], \; |\cC_n(v)| < n^{\alpha_0} \mbox{ implies } X_n(v,1) \leq 32 n^{\alpha_0} },
\end{equation}
we have
\begin{equation}\label{eqn:508}
	\frac{1}{n^{8/3}} \sum_{i \geq 1} \ind{\set{|\cC_n^{\sss(i)}| < n^{\alpha_0}}} X_{n,i}(1)^4 \leq 32^4 n \cdot \frac{n^{4\alpha_0}}{n^{8/3}} \leq \frac{32^4}{n^{4/3}}.
\end{equation}
By \eqref{eqn:comparability-higher}, the event in \eqref{eqn:2304-925} occurs with high probability and this take care of the small components. For components with $|\cC_n^{\sss(i)}| \geq n^{\alpha_0}$ we have
\begin{align}\label{eqn:509}
\ind\set{G_n(\alpha_0,\eta)} \sum\nolimits_1 \frac{X_{n,i}(1)^4}{n^{8/3}}
\leq  \ind\set{G_n(\alpha_0,\eta)} \sum\nolimits_1 \frac{\overline{A}^4 |\cC_n^{\sss(i)}|^4}{n^{8/3}}
\leq \overline{A}^4 \eta \frac{ |\cC_n^{\sss(1)}|^2}{n^{4/3}}.
\end{align}
Since $\eta$ is arbitrary and ${ |\cC_n^{\sss(1)}|}/{n^{2/3}}$ is tight, combining \eqref{eqn:508} and \eqref{eqn:509}, we conclude that
\begin{align*}
\limsup_{N\to\infty}\limsup_{n\to\infty}\pr\left(n^{-8/3}\sum_{i\geq N}X_{n,i}(1)^4>\eps\right)=0.
\end{align*}
This together with \eqref{eqn:507} yields \eqref{eqn:501} and completes the proof of Theorem \ref{thm:inhom-random-graph} (ii). \qed

\section{Tail bounds for height of $\vp$ trees:  Proof of Theorem \ref{thm:ht-p-tree}}
\label{sec:proof-height}
For the convenience of the reader, we restate the assumptions in Theorem \ref{thm:ht-p-tree} as follows.
\begin{ass}
	\label{ass:aldous-AMP-uniform-bound}
There exist $\epsilon_0  \in (0, 1/2)$ and $r_0 \in (2, \infty)$ such that
\begin{equation*}
	\sigma(\vp) \leq \frac{1}{2^{20}}, \qquad \frac{p_{\max}}{[\sigma(\vp)]^{3/2 + \epsilon_0}} \leq 1, \qquad \frac{[\sigma(\vp)]^{r_0}}{p_{\min}} \leq 1.
\end{equation*}
\end{ass}

We will prove the following lemma in this section.
\begin{lem}
	\label{lem:ht-p-tree-local-bound}
	Assume the setting of Theorem \ref{thm:ht-p-tree}. Then for any integer $r \geq  \lfloor r_0/2\epsilon_0 \rfloor + 1$, there exists a constant $K_{\ref{lem:ht-p-tree-local-bound}} = K_{\ref{lem:ht-p-tree-local-bound}}(r) > 0$ such that
\begin{equation*}
	\pr\left(\height(\cT) \geq \frac{x}{\sigma(\vp)}\right) \leq \frac{K_{\ref{lem:ht-p-tree-local-bound}}}{x^{r}}, \qquad \mbox{ for } 1 \leq x \leq [\sigma(\vp)]^{-2\epsilon_0}.
\end{equation*} 	
\end{lem}

Using Lemma \ref{lem:ht-p-tree-local-bound}, we prove Theorem \ref{thm:ht-p-tree} as follows:

\noindent\textbf{Proof of Theorem \ref{thm:ht-p-tree}:} Note that $\pr(\height(\cT) > m) = 0$. So it is enough to work with $[\sigma(\vp)]^{-2\epsilon_0} < x  \leq m \sigma(\vp)$. Take any $r \geq  \lfloor r_0/2\epsilon_0 \rfloor + 1$ and define $r' := \lfloor(r_0-1)r/(2\epsilon_0)\rfloor + 1 $. Then we have $r' \geq r \geq \lfloor r_0/2\epsilon_0 \rfloor + 1$, thus we can apply Lemma \ref{lem:ht-p-tree-local-bound} with $r'$. For $[\sigma(\vp)]^{-2\epsilon_0} < x  \leq m \sigma(\vp)$,
\begin{equation*}
	\pr\left(\height(\cT) \geq \frac{x}{\sigma(\vp)} \right) \leq \pr\left(\height(\cT) \geq \frac{[\sigma(\vp)]^{-2\epsilon_0}}{\sigma(\vp)} \right) \leq {K_{\ref{lem:ht-p-tree-local-bound}}(r') [\sigma(\vp)]^{2\epsilon_0r'}} \leq K_{\ref{lem:ht-p-tree-local-bound}}(r') \frac{ [\sigma(\vp)]^{r_0 r}}{ [\sigma(\vp)]^r}.
\end{equation*}
By Assumption \ref{ass:aldous-AMP-uniform-bound}, $[\sigma(\vp)]^{r_0} \leq p_{\min} \leq 1/m$. Then
\begin{equation*}
	\pr\left(\height(\cT) \geq \frac{x}{\sigma(\vp)} \right) \leq K_{\ref{lem:ht-p-tree-local-bound}}(r') \frac{ 1}{ [m \sigma(\vp)]^r} \leq  \frac{K_{\ref{lem:ht-p-tree-local-bound}}(r')}{x^r},
\end{equation*}
for all $ [\sigma(\vp)]^{-2\epsilon_0} < x \leq m \sigma(\vp)$. Defining $K_{\ref{thm:ht-p-tree}}(r) := K_{\ref{lem:ht-p-tree-local-bound}}( \lfloor(r_0-1)r/(2\epsilon_0)\rfloor + 1 )$, we complete the proof of Theorem \ref{thm:ht-p-tree}. \qed

The goal of the rest of this section is to prove Lemma \ref{lem:ht-p-tree-local-bound}. We will derive quantitative versions of some of the results of \cite{aldous2004exploration}.
We will also use the techniques developed in \cite{camarri2000limit}. Recall that $\cT$ is a rooted tree with vertex set labelled by $[m]$ and so given a vertex $v\in \cT$, we can let $\mathds{A}(v)$ be the set of ancestors of $v$. More precisely, writing $\height (v)$ for the height of vertex $v\in \cT$ and the path from the root $\rho$ to $v$ as $u_0 = \rho, u_1, \ldots, u_{\height (v)-1}, u_{\height (v)} = v$, then $\dA(v) = \set{u_0, \ldots, u_{\height (v)-1}}$.
Let $\cG(v) := \sum_{u\in \dA(v)} p_u$. Recall the function $F^{\exec,\vp}$ in \eqref{eqn:fexc-def} used to construct $\cT$. In particular recall that for each vertex $v\in \cT$, there is an $i$ such that we find the children of $v$ in the interval $[y^*(i-1), y^*(i))$. Define $e(v) = y^*(i)$.
Fix $x> 0$ and define the events
\begin{equation}
\label{eqn:b1-def}
	\cB_1:= \set{\max_{v\in [m]} \frac{F^{\exec,\vp}(e(v))}{\sigma(\vp)} \geq \frac{x}{8}},
\end{equation}
\begin{equation}
\label{eqn:b2-def}
	\cB_2:= \set{\max_{v\in [m]} \frac{F^{\exec,\vp}(e(v))}{\sigma(\vp)} \leq \frac{x}{8}, \max_{v\in [m]} \left( \frac{\cG(v)}{2\sigma(\vp)}-\frac{F^{\exec,\vp}(e(v))}{\sigma(\vp)} \right) \geq \frac{x}{8}},
\end{equation}
and finally
\begin{equation}
\label{eqn:b3-def}
	\cB_3:= \set{\max_{v\in [m]} \frac{\cG(v)}{2\sigma(\vp)} \leq \frac{x}{4}, \height(\cT) \geq \frac{x}{\sigma(\vp)}}.
\end{equation}
Thus
\begin{equation}
\label{eqn:ht-bd-b123}
	\pr\left(\height(\cT) \geq \frac{x}{\sigma(\vp)}\right) \leq \pr(\cB_1)+ \pr(\cB_2)+\pr(\cB_3).
\end{equation}
We will bound each one of the terms on the right individually in Lemmas \ref{lem:bound-b1}, \ref{lem:bound-b2} and \ref{lem:bound-b3},  whose proofs are given in Sections \ref{sec:proof-b1}, \ref{sec:proof-b2} and \ref{sec:proof-b3} respectively.

\noindent \textbf{Proof of Lemma \ref{lem:ht-p-tree-local-bound}:}
Combining \eqref{eqn:ht-bd-b123} and Lemmas \ref{lem:bound-b1}, \ref{lem:bound-b2} and \ref{lem:bound-b3}  completes the proof of Lemma \ref{lem:ht-p-tree-local-bound}. \qed

\subsection{Analysis of the event $\cB_1$}
\label{sec:proof-b1}
We will prove the following bound on $\pr(\cB_1)$:

\begin{lem}
	\label{lem:bound-b1}
	Under Assumption \ref{ass:aldous-AMP-uniform-bound},
	\begin{equation}\label{eqn:bound on B_1}
	\pr(\cB_1)\leq 12 e^{- {x^2}/{2^{16}}} \mbox{ for } 1 \leq x \leq  {128}{ [\sigma(\vp)]^{-2\epsilon_0}}.
	\end{equation}
\end{lem}
\noindent\textbf{Proof:} Replacing $x$ by $x\sigma(\vp)/8$ in Lemma \ref{lem:fexp-linfty-tail}, we have the same bound as in \eqref{eqn:bound on B_1}, but for all $x$ such that
\begin{equation*}
	\frac{32 p_{\max}}{\sigma(\vp)} \leq x \leq  \frac{128 \sigma(\vp)}{p_{\max}}.
\end{equation*}
Then by Assumption \ref{ass:aldous-AMP-uniform-bound}, we have ${32 p_{\max}}/{\sigma(\vp)} \leq 32 [\sigma(\vp)]^{1/2} \leq 1$ and ${128 \sigma(\vp)}/{p_{\max}} \geq {128} [\sigma(\vp)]^{-1/2 - \epsilon_0} \geq {128}{ [\sigma(\vp)]^{-2\epsilon_0}}$. This completes the proof of Lemma \ref{lem:bound-b1}. \qed

\subsection{Analysis of the event $\cB_3$}\label{sec:proof-b3}
We will prove the following bound on $\pr(\cB_3)$:
\begin{lem}
	\label{lem:bound-b3}
	Under Assumption \ref{ass:aldous-AMP-uniform-bound},  for each integer $r \geq \lfloor r_0/2\epsilon_0 \rfloor + 1$, there exists a constant $K_{\ref{lem:bound-b3}}=K_{\ref{lem:bound-b3}}(r)$ such that
	\begin{equation}\label{eqn:b3-final-bound}
		\pr(\cB_3) \leq \frac{K_{\ref{lem:bound-b3}}}{x^r} \mbox{ for }  x \geq 1.
	\end{equation}
\end{lem}

The proof of Lemma \ref{lem:bound-b3} uses the known connection between the $\vp$-tree and the first repeat time of an \emph{i.i.d.} sequence. Let $(\xi_i:i\geq 1)$ be \emph{i.i.d.} with distribution $\vp$ namely $\pr(\xi_i =j ) = p_j$ for $j\in [m]$ and $T$ is the first repeat time of this sequence namely
\begin{equation}
	\label{eqn:def-repeat-time}
	T := \min\set{j\geq 2: \xi_j = \xi_i \mbox{ for some } 1\leq i< j}.
\end{equation}
Define the random variables $X_j := p_{\xi_j}/(\sigma(\vp))^2 -1$ and set $S_j = \sum_{i=1}^j X_i $. Notice that the sequence $\set{S_j:j\geq 1}$ is a martingale. We will need the following concentration result about $S_j$.
\begin{lem}
	\label{lem:concentration-sj}
	For each integer $r \geq 1$, there exists a constant $K_{\ref{lem:concentration-sj}}=K_{\ref{lem:concentration-sj}}(r) >0$ such that for all $k \geq 1$ and $t >0$, we have
	\begin{equation*}
		 \pr\left(\max_{1\leq j \leq k} |S_j| \geq  t \right) \leq K_{\ref{lem:concentration-sj}} \cdot \frac{k^r p_{\max}^{2r}}{t^{2r} [\sigma(\vp)]^{4r}}.
	\end{equation*}
\end{lem}
 \noindent \textbf{Proof:} By Markov's inequality and the Burkholder-Davis-Gundy inequality \cite[Theorem 14.10]{dasgupta2011probability} (here it is clear by the definition of $X_j$ that the  $2r$-th moments are finite), we have for any integer $r > 0$,
\begin{align}
	\pr\left(\max_{1\leq j \leq k} |S_j| \geq  t \right) \leq \frac{\E\left[\max_{1 \leq j \leq k} |S_j |^{2r} \right] }{t ^{2r}} \leq C_1(r) \frac{ \E\left[(\sum_{j=1}^k X_j^2)^r\right] }{t^{2r}}, \label{eqn:2484}
\end{align}
where $C_1(r)$ is the constant that shows up in the Burkholder-Davis-Gundy inequality and only depends on $r$. Notice that $|X_1| \leq \max\set{\frac{p_{\xi_1}}{\sigma^2(\vp)}, 1 } \leq p_{\max}/\sigma^2(\vp)$. We have

\begin{align}
	\E\left[(\sum_{j=1}^k X_j^2)^r\right]
	\leq& k^r \E[|X_1|^{2r}] \leq k^r \frac{p_{\max}^{2r}}{[\sigma(\vp)]^{4r}}.  \label{eqn:2490}
\end{align}
Combining \eqref{eqn:2484} and \eqref{eqn:2490} proves the bound in Lemma \ref{lem:concentration-sj} with $K_{\ref{lem:concentration-sj}}(r) = C_1(r)$. \qed\\

\noindent\textbf{Proof of Lemma \ref{lem:bound-b3}:} Note that on the set $\cB_3$, there exists a vertex $v\in \cT$ such that
\begin{enumeratea}
	\item The height of this vertex satisfies $x/\sigma(\vp)\leq \height(v)\leq x/\sigma(\vp)+1$.
	\item For this $v$, since $\sigma(\vp)\height(v) \geq x$ and $\cG(v)/\sigma(\vp) \leq x/2$ (see the definition of $\cB_3$ in \eqref{eqn:b3-def}),
	\[\sigma(\vp)\height(v) - \frac{\cG(v)}{\sigma(\vp)}\geq \frac{x}{2}.\]
\end{enumeratea}
Thus
\begin{align}
\displaystyle	\pr(\cB_3) &\leq \frac{1}{p_{\min}}\sum_{v\in [m]} p_v \E\left(\ind\set{\sigma(\vp)\height(v) - \frac{\cG(v)}{\sigma(\vp)}\geq \frac{x}{2}, \height(v) \leq \frac{x}{\sigma(\vp)}+1}\right)\notag\\
	&=\frac{1}{p_{\min}} \pr\left(\sigma(\vp)\height(\dV) - \frac{\cG(\dV)}{\sigma(\vp)}\geq \frac{x}{2}, \height(\dV) \leq \frac{x}{\sigma(\vp)}+1\right) \notag\\
	&=: \frac{1}{p_{\min}} \pr(\cB_4),
\label{eqn:t3-t4-bound}	
\end{align}
where the first line uses the fact that $ p_v/p_{\min} \geq 1$ for all $v \in [m]$, and $\dV$ with distribution independent of $\cT$ is a vertex selected from $\cT$ with $\pr(\dV = j) = p_j$. Recall the definition of $T$ and $(\xi_i : i \geq 1)$ around \eqref{eqn:def-repeat-time}. By \cite[Corollary 3]{camarri2000limit}, we have
\begin{equation}
\label{eqn:ht-g-V-dist}
	(\height(\dV),\cG(\dV)) \stackrel{d}{=} \left(T-2,\sum_{i=1}^{T-1} p_{\xi_i}\right).
\end{equation}
Hence
\begin{align}
\pr(\cB_4) &= \pr\left((T-2)\sigma(\vp) -\frac{\sum_{i=1}^{T-1} p_{\xi_i}}{\sigma(\vp)}\geq \frac{x}{2},~ \sigma(\vp)(T-2) \leq x+\sigma(\vp)\right)\notag\\
&\leq \pr\left((T-1) -\frac{\sum_{i=1}^{T-1} p_{\xi_i}}{(\sigma(\vp))^2}\geq \frac{x}{2\sigma(\vp)},~ (T-1) \leq \frac{x+2\sigma(\vp)}{\sigma(\vp)}\right).
\label{eqn:b4-bound}
\end{align}
Recall $X_j = p_{\xi_j}/(\sigma(\vp))^2 -1$ and $S_j = \sum_{i=1}^j X_i $. Then we have
\begin{equation}
	\label{eqn:2475}
	\pr(\cB_4) \leq \pr \left( \max_{1 \leq j \leq 2 + x/\sigma(\vp)} |S_j| \geq \frac{x}{2\sigma(\vp)} \right).
\end{equation}
Applying Lemma \ref{lem:concentration-sj} to \eqref{eqn:2475} with $t = x/2\sigma(\vp)$ and $k = 2x/\sigma(\vp) > 2 + x/\sigma(\vp)$, we have for $r \geq 1$
\begin{equation*}
	\pr(\cB_4) \leq  K_{\ref{lem:concentration-sj}}(r) \left(\frac{2\sigma(\vp)}{x}\right)^{2r} \cdot \left(\frac{2x}{\sigma(\vp)} \right)^r  \cdot \frac{ p_{\max}^{2r}}{[\sigma(\vp)]^{4r}} =  \frac{K_{\ref{lem:concentration-sj}}(r)2^{3r}}{x^r} \cdot \frac{p_{\max}^{2r}}{[\sigma(\vp)]^{3r}}.
\end{equation*}
By Assumption \ref{ass:aldous-AMP-uniform-bound}, \eqref{eqn:t3-t4-bound} and the above bound, we have, for $x \geq 1$ and  $r \geq \lfloor r_0/2\epsilon_0 \rfloor + 1$,
\begin{equation*}
\pr(\cB_3)\leq \frac{K_{\ref{lem:concentration-sj}}(r)2^{3r}}{x^r} \cdot \frac{1}{[\sigma(\vp)]^{r_0}} \cdot {[\sigma(\vp)]^{2r\epsilon_0}} \leq \frac{K_{\ref{lem:concentration-sj}}(r)2^{3r}}{x^r}.
\end{equation*}
The proof of Lemma \ref{lem:bound-b3} is completed with $K_{\ref{lem:bound-b3}}(r) := K_{\ref{lem:concentration-sj}}(r)2^{3r}$. \qed

\subsection{Analysis of the event $\cB_2$}
\label{sec:proof-b2}
Let us now analyze $\cB_2$. In this section we will prove:
\begin{lem}
	\label{lem:bound-b2}
	Under Assumption \ref{ass:aldous-AMP-uniform-bound}, for each integer $r \geq \lfloor r_0/2\epsilon_0\rfloor + 1$, there exists a constant $K_{\ref{lem:bound-b2}}= K_{\ref{lem:bound-b2}}(r)$ such that
	\begin{equation}\label{eqn:bound on B_2}
		\pr(\cB_2) \leq \frac{K_{\ref{lem:bound-b2}}}{x^r} \mbox{ for } 1 \leq x \leq 8 [\sigma(\vp)]^{-2\epsilon_0}.
	\end{equation}
\end{lem}

We need the following tail bound on $T$ as defined in \eqref{eqn:def-repeat-time}:
\begin{lem}\label{lem:tail-bound-repeat-time}
	For $0 < t < 1/p_{\max}$,
	\begin{equation*}
		\pr\left( T \geq t \right) \leq 2 \exp \left( -\frac{t^2 \sigma^2(\vp)}{24} \right).
	\end{equation*}
\end{lem}
\noindent\textbf{Proof:} We will need an alternative construction of the random variable $T$, (see \cite[Section 4]{camarri2000limit}) where we essentially construct the sequence $\set{\xi_j:j\geq 1}$ in continuous time. The advantage of this construction is reflected in \eqref{eqn:cb-8-first-bound} below.  Using $\vp = (p_1,\ldots, p_m)$ partition the unit interval $[0,1]$ as $\set{I_{j}: j\in [m]}$ where $I_j$ has length $p_j$. Consider a rate one Poisson process $\cN$ on $\bR_+\times [0,1]$. We can represent $\cN=\set{(S_0,U_0), (S_1, U_1), \ldots}$ where $S_0< S_1 < \cdots$ are points of a rate one Poisson process on $\bR_+$ and $U_j$ are \emph{i.i.d.} uniform random variables. Abusing notation, write $\cN(t)$ for the number of points in $(0,t]\times [0,1]$ and $\cN(t^{-})$ for the number of points in $(0,t)\times[0,1]$. Now write $\xi_j = \sum_{i=1}^m i \ind\set{U_j\in I_i} $. In this continuous time construction, as before let $T$ denote the first repeat time of the sequence $\set{\xi_j:j\geq 1}$ and write $\cS$ for the actual ``time'' namely $\cS= \inf\set{s: \cN(s)> T}$. Thus $\cN(\cS^-) = \cN((0,\cS^{-})\times [0,1]) = T$. Then we have
\begin{align}
\pr(T \geq t) \leq \pr\left(\cS \leq {t}/{2} ,T \geq t \right)+ \pr\left(\cS \geq {t}/{2}\right) \label{eqn:cb-7-cb-8}.
\end{align}
Let us analyze $\pr\left(\cS \leq {t}/{2} ,T \geq t \right)$. Note that this event implies that $\cN(t/2) \geq t$. Standard tail bounds for the Poisson distribution then give
\begin{equation}
\label{eqn:cb-7-exp}
	\pr(\cN(t/2) \geq t) \leq  \exp\left(- \frac{t}{2} (2 \log 2-1)\right) < e^{-t/6}.
\end{equation}
Next, we bound $\pr\left(\cS \geq {t}/{2}\right)$. By \cite[Equations (26) and (29)]{camarri2000limit}, for $0 < t < 1/2p_{\max}$ we have
\begin{equation}
\label{eqn:cb-8-first-bound}
\log \pr(\cS > t) \leq  -\frac{t^2}{2}\sigma^2(\vp) + \frac{t^3 p_{\max} \sigma^2(\vp)}{3(1- tp_{\max})} \leq  - \frac{t^2 \sigma^2(\vp)}{6}.
\end{equation}
Replacing $t$ by $t/2$ in the above expression, we have for all $0 < t < 1/p_{\max}$, $\pr(\cS > t/2) \leq \exp( -t^2 \sigma^2(\vp)/24)$.  We complete the proof of Lemma \ref{lem:tail-bound-repeat-time} by combining the last bound, \eqref{eqn:cb-7-exp} and the fact $t \sigma^2(\vp) \leq \sigma^2(\vp)/p_{\max} \leq 1$. \qed

In order to prove Lemma \ref{lem:bound-b2}, we start with the following proposition:
\begin{prop}
\label{prop:gv-bound}
Under Assumption \ref{ass:aldous-AMP-uniform-bound}, for each integer $r \geq r_0$, there exists a constant $K_{\ref{prop:gv-bound}} = K_{\ref{prop:gv-bound}}(r)$ such that
\[\pr\left(\max_{v\in [m]} \cG^2(v) \geq x \sigma(\vp)\right) \leq \frac{K_{\ref{prop:gv-bound}}}{x^r} \mbox{ for all } \frac{1}{8} \leq x \leq [\sigma(\vp)]^{-2\epsilon_0}.\]
\end{prop}
\noindent{\bf Proof:}
We have
\begin{equation}
\label{eqn:t5-t6-bd}
	\pr\left(\max_{v\in [m]} \cG^2(v) \geq x \sigma(\vp)\right)\leq \pr(\cB_5) + \pr(\cB_6),
\end{equation}
where
\[\cB_5:= \set{\max_{v\in [m]} \cG(v)\geq \sqrt{x\sigma(\vp)},~\height(\cT)\leq \frac{\sqrt{x}}{2[\sigma(\vp)]^{3/2}}},\]
and
\[\cB_6:= \set{\height(\cT) \geq \frac{\sqrt{x}}{2[\sigma(\vp)]^{3/2}} }.\]

Arguing as in \eqref{eqn:t3-t4-bound}, we see that,
\begin{align}
	\pr(\cB_5) &\leq \frac{1}{p_{\min}}\pr\left(\cG(\dV)\geq \sqrt{x\sigma(\vp)},~\height(\dV)\leq \frac{\sqrt{x}}{2[\sigma(\vp)]^{3/2}}\right) \label{eqn:t5-bound},
\end{align}
where as before $\dV$ is selected from $[m]$ independent of $\cT$ using the probability vector $\vp$. Using the distributional representation \eqref{eqn:ht-g-V-dist} we conclude that for $x \geq 1/8$
\begin{align*}
\pr(\cB_5) &= \frac{1}{p_{\min}}\pr\left(\sum_{i=1}^{T-1} \frac{p_{\xi_i}}{(\sigma(\vp))^2} \geq \frac{\sqrt{x}}{[\sigma(\vp)]^{3/2}},\ T-2 \leq \frac{\sqrt{x}}{2[\sigma(\vp)]^{3/2}}\right)\\
&\leq \frac{1}{p_{\min}}\pr\left(\sum_{i=1}^{T-1} \left(\frac{p_{\xi_i}}{(\sigma(\vp))^2}-1\right) \geq \frac{\sqrt{x}}{4[\sigma(\vp)]^{3/2}},\ T-1 \leq \frac{\sqrt{x}}{[\sigma(\vp)]^{3/2}}\right)\\
&\leq \frac{1}{p_{\min}}
\pr\left(\max_{1\leq k\leq\sqrt{x}/[\sigma(\vp)]^{3/2}}S_k\geq\frac{\sqrt{x}}{4[\sigma(\vp)]^{3/2}}\right),
\end{align*}
where the second line uses the facts that $x \geq 1/8$ and $\sigma(\vp)\leq 2^{-10}$; and $S_k$ in the third line is as defined after \eqref{eqn:def-repeat-time}.
Using Lemma \ref{lem:concentration-sj} with $k=\sqrt{x}/[\sigma(\vp)]^{3/2}$ and $t = {\sqrt{x}}/{(4[\sigma(\vp)]^{3/2})}$ in the last display, for $x \geq 1/8$ and $r' \geq 2r_0$,
\begin{align*}
\pr(\cB_5)
& \leq \frac{K_{\ref{lem:concentration-sj}}(r')}{p_{\min}} \left( \frac{4 [\sigma(\vp)]^{3/2}}{ \sqrt{x}} \right)^{2r'} \cdot  \left(\frac{\sqrt{x}}{[\sigma(\vp)]^{3/2}}\right)^{r'}  \cdot \frac{p_{\max}^{2r'}}{[\sigma(\vp)]^{4r'}}
=  \frac{K_{\ref{lem:concentration-sj}}(r') 2^{4r'}}{p_{\min}} \cdot \frac{1}{x^{r'/2}}  \cdot \frac{p_{\max}^{2r'}}{[\sigma(\vp)]^{5r'/2}}\\
&\leq \frac{K_{\ref{lem:concentration-sj}}(r') 2^{4r'}}{[\sigma(\vp)]^{r_0}} \cdot \frac{1}{x^{r'/2}}  \cdot \frac{[\sigma(\vp)]^{3r'}}{[\sigma(\vp)]^{5r'/2}}
\leq \frac{K_{\ref{lem:concentration-sj}}(r') 2^{4r'}}{x^{r'/2}},
\end{align*}
where the second line uses $[\sigma(\vp)]^{r_0}/p_{\min} \leq 1$ and $p_{\max} \leq [\sigma(\vp)]^{3/2}$ in Assumption \ref{ass:aldous-AMP-uniform-bound}.
Letting $r' = 2r$ in the above display, we have for $x \geq 1/8$, $r \geq r_0$,
\begin{equation}\label{eqn:b5-final-bound}
\pr(\cB_5)\leq \frac{K_{\ref{lem:concentration-sj}}(2r) 2^{4r}}{x^{r}}.
\end{equation}
To finish the proof for Proposition \ref{prop:gv-bound}, we need to bound $\pr(\cB_6)$. Arguing as before and using the distributional
representation in \eqref{eqn:ht-g-V-dist} we first get
\begin{equation}
\label{eqn:b-6-first-bound}
\pr(\cB_6)\leq \frac{1}{p_{\min}}\pr\left(\height(\dV) \geq \frac{\sqrt{x}}{2[\sigma(\vp)]^{3/2}}\right)
\leq \frac{1}{p_{\min}}\pr\left(T \geq \frac{\sqrt{x}}{2[\sigma(\vp)]^{3/2}}\right).
\end{equation}
Applying Lemma \ref{lem:tail-bound-repeat-time} with $t = \sqrt{x}/(2[\sigma(\vp)]^{3/2})$ to \eqref{eqn:b-6-first-bound}, when $1/8 \leq x \leq [\sigma(\vp)]^{-2\epsilon_0}$, we have $t p_{\max} \leq p_{\max}/(2[\sigma(\vp)]^{3/2 + \epsilon_0}) < 1$ and $x > x/2 + 1/16$ and therefore
\begin{equation*}
	\pr(\cB_6) \leq \frac{2}{p_{\min}} \exp\left( -\frac{x}{96 \sigma(\vp)} \right) \leq \frac{2}{[\sigma(\vp)]^{r_0}} \exp\left(-\frac{1}{2^{11} \sigma(\vp)}\right ) \exp\left(-\frac{x}{192 }\right ) \leq C e^{-x/192},
\end{equation*}
where  $C := \sup_{y \geq 0} 2 y^{r_0} e^{-y/2^{11}}$. This combined with \eqref{eqn:b5-final-bound} finishes the proof of Proposition \ref{prop:gv-bound}.\qed

Recall the depth-first-exploration of the $\vp$-tree associated with $F^{\exec,\vp}(\cdot)$ in Section \ref{sec:ptree-dfs-connected-irg}. Recall also that if $v = v(i)$ is the $i$-th explored vertex, then we define $e(v):=y^*(i)$. The last crucial ingredient in the proof is the proposition stated below.
\begin{prop}\label{prop:qv-weird-bound}
For every $r\geq 1$, there exists a constant $K_{\ref{prop:qv-weird-bound}}=K_{\ref{prop:qv-weird-bound}}(r)$ such that under Assumption \ref{ass:aldous-AMP-uniform-bound}, for $v\in[m]$ and $x\geq 1$,
\begin{align*}\label{eqn:qv-weird-bound}
\pr\left(\frac{\cG(v)}{2}-\frac{\cG(v)^2}{2}-F^{\exec, \vp}(e(v))\geq\frac{x\sigma(\vp)}{16},\ \cG(v)\leq x\sigma(\vp)\right)
\leq K_{\ref{prop:qv-weird-bound}}\bigg(\frac{\sigma(\vp)}{x}\bigg)^{r}.
\end{align*}

\end{prop}

The proof of this proposition requires the next three lemmas. First, we setup some notation. Let $B\subset [m]$.
Let $\vq$ be the probability distribution obtained by merging the elements of $B$ into a single point. More precisely $\vq = (q_1, \ldots, q_{m - |B|+1})$ where $q_1 = \sum_{v\in B} p_v:=p(B)$ and $\set{q_i:i\geq 2}$ corresponds to the set $\set{p_i: i\in [m]\setminus B}$. Let $\cT(1,\vq)$ be a $\vq$-tree constructed as in \eqref{eqn:p-tree-def} with the probability mass function $\vq$, conditional on vertex $1$ being the root. Denote by $\cH_1$ the children of vertex $1$ in $\cT(1,\vq)$
and let
\begin{equation}\label{eqn:def-X}
X = \sum_{v\in \cH_1} q_v. 
\end{equation}
\begin{lem}[Bounds on $X$]
\label{lem:bound-on-sum-children}
Let $B\subseteq [m]$ and $q_1 = p(B)$ and define $K_{\ref{lem:bound-on-sum-children}} := \sup_{y \geq 0}ye^{-y/2^{12}} $. Then
\begin{align}\label{eqn:lower-bound}
\pr\bigg(q_1 - q_1^2 - X \geq \frac{x\sigma(\vp)}{32}\bigg)
\leq K_{\ref{lem:bound-on-sum-children}} \exp\left(-\frac{x^2}{2^{12}q_1}\right) \mbox{ for } x \geq 1.
\end{align}
Further, under Assumption \ref{ass:aldous-AMP-uniform-bound}, for every positive integer $r$, there exists $K'_{\ref{lem:bound-on-sum-children}}=K'_{\ref{lem:bound-on-sum-children}}(r)$ such that
\begin{align}\label{eqn:lower-bound-1}
\pr\bigg(X \geq q_1+x\sigma(\vp)\bigg)
\leq K'_{\ref{lem:bound-on-sum-children}}\frac{\sigma(\vp)^{r}}{x^{r}} \mbox{ for } x \geq 1,
\end{align}
whenever $q_1\leq 2x\sigma(\vp)$.
\end{lem}
\noindent {\bf Proof:} Let $\set{U_i: i\in[m]}$ be \emph{i.i.d.} Uniform$[0,1]$ random variables (independent of the tree and the coin tosses as well). Define the random variables
\[Y:=\sum_{i\notin B} p_i \ind\set{U_i\leq q_1}.\]
By the argument given below \cite[Equation 40]{aldous2004exploration},
\begin{equation}
\label{eqn:x-y-bound}
	\pr(X\in \cdot)\leq \frac{1}{q_1}\pr(Y\in \cdot).
\end{equation}
Consider the centered version
\[\tilde Y = \sum_{i\notin B} p_i\biggl(\ind\set{U_i\leq q_1} - q_1\biggr).\]
Then note that
\begin{equation}
\label{eqn:y-tildy}
	\pr\bigg(q_1 - q_1^2 - Y \geq x \sigma(\vp)/32\bigg) = \pr\left(-\tilde Y\geq x\sigma(\vp)/32\right).
\end{equation}
A Chernoff bound gives, for any $\lambda > 0$,
\begin{align*}
	\pr\left(-\frac{\tilde Y}{\sigma^2(\vp)}\geq \frac{x}{32\sigma(\vp)}\right)
\leq \exp\left(-\frac{\lambda x}{32\sigma(\vp)}\right)\prod_{i\notin B}&\exp\left(\frac{\lambda q_1 p_i}{\sigma^2(\vp)}\right)\\
	& \times \prod_{i\notin B}\left[1-q_1\left(1-\exp\left(-\frac{\lambda p_i}{\sigma^2(\vp)}\right)\right)\right].
\end{align*}
The simple inequality $1-u \leq \exp(-u)$ for $u\geq 0$ and some algebra gives
\begin{align*}
\pr\left(-\tilde Y \geq \frac{x \sigma(\vp)}{32}\right)
\leq \exp\left(-\frac{\lambda x}{32\sigma(\vp)}\right)
\prod_{i\notin B}\exp\left[q_1\left(\frac{\lambda p_i}{\sigma^2(\vp)} - 1 +\exp\left(-\frac{\lambda p_i}{\sigma^2(\vp)}\right)\right)\right].
\end{align*}
Since $e^{-u}-1+u \leq u^2/2$ for all $u \geq 0$, we finally get
\[\pr\left(-\tilde Y \geq \frac{x \sigma(\vp)}{32}\right)
\leq \exp\left(-\frac{\lambda x}{32\sigma(\vp)}\right) \exp\left(\frac{q_1 \lambda^2}{2\sigma^2(\vp)}\right).\]
Taking $\lambda= x\sigma(\vp)/(32 q_1)$, we get
\[\pr\left(-\tilde Y \geq {x \sigma(\vp)/32}\right) \leq \exp\big(-x^2/(2^{11}q_1)\big).\]
Using \eqref{eqn:x-y-bound}, \eqref{eqn:y-tildy} and $x^2 \geq x^2/2  + 1/2\ $ for $x\geq 1$, we arrive at
\begin{align*}
\pr\bigg(q_1 - q_1^2 - X\geq \frac{x\sigma(\vp)}{32}\bigg)
&\leq \frac{1}{q_1} \exp\left(-\frac{x^2}{2^{11}q_1}\right)
\leq \frac{1}{q_1} \exp \left(-\frac{1}{2^{12} q_1}\right)  \exp\left(-\frac{x^2}{2^{12}q_1}\right)\\
& \leq K_{\ref{lem:bound-on-sum-children}} \exp\left(-x^2/(2^{12}q_1)\right),
\end{align*}
which proves \eqref{eqn:lower-bound}. Next, note that
\begin{align}\label{eqn:89}
\pr\bigg(X\geq q_1+x\sigma(\vp)\bigg)
&\leq \pr\bigg(X-q_1(1-q_1)\geq x\sigma(\vp)\bigg)\\
&\leq \frac{1}{q_1} \pr\bigg(Y-q_1(1-q_1)\geq x\sigma(\vp)\bigg)=\frac{1}{q_1} \pr\bigg(\tilde{Y}\geq x\sigma(\vp)\bigg),\notag
\end{align}
where the second step follows from \eqref{eqn:x-y-bound}. Using Markov's inequality and Burkholder-Davis-Gundy inequality, for any positive integer $r$,
\begin{align}\label{eqn:90}
\pr\bigg(\tilde{Y}\geq x\sigma(\vp)\bigg)
&\leq \big(x\sigma(\vp)\big)^{-2r}\E(\tilde Y^{2r})\\
&\leq \frac{C_1(r)}{\big(x\sigma(\vp)\big)^{2r}}\E\left[\left(\sum_{j\in[m]}p_j^2\bigg(\ind\set{U_j\leq q_1}-q_1\bigg)^2\right)^r\right].\notag
\end{align}
A direct expansion yields
\begin{align}\label{eqn:91}
&\left(\sum_{j\in[m]}p_j^2\bigg(\ind\set{U_j\leq q_1}-q_1\bigg)^2\right)^r\\
&\hskip50pt=\sum_{\substack{s, r_1,\ldots,r_s\geq 1:\\r_1+\ldots+r_s=r}}\frac{r!}{r_1!\ldots r_s!}
\sum_{\substack{\text{ distinct}\\j_1,\ldots,j_s}}\
\prod_{\ell=1}^s \bigg(p_{j_{\ell}}\big(\ind\set{U_{j_{\ell}}\leq q_1}-q_1\big)\bigg)^{2r_{\ell}}.\notag
\end{align}
Since $\E\bigg(\ind\set{U_j\leq q_1}-q_1\bigg)^{r}\leq 2q_1$ for all $r\geq 1$, we get
\begin{align}\label{eqn:92}
&\E\bigg[\sum_{\substack{\text{ distinct}\\j_1,\ldots,j_s}}\
\prod_{\ell=1}^s \bigg(p_{j_{\ell}}\big(\ind\set{U_{j_{\ell}}\leq q_1}-q_1\big)\bigg)^{2r_{\ell}}\bigg]\\
&\hskip70pt\leq (2q_1)^s\sum_{\substack{\text{ distinct}\\j_1,\ldots,j_s}}\
\prod_{\ell=1}^s p_{j_{\ell}}^{2r_{\ell}}
\leq (2q_1)^s \prod_{\ell=1}^s\bigg(\sum_{j\in[m]} p_{j}^{2r_{\ell}}\bigg)
\leq (2q_1)^s p_{\max}^{2r-2s}\sigma(\vp)^{2s},\notag
\end{align}
where the last step uses the inequality: $\sum_{j\in[m]} p_{j}^{2r_{\ell}}\leq p_{\max}^{2r_{\ell}-2}\sigma(\vp)^2$. As a consequence of Assumption \ref{ass:aldous-AMP-uniform-bound} and $q_1\leq 2x\sigma(\vp)$,
\begin{align}\label{eqn:93}
(2q_1)^s p_{\max}^{2r-2s}\sigma(\vp)^{2s}
&\leq 2^r q_1\times\big(2x\sigma(\vp)\big)^{s-1}\times\sigma(\vp)^{3r-3s}\times \sigma(\vp)^{2s}\\
&\leq 2^{2r} q_1 x^{r-1}\sigma(\vp)^{3r-1}\notag.
\end{align}
We get \eqref{eqn:lower-bound-1} upon combining \eqref{eqn:89}, \eqref{eqn:90}, \eqref{eqn:91}, \eqref{eqn:92} and \eqref{eqn:93}.\qed

$X$ in Lemma \ref{lem:bound-on-sum-children} represents the sum of $p$-values associated with all children hanging off of the ancestral line of a vertex. We need a result that connects $X$ to the sum of $p$-values associated with those children that appear on the ``right side" of the ancestral line in the ordered tree. The next lemma serves this purpose.
\begin{lem}\label{lem:sum-children-to-right-children}
Let $y_{ij}$, $1\leq j\leq m_i$, $1\leq i\leq k$ be given positive numbers. Further, let $U_{ij}$, $1\leq j\leq m_i$, $1\leq i\leq k$ and $V_i$, $1\leq i\leq k$ be \emph{i.i.d.} Uniform$[0, 1]$ random variables. Define
\[Y:=\sum_{i=1}^k\sum_{j=1}^{m_i}y_{ij}\ind\set{U_{ij}\leq V_i}.\]
Then, for every positive integer $r$, there exists $K_{\ref{lem:sum-children-to-right-children}}=K_{\ref{lem:sum-children-to-right-children}}(r)$ such that
\begin{align*}
\pr\left(\bigg|Y-\sum_{i=1}^k\sum_{j=1}^{m_i}\frac{y_{ij}}{2}\bigg|\geq z\right)
\leq \frac{K_{\ref{lem:sum-children-to-right-children}}}{z^{2r}}
\max_{1\leq i\leq k}\left(\sum_{j=1}^{m_i}y_{ij}\right)^r\left(\sum_{i=1}^k\sum_{j=1}^{m_i}y_{ij}\right)^r
\end{align*}
for $z>0$.
\end{lem}
\noindent{\bf Proof:}
Note that
\begin{align}\label{eqn:97}
&\pr\left(\bigg|\sum_{i=1}^k\sum_{j=1}^{m_i}y_{ij}\ind\set{U_{ij}\leq V_i}-\sum_{i=1}^k\sum_{j=1}^{m_i}\frac{y_{ij}}{2}\bigg|\geq z\right)\\
&\hskip30pt\leq \pr\left(\bigg|\sum_{i=1}^k\sum_{j=1}^{m_i}y_{ij}\left(\ind\set{U_{ij}\leq V_i}-V_i\right)\bigg|\geq\frac{z}{2}\right)
+\pr\left(\bigg|\sum_{i=1}^k\sum_{j=1}^{m_i}y_{ij}\left(V_i-\frac{1}{2}\right)\bigg|\geq\frac{z}{2}\right)\notag\\
&\hskip30pt=:T_1+T_2.\notag
\end{align}

By Markov's inequality and Burkholder-Davis-Gundy inequality,
\begin{align*}
T_1
&\leq\frac{2^{2r}}{z^{2r}}\E\left(\sum_{i=1}^k\sum_{j=1}^{m_i}y_{ij}\left(\ind\set{U_{ij}\leq V_i}-V_i\right)\right)^{2r}\\
&\leq\frac{2^{2r}C_1(r)}{z^{2r}}\E\bigg[\sum_{i=1}^k\bigg(\sum_{j=1}^{m_i}y_{ij}\left(\ind\set{U_{ij}\leq V_i}-V_i\right)\bigg)^2\bigg]^{r}\notag\\
&\leq\frac{2^{2r}C_1(r)}{z^{2r}}\left(\sum_{i=1}^k\bigg(\sum_{j=1}^{m_i}y_{ij}\bigg)^2\right)^{r}
\leq \frac{2^{2r}C_1(r)}{z^{2r}}\max_{1\leq i\leq k}\left(\sum_{j=1}^{m_i}y_{ij}\right)^r\left(\sum_{i=1}^k\sum_{j=1}^{m_i}y_{ij}\right)^r.
\end{align*}
A similar bound on $T_2$ can be obtained in an identical fashion. Combined with \eqref{eqn:97}, this completes the proof.\qed

We will use the term $\sum_{j=1}^{m_i}y_{ij}$ in the above lemma to represent the sum over $p$-values of the children of a vertex on the ancestral line. The next lemma gives us control over this quantity.
\begin{lem}\label{lem:bound-max-children-weight}
For $i\neq j\in[m]$, let $\set{i\leadsto j}$ denote the event ``$i$ is the parent of $j$." Then
\begin{align*}
\pr\left(\sum_{j\in[m]\setminus\set{i}}p_j\ind\set{i\leadsto j}\geq p_{\max}^{3/4}\right)
\leq K_{\ref{lem:bound-max-children-weight}}\exp\left(-p_{\max}^{-1/4}\right),
\end{align*}
where $K_{\ref{lem:bound-max-children-weight}}=e^{e}$.
\end{lem}
\noindent{\bf Proof:}
By standard properties of $\vp$-trees (see, e.g., \cite[Section 6.1]{pitman2001random}), $\ind\set{i\leadsto j}$, $j\in[m]\setminus\set{i}$ are \emph{i.i.d.} Bernoulli$(p_i)$ random variables. Hence, by a Chernoff bound,
\begin{align}\label{eqn:98}
\pr\left(\sum_{j\in[m]\setminus\set{i}}p_j\ind\set{i\leadsto j}\geq p_{\max}^{3/4}\right)
&=\pr\left(\sum_{j\in[m]\setminus\set{i}}\frac{p_j}{p_{\max}}\ind\set{i\leadsto j}\geq p_{\max}^{-1/4}\right)\\
&\leq\exp\bigg(-p_{\max}^{-1/4}\bigg)\prod_{j\in[m]\setminus\set{i}}\bigg[1+p_i\bigg(\exp\left(\frac{p_j}{p_{\max}}\right)-1\bigg)\bigg].\notag
\end{align}
Since $p_j/p_{\max}\leq 1$,
\[\exp\left(\frac{p_j}{p_{\max}}\right)-1\leq \frac{e p_j}{p_{\max}},\]
and hence
\[1+p_i\bigg(\exp\left(\frac{p_j}{p_{\max}}\right)-1\bigg)
\leq 1+\frac{p_i e p_j}{p_{\max}}
\leq\exp\left(\frac{e p_i p_j}{p_{\max}}\right)
\leq\exp\bigg(ep_j\bigg).\]
We plug this into \eqref{eqn:98} to get the desired bound.
\qed

\noindent{\bf Proof of Proposition \ref{prop:qv-weird-bound}:}
Note that
\begin{align}\label{eqn:timmeeh1}
\cG(v)\leq \cG(v)+p_v=:\bar\cG(v)\text{ and }
\end{align}
\begin{align}\label{eqn:timmeeh2}
\bar\cG^2(v)\leq \cG^2(v)+3p_{\max}\leq \cG^2(v)+x\sigma(\vp)/16,
\end{align}
where the last inequality is a consequence of Assumption \ref{ass:aldous-AMP-uniform-bound} and the lower bound $x\geq 1$. Similarly, under Assumption \ref{ass:aldous-AMP-uniform-bound},
\begin{align}\label{eqn:timmeeh3}
\bar\cG(v)\leq\cG(v)+p_{\max}\leq 2x\sigma(\vp)
\end{align}
whenever $\cG(v)\leq x\sigma(\vp)$ for some $x\geq 1$. Next recall that $\dA(v)$ denotes the set of ancestors of $v$ and let $\dB(v)=\{v\}\cup\dA(v)$. Let $U_{ij}, V_i;\ i, j\in[m]$ be i.i.d. uniform$[0, 1]$ random variables. For $i, j\in[m]$, let $\set{i\leadsto j}$ denote the event ``$j$ is a child of $i$.'' Then note that
\begin{align}\label{eqn:timmeeh4}
F^{\exec, \vp}(e(v))\stackrel{d}{=}\sum_{i\in\dA(v)}\sum_{j\in[m]\setminus\dB(v)} p_j\ind\set{i\leadsto j,\ U_{ij}\leq V_i}+\sum_{j\in[m]}p_j\ind\set{v\leadsto j}.
\end{align}
Write $\sum_1$ for the sum $\sum_{i\in\dB(v)}\sum_{j\in[m]\setminus\dB(v)}$. Combining \eqref{eqn:timmeeh1} through \eqref{eqn:timmeeh4}, we get
\begin{align}\label{eqn:timmeeh5}
&\pr\left(\frac{\cG(v)}{2}-\frac{\cG(v)^2}{2}-F^{\exec, \vp}(e(v))\geq\frac{x\sigma(\vp)}{16},\ \cG(v)\leq x\sigma(\vp)\right)\\
&\hskip30pt\leq \pr\left(\frac{\bar\cG(v)}{2}-\frac{\bar\cG(v)^2}{2}-\sum\nolimits_1 p_j\ind\set{i\leadsto j,\ U_{ij}\leq V_i}\geq\frac{x\sigma(\vp)}{32},\ \bar\cG(v)\leq 2x\sigma(\vp)\right)\notag\\
&\hskip30pt\leq \pr\left(\bar\cG(v)-\bar\cG(v)^2-\sum\nolimits_1 p_j\ind\set{i\leadsto j}\geq\frac{x\sigma(\vp)}{32},\ \bar\cG(v)\leq 2x\sigma(\vp)\right)\notag\\
&\hskip30pt +\pr\left(\frac{1}{2}\sum\nolimits_1 p_j\ind\set{i\leadsto j}-\sum\nolimits_1 p_j\ind\set{i\leadsto j, U_{ij}\leq V_i}\geq\frac{x\sigma(\vp)}{64},\ \bar\cG(v)\leq 2x\sigma(\vp)\right)\notag\\
&\hskip30pt =: T_1+T_2.\notag
\end{align}

Suppose $A\subset [m]\setminus\set{v}$. Let $q_1:=p(A)+p_v$ and $\vq:=(q_1,q_2,\hdots, q_{m-|A|})$, where $(q_i: i\geq 2)$ corresponds to the set $\set{p_i: i\notin A\cup\set{v}}$. Now, the distribution of a $\vp$-tree conditional on $\dA(v)=A$ is same as the distribution of a $\vq$-tree conditional on vertex $1$ being the root. Since
\begin{align*}
T_1\leq \max_{\substack{A\subset [m]\setminus\set{v}:\\p(A)+p_v\leq 2x\sigma(\vp)}}
\pr\left(\bar\cG(v)-\bar\cG(v)^2-\sum\nolimits_1 p_j\ind\set{i\leadsto j}\geq\frac{x\sigma(\vp)}{32}\ \bigg|\ \dA(v)=A\right),
\end{align*}
an application of \eqref{eqn:lower-bound} yields
\begin{align}\label{eqn:66}
T_1\leq K_{\ref{lem:bound-on-sum-children}}\exp\left(-\frac{x}{2^{13}\sigma(\vp)}\right).
\end{align}

To bound the term $T_2$, we make note of the following two inequalities:
\begin{align}\label{eqn:67}
\pr\left(\max_{i\in\dB(v)}\sum_j p_j\ind\set{i\leadsto j}\geq p_{\max}^{3/4}\right)\leq mK_{\ref{lem:bound-max-children-weight}}\exp\left(-p_{\max}^{-1/4}\right),
\end{align}
which follows from Lemma \ref{lem:bound-max-children-weight}, and for every $x\geq 1$ and positive integer $r$,
\begin{align}\label{eqn:68}
&\pr\left(\sum\nolimits_1 p_j\ind\set{i\leadsto j}\geq \bar\cG(v)+x\sigma(\vp),\ \bar\cG(v)\leq 2x\sigma(\vp)\right)\\
&\leq \max_{\substack{A\subset [m]\setminus\set{v}:\\p(A)+p_v\leq 2x\sigma(\vp)}}
\pr\left(\sum\nolimits_1 p_j\ind\set{i\leadsto j}\geq \bar\cG(v)+x\sigma(\vp)\ \bigg|\ \dA(v)=A\right)
\leq K'_{\ref{lem:bound-on-sum-children}}\sigma(\vp)^{r}/x^{r},\notag
\end{align}
where the last inequality follows from Lemma \ref{lem:bound-on-sum-children}. From \eqref{eqn:67} and \eqref{eqn:68}, we conclude that
\begin{align}\label{eqn:69}
T_2
&\leq \pr\bigg(\frac{1}{2}\sum\nolimits_1 p_j\ind\set{i\leadsto j}-\sum\nolimits_1 p_j\ind\set{i\leadsto j, U_{ij}\leq V_i}\geq\frac{x\sigma(\vp)}{64},\ \bar\cG(v)\leq 2x\sigma(\vp),\\
&\hskip30pt \max_{i\in\dB(v)}\sum_j p_j\ind\set{i\leadsto j}\leq p_{\max}^{3/4},\
\sum\nolimits_1 p_j\ind\set{i\leadsto j}\leq \bar\cG(v)+x\sigma(\vp)\bigg)\notag\\
&\hskip30pt +mK_{\ref{lem:bound-max-children-weight}}\exp\left(-p_{\max}^{-1/4}\right)
+K'_{\ref{lem:bound-on-sum-children}}\sigma(\vp)^{r}/x^{r},\notag
\end{align}
for every $x\geq 1$ and positive integer $r$. Observe that the first term on the right side of \eqref{eqn:69} can be bounded by first conditioning on $\dA(v)$ and the children of $\dB(v)$, and then applying Lemma \ref{lem:sum-children-to-right-children}. Thus, for every $x\geq 1$ and positive integer $r$,
\begin{align}\label{eqn:70}
T_2
&\leq K_{\ref{lem:sum-children-to-right-children}}\left(\frac{64}{x\sigma(\vp)}\right)^{2r}\times(p_{\max})^{3r/4}\times(3x\sigma(\vp))^r\\
&\hskip30pt +mK_{\ref{lem:bound-max-children-weight}}\exp\left(-p_{\max}^{-1/4}\right)
+K'_{\ref{lem:bound-on-sum-children}}\sigma(\vp)^{r}/x^{r}.\notag
\end{align}
Under Assumption \ref{ass:aldous-AMP-uniform-bound}, $p_{\max}\leq\sigma(\vp)^{3/2}$. Further,
$\sigma(\vp)^{r_0}\leq p_{\min}\leq 1/m$. Note also that in Proposition \ref{prop:qv-weird-bound}, it is enough to consider $x\leq 16/\sigma(\vp)$. We combine these observations to get
\begin{align}\label{eqn:71}
m\exp\left(-p_{\max}^{-1/4}\right)&\leq m\exp(-\sigma(\vp)^{-3/8})\\
&\leq  m\exp\bigg(-\frac{m^{\frac{3}{8r_0}}}{2}\bigg)
\times \exp\bigg(-\frac{\sigma(\vp)^{-3/8}}{2}\bigg)\times \left(\frac{x}{\sigma(\vp)}\right)^r\times \left(\frac{\sigma(\vp)}{x}\right)^r\notag\\
&\leq\bigg[\sup_{m}\ m\exp\bigg(-\frac{m^{\frac{3}{8r_0}}}{2}\bigg)\bigg]
\times\bigg[\sup_{0<y<1}\ \exp\bigg(-\frac{y^{-3/8}}{2}\bigg)\times \left(\frac{16}{y^2}\right)^r\bigg]\times \left(\frac{\sigma(\vp)}{x}\right)^r,\notag
\end{align}
where the last step makes use of the bound $x\leq 16/\sigma(\vp)$. We now combine \eqref{eqn:timmeeh5} with \eqref{eqn:66}, \eqref{eqn:70} and \eqref{eqn:71} to complete the proof.
\qed

\noindent\textbf{Proof of Lemma \ref{lem:bound-b2}: }
By Proposition \ref{prop:qv-weird-bound}, for every $r\geq 1$, $x\geq 1$ and $v\in[m]$,
\[\pr\left(\frac{\cG(v)}{2}-\frac{\cG^2(v)}{2}-F^{\exec,\vp}(e(v))
\geq \frac{x\sigma(\vp)}{16}, \cG(v)\leq x\sigma(\vp)\right) \leq K_{\ref{prop:qv-weird-bound}} \left(\frac{\sigma(\vp)}{x}\right)^r.\]
Since $x \geq 1$ and $m \leq 1/p_{\min} \leq [\sigma(\vp)]^{-r_0}$, this yields, for $r\geq r_0$,
\begin{align*}
&\pr\left(\frac{\cG(v)}{2}-\frac{\cG^2(v)}{2}-F^{\exec,\vp}(e(v))
\geq \frac{x\sigma(\vp)}{16}\text{ and }\cG(v)\leq x\sigma(\vp)\text{ for some } v\in[m]\right)\\
\leq& K_{\ref{prop:qv-weird-bound}} m \left(\frac{\sigma(\vp)}{x}\right)^r
\leq \frac{K_{\ref{prop:qv-weird-bound}}}{x^r}.
\end{align*}
By the above bound and Proposition \ref{prop:gv-bound}, for $1 \leq x \leq 8 [\sigma(\vp)]^{-2\epsilon_0}$ and $r \geq r_0$,
\begin{align}\label{eqn:final1}
&\pr\left(\frac{\cG(v)}{2}-F^{\exec,\vp}(e(v))
\geq \frac{x\sigma(\vp)}{8}\text{ and }\cG(v)\leq x\sigma(\vp)\text{ for some }v\in[m]\right)\\
\leq& \frac{K_{\ref{prop:qv-weird-bound}}}{x^r} + \pr\left(\frac{1}{2}\max_{v\in[m]}\cG(v)^2\geq \frac{x\sigma(\vp)}{16}\right)
\leq \frac{K_{\ref{prop:qv-weird-bound}}}{x^r} + \frac{8^r K_{\ref{prop:gv-bound}}}{x^r}. \nonumber
\end{align}
Define
$$ E:=\left\{\frac{\|F^{\exec,\vp}\|_{L^{\infty}}}{\sigma(\vp)}\leq \frac{x}{8}\right\}\bigcap\left\{ \frac{\cG(v)}{2\sigma(\vp)}- \frac{F^{\exec,\vp}(e(v))}{\sigma(\vp)}
\geq \frac{x}{8}\text{ and }\frac{\cG(v)}{\sigma(\vp)}\leq x\text{ for some }v\in[m]\right\}.$$
Then $E$ and $\cB_2$ are the same provided $x\geq 2 p_{\max}/\sigma(\vp)$. Indeed, if
$\left\{\|F^{\exec,\vp}\|_{L^{\infty}}\leq x\sigma(\vp)/8\right\}$ holds and
$\{\cG(v_0)/2-F^{\exec,\vp}(e(v_0))\geq x\sigma(\vp)/8\}$ holds for some $v_0$ with
$\cG(v_0)\leq x\sigma(\vp)$, then $E$ is true. On the other hand, if $\cG(v_0)> x\sigma(\vp)$,
then there is an ancestor $v_1$ of $v_0$ satisfying $x\sigma(\vp)/2\leq \cG(v_1)\leq x\sigma(\vp)$
(this is true since $x\sigma(\vp)\geq 2 p_{\max}$). For this $v_1$, we have
$$\frac{\cG(v_1)}{2}-F^{\exec, \vp}(e(v_1))\geq \frac{x\sigma(\vp)}{4}-\frac{x\sigma(\vp)}{8}=\frac{x\sigma(\vp)}{8}.$$
Thus, the event $E$ is still true. Since $2p_{\max}/\sigma(\vp)\leq 1$ under Assumption \ref{ass:aldous-AMP-uniform-bound}, we conclude
from \eqref{eqn:final1} that for $r \geq r_0$, there exists some constant $K_{\ref{lem:bound-b2}}(r)$ depending only on $r$ such that for $1 \leq x \leq 8 [\sigma(\vp)]^{-2\epsilon_0}$,
\begin{equation*}
\pr(\cB_2) = \pr(E) \leq \frac{K_{\ref{lem:bound-b2}}}{x^r}.
\end{equation*}
This completes the proof of Lemma \ref{lem:bound-b2}.\qed\\

\appendix
\section{} \label{Appendix}
The aim of this section is to outline a proof of Lemma \ref{lem:AMP}.

\noindent{\bf Proof of Lemma \ref{lem:AMP}:} By imitating the proof of \eqref{eqn:lower-bound}, we can show that for any $\eps>0$,
\begin{align*}
\pr\bigg(q_1-q_1^2-X\geq \eps\sigma(\vp)/2\bigg)\leq \frac{4K_{\ref{lem:bound-on-sum-children}}}{\eps^2}\exp\left(-\frac{\eps^2}{2^{14}q_1}\right),
\end{align*}
where the meaning of the symbols is as in Lemma \ref{lem:bound-on-sum-children}. Now if $q_1\leq K\sigma(\vp)$ for some fixed $K>0$, then $q_1^2\leq\eps\sigma(\vp)/2$ for large $n$ since $\sigma(\vp)\to 0$ under Assumption \ref{ass:aldous-AMP}. Thus,
\begin{align}\label{eqn:122}
\pr\bigg(q_1-X\geq \eps\sigma(\vp)\bigg)\leq \frac{4K_{\ref{lem:bound-on-sum-children}}}{\eps^2}\exp\left(-\frac{\eps^2}{2^{14}K\sigma(\vp)}\right)
\end{align}
for large $n$ provided $q_1\leq K\sigma(\vp)$. Similarly, by imitating the proof of \eqref{eqn:lower-bound-1}, we can show that for any $\eps>0$ and integer $r\geq 2$,
\begin{align}\label{eqn:123}
\pr\bigg(X \geq q_1+\eps\sigma(\vp)\bigg)
\leq K'_{\ref{lem:bound-on-sum-children}}\left(\frac{K^r}{\eps^{2r}}\right)\sigma(\vp)^{r-1}
\end{align}
if $q_1\leq K\sigma(\vp)$ for some $K>0$. Now, with notation as in the proof of Proposition \ref{prop:qv-weird-bound},
\begin{align}\label{eqn:124}
&\pr\left(\bigg|\frac{\cG(v)}{2}-F^{\exec, \vp}(e(v))\bigg|\geq\eps\sigma(\vp),\ \cG(v)\leq K\sigma(\vp)\right)\\
&\hskip30pt\leq\pr\left(\bigg|\frac{\cG(v)}{2}-\sum\nolimits_1 p_j\ind\set{i\leadsto j,\ U_{ij}\leq V_i}\bigg|\geq\frac{\eps\sigma(\vp)}{2},\ \cG(v)\leq K\sigma(\vp) \right)\notag\\
&\hskip40pt +\pr\left(\bigg|F^{\exec, \vp}(e(v))-\sum\nolimits_1 p_j\ind\set{i\leadsto j,\ U_{ij}\leq V_i}\bigg|\geq\frac{\eps\sigma(\vp)}{2}\right)\notag\\
&\hskip40pt =: Z_1+Z_2.\notag
\end{align}
Under Assumption \ref{ass:aldous-AMP}, $p_{\max}^{3/4}\leq \sigma(\vp)^{9/8}\leq\eps\sigma(\vp)/2$ for large $n$. In view of \eqref{eqn:timmeeh4} and \eqref{eqn:67},
\begin{align}\label{eqn:125}
Z_2\leq mK_{\ref{lem:bound-max-children-weight}}\exp\left(-p_{\max}^{-1/4}\right).
\end{align}
As in the proof of Proposition \ref{prop:qv-weird-bound}, we use the following bound for $Z_1$:
\begin{align}\label{eqn:126}
Z_1
&\leq \pr\left(\big|\cG(v)-\sum\nolimits_1 p_j\ind\set{i\leadsto j}\big|\geq\frac{\eps\sigma(\vp)}{2},\ \cG(v)\leq K\sigma(\vp) \right)\\
&\hskip10pt+\pr\bigg(\bigg|\sum\nolimits_1 \frac{p_j}{2}\ind\set{i\leadsto j}-\sum\nolimits_1 p_j\ind\set{i\leadsto j,\ U_{ij}\leq V_i}\bigg|\geq\frac{\eps\sigma(\vp)}{4},\ \cG(v)\leq K\sigma(\vp),\notag\\
&\hskip30pt \sum\nolimits_1 p_j\ind\set{i\leadsto j}\leq\cG(v)+\frac{\eps\sigma(\vp)}{2},\ \max_{i\in\dB(v)}\sum_j p_j\ind\set{i\leadsto j}\leq p_{\max}^{3/4} \bigg)\notag\\
&\hskip10pt+\pr\bigg(\max_{i\in\dB(v)}\sum_j p_j\ind\set{i\leadsto j}\geq p_{\max}^{3/4}\bigg)
=:Z_{11}+Z_{12}+Z_{13}.\notag
\end{align}

From \eqref{eqn:122} and \eqref{eqn:123}, for every $r\geq 2$, there exists $C_r>0$ such that
\begin{align}\label{eqn:127}
Z_{11}\leq \max_{\substack{A\subset [m]\setminus\set{v}:\\ p(A)\leq K\sigma(\vp)}}
\pr\left(\big|\cG(v)-\sum\nolimits_1 p_j\ind\set{i\leadsto j}\big|\geq\frac{\eps\sigma(\vp)}{2}\ \bigg|\ \dA(v)=A\right)
\leq C_r\left(\frac{K^r}{\eps^{2r}}\right)\sigma(\vp)^{r-1}
\end{align}
for sufficiently large $n$. Similarly, conditioning on $\dA(v)$ and the children of $\dB(v)$ and using Lemma \ref{lem:sum-children-to-right-children},
\begin{align}\label{eqn:128}
Z_{12}&\leq K_{\ref{lem:sum-children-to-right-children}}\left(\frac{4}{\eps\sigma(\vp)}\right)^{2r}p_{\max}^{3r/4}
\left[\left(K+\frac{\eps}{2}\right)\sigma(\vp)\right]^r
\leq K_{\ref{lem:sum-children-to-right-children}}\left[\frac{16}{\eps^2}\left(K+\frac{\eps}{2}\right)\right]^{r}\sigma(\vp)^{r/8},
\end{align}
where the last step uses the bound $p_{\max}\leq\sigma(\vp)^{3/2}$ for large $n$. Finally, from \eqref{eqn:67},
\begin{align}\label{eqn:129}
Z_{13}\leq mK_{\ref{lem:bound-max-children-weight}}\exp\left(-p_{\max}^{-1/4}\right).
\end{align}
Combining \eqref{eqn:124}--\eqref{eqn:129} and using a union bound, we conclude that under Assumption \ref{ass:aldous-AMP},
\begin{align*}
\pr\left(\bigg|\frac{\cG(v)}{2}-F^{\exec, \vp}(e(v))\bigg|\geq\eps\sigma(\vp)\text{ for some }v\text{ with }\cG(v)\leq K\sigma(\vp)\right)=o(1).
\end{align*}
This shows that \cite[Equation (46)]{aldous2004exploration} remains valid under Assumption \ref{ass:aldous-AMP}. Now the claim follows from the arguments given after Equation (46) in \cite{aldous2004exploration}.
\qed

\section*{Acknowledgements}
We thank an anonymous referee for a thorough review of the paper. This significantly improved both the proofs as well as presentation of the paper.  We thank Amarjit Budhiraja for many stimulating conversations.
We also thank Gr\'egory Miermont for insightful discussions about the results of \cite{aldous2004exploration}.
SS thanks UNC Chapel Hill for hospitality and support during visits.
SB has been partially supported by NSF-DMS grants 1105581, 1310002, 160683, 161307 and SES grant 1357622.
SS has been supported in part by NSF grant DMS-1007524 and the Netherlands Organization for Scientific Research (NWO) through the Gravitation Networks grant 024.002.003.
XW has been supported in part by the National Science Foundation (DMS-1004418, DMS-1016441), the Army Research Office (W911NF-0-1-0080, W911NF-10-1-0158) and the US-Israel Binational Science Foundation (2008466).

\bibliographystyle{plain}
\bibliography{scaling}
\end{document}